\documentclass[10pt]{amsart}
\usepackage{amsmath}
\usepackage{xypic}
\usepackage{amssymb}
\usepackage{amsfonts}
\usepackage{amsthm}
\usepackage{enumerate}
\usepackage[all]{xy}
\usepackage[latin1]{inputenc}
\usepackage{graphicx}
\usepackage{latexsym}
\usepackage{accents}
\usepackage{mathdots}
\usepackage{mathrsfs}
\input xy

\makeatletter
    \newtheoremstyle{definition}
        {5pt}
        {3pt}
        {}
        {0pt}
        {\scshape}
        {.}
        {5pt}
        {\thmname{#1} \thmnumber{#2} \thmnote{[#3]}} 

\newtheoremstyle{theorems}
        {5pt}
        {3pt}
        {\itshape}
        {0pt}
        {\scshape}
        {.}
        {5pt}
        {\thmname{#1} \thmnumber{#2}\thmnote{[#3]}} 

\swapnumbers \theoremstyle{theorems}

\newtheorem{Theo}{Theorem}[section]
\newtheorem{Prop}[Theo]{Proposition}
\newtheorem{Cor}[Theo]{Corollary}
\newtheorem{Lemma}[Theo]{Lemma}
\newtheorem{Prop(BG)}[Theo]{Proposition (Bongartz-Gabriel)}
\newtheorem{Lemma(Asashiba)}[Theo]{Lemma(Asashiba)}
\newtheorem{Lemma(Gab)}[Theo]{Lemma(Gabriel)}
\newtheorem{Theo(Mil)}[Theo]{Theorem (Milicic)}

\theoremstyle{definition}
\newtheorem{Defn}[Theo]{Definition}

\newtheorem{Defn(Asashiba)}[Theo]{Definition (Asashiba)}

\newcommand{\Z}{\mathbb{Z}}
\newcommand{\T}{\mathbb{T}}

\newcommand{\cA}{\mathcal{A}}
\newcommand{\cB}{\mathcal{B}}
\newcommand{\cC}{\mathcal{C}}

\newcommand{\fF}{\mathfrak{F}}

\def\Sa{\hbox{${\mathit\Sigma}$}}
\def\Oa{\hbox{${\mathit\Omega}$}}

\def\La{\hbox{$\it\Lambda$}}

\newcommand{\Hom}{{\rm Hom}}
\newcommand{\Ext}{{\rm Ext}}
\newcommand{\End}{{\rm End}}

\def\Mod{\hbox{{\rm Mod}{\hskip 0.3pt}}}

\def\ModLa{\hbox{{\rm Mod}{\hskip 0.5pt}\La}}
\def\ModbLa{\hbox{{\rm Mod\hspace{0.6pt}}$^b$\hspace{-2pt}\La}}

\newcommand{\dt}{{\accentset{\hspace{2pt}\mbox{\large\bfseries .}}{}}}

\newcommand{\cdt}{\dt\hspace{2.5pt}}

\newcommand{\pdt}{{\hspace{1pt}\dt\hspace{1.5pt}}}
\newcommand{\ydt}{{\hspace{0.5pt}\dt\hspace{1.5pt}}}

\newcommand{\id}{{\rm 1\hspace*{-0.5ex}\rule{0.1ex}{1.52ex}\hspace*{0.2ex}}}

\begin{document}

\title[Derived categories]{\sc Koszul duality for non-graded derived categories}

\author[Ales Bouhada]{Ales Bouhada}

\author[Min Huang]{Min Huang}

\author[Shiping Liu]{Shiping Liu \vspace{-10pt}}



\thanks{The last named author is supported in part by the Natural Sciences and Engineering Research Council of Canada.}

\address{Ales Bouhada \\ D\'epartement de math\'ematiques, Universit\'e de Sherbrooke, Sherbrooke, Qu\'ebec, Canada.}
\email{mohammed.bouhada@usherbrooke.ca}

\address{Min Huang \\ D\'epartement de math\'ematiques, Universit\'e de Sherbrooke, Sherbrooke, Qu\'ebec, Canada.}
\email{min.huang@usherbrooke.ca}

\address{Shiping Liu\\ D\'epartement de math\'ematiques, Universit\'e de Sherbrooke, Sherbrooke, Qu\'ebec, Canada.}
\email{shiping.liu@usherbrooke.ca}

\maketitle

\begin{abstract}

We are concerned with relating derived categories of all modules of two dual Koszul algebras defined by a locally bounded quiver. We first gene\-ralize the well known Acyclic Assembly Lemma and formalize an old method of extending a functor from an additive category into a complex category to its complex category. Applying this to the Koszul functor associated with a Koszul algebra defined by a gradable quiver, we obtain a Koszul complex functor, that descends to an equivalence of a continuous family of pairs of triangulated subcategories of doubly unbounded complexes of the respective derived categories of all modules of the Koszul algebra and its Koszul dual. Under this special setting, this extends Beilinson, Ginzburg and Soegel's Koszul duality. In case the Koszul algebra is right or left locally bounded and its Koszul dual is left or right locally bounded respectively (for instance, the quiver has no right infinite path or no left infinite path), our Koszul duality restricts to an equivalence of the bounded derived categories of finitely supported modules, and an equivalence of the bounded derived categories of finite dimensional modules.

\vspace{-2pt}

\end{abstract}

\medskip

\section*{Introduction}

\medskip

The history of Koszul theory traces back to Cartan and Eilenberg's computing the cohomology groups of a Lie algebra using the Koszul resolution; see \cite[Chapter 8, Section 7]{CEi}. Later, various Koszul resolutions were used to compute the homology and the cohomology of Hopf algebras, restricted Lie algebras in characteristic 2, and Steenrod algebra; see \cite{BCKQRS, May}. In dealing with graded algebras arising from algebraic topology, Priddy formalized the Koszul theory of Koszul algebras and Koszul complexes and discovered a duality among homology algebras of certain Koszul algebras; see \cite{Pr1}. This beautiful theory has applications in many branches of mathematics such as algebraic topology; see \cite{GKM,Pr2}, algebraic geometry; see \cite{BGG, BGS}, quantum group; see \cite{Man}, commutative algebra; see \cite{E}, and representation theory of Lie algebras; see \cite{BGS,RH} and of associative algebras; see \cite{GMV,GMV2,Mar,Ma2}.

\medskip

Beilinson, Ginzburg and Soegel described the Koszul duality in terms of graded derived categories of two dual Koszul algebras; see \cite{BGS}, and for similar considera\-tion, see \cite{BGG,Fl,Kel,MvS}. More precisely, they established an equivalence between a triangulated subcategory of bounded-below complexes of the graded derived category of a left finite Koszul algebra and one of bounded-above complexes of the graded derived category of its Koszul dual, which induces an equivalence of the bounded derived categories of finitely generated graded modules in case the Koszul algebra is of finite length and its Koszul dual is noetherian. In \cite{MOS}, Mazorchuk, Ovsienko and Stroppel generalized the Koszul theory to positively graded categories; see also \cite{MSo}. In particular, they obtained a Koszul duali\-ty between a pair, the same as described in \cite{BGS}, of triangulated subcategories of the respective derived categories of graded modules of two dual Koszul categories; see \cite[Theorem 30]{MOS}.

\medskip

The Koszul duality in the existing literature always deals with derived categories of graded modules.
It is, however, also important to study all modules over Koszul algebras, for instance, those arising from mixed geo\-metry; see \cite[(1.4.2)]{BGS}. This is even more interesting from the representation theoretic point of view; see, for example, \cite{BaL,BaL2}. Indeed, given an algebra defined by a gradable quiver with radical squared zero, which turns out to be a Koszul algebra, Bautista and Liu established an equivalence between its bounded derived category of finitely supported modules and the bounded derived category of finitely co-presented modules of the path algebra of the opposite quiver, that is the Koszul dual; see \cite[(3.9)]{BaL2}.

\medskip

Our objective is to present a self-contained combinatorial account of Koszul theo\-ry for a graded algebra $\La$ defined by a locally finite quiver $Q$ and for all $\La$-modules; compare \cite{Mar}. The Koszul duality in this paper, however, is stated for the case where $Q$ is gradable. It extends, under this setting, the classic Koszul duality by establishing an equivalence between a continuous family of pairs of triangulated subcate\-gories, all but the classic one contain doubly unbounded complexes, of the non-graded derived categories of $\La$ and its Koszul dual $\La^!$. In contrast to the highly sophisticated technique of spectral sequences or adjunction isomorphisms used in \cite{BGS,MOS}, ours is elementary and consists of a generalization of the well known Acyclic Assembly Lemma; see (\ref{AAL}), a method of extending a functor from an additive category into a complex category to its complex category with a compatibility property; see (\ref{F-extension}) and (\ref{F-composition}), and an extension of a functorial quasi-isomorphism to a functorial quasi-isomorphism; see (\ref{main-lemma2}).

\medskip

Let $\La$ be quadratic, or equivalently, every simple module has a linear projective $2$-resolution; see (\ref{qua-alg}). We shall associate a local Koszul complex with each simple module $S$; see (\ref{K-cplx}) and compare \cite[(2.6)]{BGS}, that contains always a linear projective $2$-resolution of $S$ and is a projective resolution if and only if $S$ has a linear projective resolution; see (\ref{K-cplx}) and (\ref{Koz-cplx-dual}). More importantly, we shall construct two Koszul functors defined on the respective module categories of $\La$ and its quadratic dual $\La^!$; see (\ref{F-property}) and compare \cite[Page 489]{BGS}, which are extended to two complex functors between the complex categories of all modules of $\La$ and $\La^!$. Our genera\-lized Acyclic Assembly Lemma ensures that the complex functors descend to two derived Koszul functors between a continuous family of pairs of triangulated subcategories, all but one contain doubly unbounded complexes, of the respective derived categories of all modules of $\La$ and $\La^!$; see (\ref{F-diag}) and compare \cite[Proposition 20]{MOS}.

\medskip

Let $\La$ be a Koszul algebra, that is, every simple module has a linear projective resolution. In view of an alternative description of a local Koszul complex; see (\ref{k-cplx-iso}) and its dual; see (\ref{Koz-cplx-dual}), we shall see that the derived Kozul functors send the respective indecomposable injective modules and the indecomposable projective modules to simple modules; see (\ref{inj-im}). As a consequence, the derived Koszul functors are mu\-tually quasi-inverse on the modules. This yields in particular an explicit projective resolution of any module; see (\ref{natural-tran}) and an explicit injective co-resolution of any bounded-above module; see (\ref{natural-tran-2}). Finally, using some properties of the extension of functors and that of functorial morphisms, we conclude that the derived Koszul functors are mutually quasi-inverse; see (\ref{Main}). In case $\La$ is left (respectively, right) locally bounded and $\La^!$ is right (respectively, left) locally bounded (for instance, the quiver has no right infinite path or no left infinite path); see Section 1(4), our Koszul duality restricts to an equivalence of the bounded derived categories of finitely supported modules (respectively, of finite dimensional mo\-dules); see (\ref{Main-2}).
This generalizes the above-mentioned results in \cite[(2.12.6)]{BGS} and \cite[(3.9)]{BaL2}.

\smallskip

\section{Preliminaries}

\medskip

\noindent The objective of this section is to recall some background and collect some prelimi\-nary results. We shall also introduce some new classes of algebras, whose representation theory is worth further investigation in the future. The terminology and notation introduced in this section will be used for the rest of the paper.

\medskip

\noindent{\sc (1) Linear algebra.} Throughout, $k$ denotes a commutative field, and all tensor products will be over $k$. Let ${\rm Mod}\,k$ be the category of all $k$-vector spaces, and ${\rm mod}\,k$ the category of finite dimensional $k$-vector spaces. We shall use often the exact contravariant functor $D=\Hom_k(-, k): {\rm Mod}\, k\to {\rm Mod} k.$ The $k$-vector space freely spanned by a set $S$ will be written as $k \,S$. 
%
%
%
%
%
%
The following result is well known.

\medskip

\begin{Lemma}\label{Tensor}

If $U, V\in {\rm mod}\,k$ and $M, N\in {\rm Mod}\,k$, then we have a $k$-isomorphism
$$\varphi: \Hom_k(U, V) \otimes \Hom_k(M, N)\to \Hom_k(U\otimes M, V\otimes N): f\otimes g\mapsto \varphi(f\otimes g)$$ such that $\varphi(f\otimes g)(u\otimes m)=f(u) \otimes g(m)$ for all $u\in U$ and $m\in M.$

\end{Lemma}

\medskip

\noindent{\sc Remark.} In the sequel, for the sake of simplicity, we shall identify the map $\varphi(f\otimes g)$ stated in Lemma \ref{Tensor} with $f\otimes g$.

\medskip

As a consequence of Lemma \ref{Tensor}, we obtain the following well known result.

\medskip

\begin{Cor}\label{Cor 1.2}

Let $U, V$ be $k$-vector spaces.

\begin{enumerate}[$(1)$]

\item If $U\in {\rm mod}\hspace{.5pt}k$, then we obtain a natural $k$-isomorphism
$$\theta: \Hom_k(U, k) \otimes V \to \Hom_k(U, V): f\otimes v \mapsto \theta(f\otimes v)$$ such that $\theta(f\otimes v)(u)=f(u) v,$ for all $u\in U$ and $v\in V$.

\vspace{1.5pt}

\item If $U, V\in {\rm mod}\hspace{.5pt}k$, then we obtain a natural $k$-isomorphism
 $$\varphi: DU\otimes DV \to D(V\otimes U): f\otimes g\mapsto \varphi(f\otimes g),$$ where $\varphi(f\otimes g)(v\otimes u)= g(v)f(u)$, for all $u\in U$ and $v\in V$.

\end{enumerate}\end{Cor}

\medskip

The following statement will be needed later.

\medskip

\begin{Lemma}\label{dual-mix}

Let $f:U\to M$ and $g: N\to V$ be morphisms in ${\rm mod}\hspace{.5pt}k$. We obtain a commutative diagram with
vertical isomorphisms $$\xymatrix{
U\otimes DV \ar[d]_{\theta_{U, V}} \ar[rr]^-{f\otimes Dg} &&  M\otimes DN \ar[d]^{\theta_{M, N}} \\
D(V\otimes DU) \ar[rr]^-{D(g\otimes Df)} &&  D(N\otimes DM).}
$$

\end{Lemma}

\noindent{\it Proof.} Composing the isomorphism $U\otimes DV \to D^2U\otimes DV$, induced from the isomorphism $U\to D^2U$, with the isomorphism $D^2U\otimes DV\to D(V\otimes DU)$; see (\ref{Cor 1.2}), we obtain an isomorphism $\theta_{U, V}$ such that $\theta_{U, V}(u\otimes \zeta)(v\otimes \xi)=\zeta(v) \xi(u)$, for all $u\in U$, $v\in V$, $\zeta\in DV$ and $\xi\in DU$.
The proof of the lemma is completed.

\medskip

Let $U\in {\rm Mod}\,k$. Given a subspace $V$ of $U$, 
we shall write $V^\perp$ for the subspace of $DU$ of linear forms vanishing on $V$.
The following statement is well known.

\medskip

\begin{Lemma} \label{orthogonal}

Let $U\in {\rm mod}\,k$, and consider subspaces $V, W$ of $U$.

\begin{enumerate}[$(1)$]

\item We obtain $(V+W)^\perp =V^\perp \cap W^\perp$ and $(V\cap W)^\perp=V^\perp+W^\perp.$

\item If $\{u_1, \ldots, u_n\}$ and $\{v_1, \ldots, v_n\}$ are bases of $U$ with dual bases $\{u_1^\star, \ldots, u_n^\star\}$ and $\{v_1^\star, \ldots, v_n^\star\}$ respectively, then ${\sum}_{i=1}^nu_i\otimes u_i^\star={\sum}_{i=1}^n v_i\otimes v_i^\star$ in $U\otimes DV.$

\end{enumerate}

\end{Lemma}

\noindent{\it Proof.} Statement (1) is evident.
By Corollary \ref{Cor 1.2}, we obtain an isomorphism
$\theta: U\otimes \Hom_k(U, k)\to \End_k(U): u\otimes f\to \theta(u\otimes f).$ Given a basis $\{u_1, \ldots, u_n\}$ of $U$, it is easy to see that
$\theta\left({\sum}_{i=1}^nu_i\otimes u_i^\star\right)=\id_{\hspace{.4pt}U}.$ The proof of the lemma is completed.

%

\medskip

\noindent{\sc (2) Quivers.} Let $Q=(Q_0, Q_1)$ be a locally finite quiver, where $Q_0$ is a set of vertices and $Q_1$ is a set of arrows. Given an arrow $\alpha: x\to y$, we write $x=s(\alpha)$ and $y=e(\alpha)$. Given $x\in Q_0$, one has a {\it trivial path} $\varepsilon_x$ of length $0$ with $s(\varepsilon_x)=e(\varepsilon_x)=x$. A {\it path} of length $n>0$ is a sequence $\rho=\alpha_n \cdots \alpha_1$, with $\alpha_i\in Q_1$, such that $s(\alpha_{i+1})=e(\alpha_i)$, for $i=1, \ldots, n-1$; and in this case, we write $s(\rho)=s(\alpha_1)$ and $t(\rho)=t(\alpha_n)$, and call 
$\alpha_n$ the {\it terminal arrow} of $\rho$.

\medskip

Let $x, y\in Q_0$ and $n\ge 0$. Write $Q(x, y)$ for the set of paths from $x$ to $y$, and $Q_n$ for the set of paths of length $n$. Let $Q_{\le n}(x, y)$ be the subset of $Q(x, y)$ of paths of length $\le n$, and put $Q_n(x, y)=Q_n\cap Q(x, y)$. Moreover,
$Q_n(x, -)=\cup_{z\in Q_0} Q_n(x, z)$, $Q_n(-, x)=\cup_{z\in Q_0} Q_n(z, x)$ and $Q_{\le n}(x, -)=\cup_{z\in Q_0}\, Q_{\le n}(x, z)$.

\medskip

We say that $Q$ is {\it strongly locally finite} if $Q(x, y)$ is finite for all $x, y\in Q_0$; see \cite{BLP}, and {\it gradable} if $Q_0=\cup_{n\in \Z}\, Q^n$, a disjoint union called a {\it grading}, such that every arrow is of form $x\to y$, where $x\in Q^n$ and $y\in Q^{n+1}$ for some $n\in \Z;$ see \cite[(7.1)]{BaL}.

\medskip

The {\it opposite} quiver of $Q$ is a quiver $Q^{\rm o}$ defined in such a way that $(Q^{\rm o})_0=Q_0$ and $(Q^{\rm o})_1=\{\alpha^{\rm o}: y\to x \mid \alpha: x\to y \in Q_1\}$. A non-trivial path $\rho=\alpha_n\cdots \alpha_1$ in $Q(x, y)$, where $\alpha_i\in Q_1$, corresponds to a path $\rho^{\rm o}=\alpha_1^{\rm o} \cdots \alpha_n^{\rm o}$ in $Q^{\rm o}(y,x)$. However, the trivial in $Q$ at a vertex $x$ will be identified with the one in $Q^{\rm o}$ at $x$.

\medskip

\noindent{\sc (3) Path algebras.} Let $Q=(Q_0, Q_1)$ be a locally finite quiver. We shall write $kQ$ for the path algebra of $Q$ over $k$. Let $R$ be an ideal in $kQ$. Given $x, y\in Q_0$, we put $R(x, y)=R\cap kQ(x,y)$ and $R_n(x, y)=R\cap kQ_n(x,y)$ for all $n\ge 0$. An element $\rho\in R(x, y)$ is called a {\it relation} in $R$ from $x$ to $y$. Such a relation $\rho$ is called {\it quadratic} if $\rho\in kQ_2(x,y);$
{\it homogeneous} if $\rho\in kQ_n(x,y)$ for some $n\ge 2;$ {\it monomial} if $\rho\in Q(x, y)$; and {\it primitive} if $\rho=\sum_{i=1}^s \lambda_i\rho_i$, with $\lambda_i\in k$ and $\rho_i\in Q(x,y)$, such that $\sum_{i\in \it\Sigma}\lambda_i\rho_i\not\in R$ for any $\Sa\subset \{1, \ldots, s\}$. We shall say that $R$ is {\it quadratic} (respectively, {\it homogeneous}, {\it monomial}) if it is generated by a set of quadratic (respectively, homogeneous, monomial) relations. Finally, a {\it minimal generating set} of $R$ is a set $\Oa$ of primitive relations in $R$ such that $R$ is generated by $\Oa$ but not by any proper subset of $\Oa$, and in this case, we put $\Oa(x, y)=\Oa\cap kQ(x,y)$ and $\Oa(x,-)=\cup_{z\in Q_0}\, \Oa(x, z)$.

\medskip

\begin{Lemma}\label{min-rel}

Let $R$ be a homogenous ideal in $kQ$ with a minimal generating set $\Oa$. The classes of $\rho$ modulo $(kQ^+)R+R(kQ^+)$, with $\rho\in \Oa$, are $k$-linearly independent.

\end{Lemma}

\noindent{\it Proof.} Let $\rho_1, \ldots, \rho_r\in \Oa(x, y)$, with $x, y\in Q_0$, be such that $\lambda_1\rho_1+\cdots +\lambda_r\rho_r$ lies in $(kQ^+)R+R(kQ^+)$, where $\lambda_1, \ldots, \lambda_r\in k$ are non-zero. It is easy to see that $\rho_1={\sum}_{i=1}^s \gamma_i \rho_1\delta_i+ {\sum}_{j=1}^t \xi_j \hspace{.5pt} \sigma\hspace{-2pt}_j \zeta_j,$
where $\sigma_1, \ldots, \sigma_q\in\Oa\backslash \{\rho_1\}$, and $\gamma_i, \delta_i, \xi_j, \zeta_j\in kQ$ are homogeneous such that
$\gamma_i$ or $\delta_i$ is of positive degree for every $1\le i\le s$. Since $\rho_1$ and the $\sigma_j$ are homogeneous, $\rho_1={\sum}_{j\in \it\Theta} \,\xi_j \sigma_j \zeta_j,$ where $\it\Theta$ is the set of indices $j$ such that $\xi_j \hspace{.5pt} \sigma\hspace{-2pt}_j \zeta_j$ and  $\rho_1$ are of the same degree, a contradiction. The proof of the lemma is completed.

\medskip

In the sequel, the ideal $R$ is called {\it weakly admissible} if $R\subseteq (kQ^+)^2$, where $kQ^+$ is the ideal in $kQ$ generated by the arrows. A weakly admissible ideal $R$ is called {\it locally admissible} if, for any $(x, y)\in Q_0\times Q_0$, there exists $n_{xy}>0$ such that $kQ_n(x, y)\subseteq R$ for all $n\ge n_{xy}$; {\it right} (respectively, {\it left}) {\it admissible} if, for any $x\in Q_0$, there exists $n_{x}>0$ such that $kQ_n(x, -)\subseteq R$ (respectively  $kQ_n(-, x)\subseteq R$) for all $n\ge n_x$; and {\it admissible} if it is right and left admissible; compare \cite[(2.1)]{BoG}.

\medskip

The opposite algebra of $kQ$ is the path algebra $kQ^{\rm o}$. Given $\gamma=\sum_{i=1}^s\lambda_i \rho_i\in kQ$, where $\lambda_i\in k$ and $\rho_i$ are paths, we shall write
$\gamma^{\rm o}=\sum_{i=1}^s\lambda_i \rho_i^{\rm o}.$ This yields an algebra anti-isomorphism $kQ\to kQ^{\rm o}: \gamma \mapsto \gamma^{\rm o}$.

\medskip

\noindent{\sc (4) Algebras given by a quiver with relations.} Let $\La=kQ/R$, where $Q$ is a locally finite quiver and $R$ is a weakly admissible ideal in $kQ$. The pair $(Q, R)$ is called a {\it bound quiver}. Given $\gamma\in kQ$, we shall write $\bar{\gamma}=\gamma+R\in \La$. Writing $e_x= \bar{\varepsilon}_x$, we obtain a complete set $\{e_x \mid x\in Q_0\}$ of pairwise orthogonal primitive idempotents of $\La$. We shall say that $\La$ is {\it locally finite dimensional} if $e_y\La e_x$ is finite dimensional for all $x, y\in Q_0$; {\it strongly locally finite dimensional} if $R$ is locally admissible; {\it right} (respectively, {\it left}) {\it locally bounded} if $R$ is right (respectively,  left) admissible; and {\it locally bounded} if $R$ is admissible; compare \cite[(2.1)]{BoG}. It is evident that $\La$ is strongly locally finite dimensional if it is right or left locally bounded.

\medskip

Throughout, we shall write $J$ for the ideal in $\La$ generated by $\bar \alpha$ with $\alpha\in Q_1$.
We shall say that $J$ is {\it locally nilpotent} if, for any $(x, y)\in Q_0\times Q_0$, there exists $n_{xy}$ such that $e_y J^{n_{xy}} e_x=0$. We obtain the following easy statement.

\medskip

\begin{Prop} \label{Alg-dec}

Let $\La=kQ/R$, where $Q$ is a locally finite quiver and $R$ is a weakly admissible ideal in $kQ$.

\begin{enumerate}[$(1)$]

\item As a $k$-vector space, $\La=\La_0\oplus \La_1 \oplus J^2$, where $\La_0$ has a $k$-basis $\{e_x \mid x\in Q_0\}$ and $\La_1$ has a $k$-basis $\{\bar{\alpha} \mid \alpha\in Q_1\}.$

\vspace{1pt}

\item The ideal $R$ is locally admissible if and only if $J$ is locally nilpotent$\hspace{.8pt};$ and in this case, $J$ contains only nilpotent elements.

\vspace{1pt}

\item If $R$ is homogenous, then $\La$ is locally finite dimensional if and only if $\La$ is strongly locally finite dimensional.

\end{enumerate}

\end{Prop}

\noindent{\it Proof.}  Statement (1) holds because $J\subseteq (kQ^+)^2$. The first part of Statement (2) is evident. Given $u\in J$, we may write $u=\sum_{i=1}^s u_i$, where $u_i\in e_{y_i} J e_{x_i}$ with $x_i, y_i\in Q_0$. If $J$ is locally nilpotent, then there exists $n$ for which $e_{y_i} J^n e_{x_i}=0$ for all $1\le i, j\le s$, and consequently, $u^n=0$. This establishes Statement (2).

Assume that $R$ is homogeneous but not locally admissible such that $\La$ is locally finite dimensional. Then, $Q(x, y)$ for some $x, y\in Q_0$ has arbitrarily long paths not lying in $R$. Since $e_y\La e_x$ is finite dimensional, $\lambda_1 \delta_1+\cdots + \lambda_n \delta_n\in R(x, y)$, where $\lambda_1, \ldots, \lambda_n\in k$ are non-zero and $\delta_1, \ldots, \delta_n\in Q(x, y)\backslash R$ are of pairwise different lengths. Since $R$ is homogeneous, $\lambda_1 \delta_1+\cdots + \lambda_n \delta_n=\rho_1+\cdots+\rho_s,$ where $\rho_1, \ldots, \rho_s \in R(x, y)$ are homogeneous of pairwise different degrees. Then, each $\delta_i$ is a summand of a unique $\rho_j$, say $\rho_i$. This yields $\sum_{i=1}^n(\rho_i-\lambda_i\delta_i)+(\sum_{j>n}\rho_j)=0$, and $\lambda_i \delta_i=\rho_i$ for $i=1, \ldots, n$, a contradiction. The proof of the proposition is completed.

\medskip

\noindent{\sc Example.} (1) If $Q$ is a strongly locally finite quiver, then $kQ$ is strongly locally finite dimensional.

(2) If $\La=kQ/R$, where $Q$ consists of a loop $\alpha$ and $R$ is generated by $\alpha^2=\alpha^3$, then $\La$ is locally finite dimensional, but not strongly locally finite dimensional.

\medskip

The opposite algebra of $\La$ is $\La^{\rm o}=kQ^{\rm o}/R^{\rm o},$ where $R^{\rm o}= \{\rho^{\rm o} \mid \rho\in R\}$. Writing $\bar{\gamma\hspace{1.8pt}}^{\hspace{-1pt}\rm o}=\gamma^{\rm o}+R^{\rm o}$ for $\gamma\in kQ$, we obtain an anti-isomorphism $\La\to \La^{\rm o}: \bar{\gamma\hspace{1.8pt}} \to \bar{\gamma\hspace{1.8pt}}^{\hspace{-1pt}\rm o}$. For the sake of simplicity, $\varepsilon_a+R^{\rm o}$ will be written again as $e_a$ for any $a\in Q_0$.

\section{Double complexes and extension of functors}

\noindent Throughout this paper, all functors between linear categories are additive, and morphisms in a category are composed from the right to the left. An abelian category is called {\it concrete} if the objects admit a underlying abelian group structure and the composition of morphisms are compatible with these abelian group structures. 

\medskip

The objective of this section is to provide tools which will be used in Section 5. We shall start with some kind of homotopy theory of double complexes over a concrete abelian category. In particular, we shall generalize the well known Acyclic Assembly Lemma; see, for example, \cite{Wei}. This ensures the acyclicity of the total complex of a substantially large family of double complexes. Secondly, we shall formalize a method of extending a functor from an additive category into the complex category of another additive category to its complex category. Observe that this method has already been used in various special circumstances; see, for example, \cite{BaL2, BGS, Ke2, Ric}. We shall also study how to extend a functorial morphism between functors to a functorial morphism between the extended complex functors, and find sufficient conditions for a functorial quasi-morphism to be extended to a functorial quasi-morphism.

\medskip

Let $\cA$ be a concrete abelian category with countable direct sums. We denote by $C(\cA)$, $K(\cA)$ and $D(\cA)$ the {\it complex category}, the {\it homotopy category} and the {\it derived category} of $\cA$, respectively. We shall regard $\cA$ as a full subcategory of $C(\cA)$, by identifying an object $X$ with the stalk complex $X[0]$ concentrated in degree zero. Given a complex $(X^\cdt, d\hspace{-1.5pt}_X^{\,\pdt})$ in $C(\cA)$, we shall write $X^\ydt[1]$ for its {\it shift} and ${\rm H}^n(X^\ydt)$ for its {\it $n$-th cohomology group}. One says that $X^\ydt$ is {\it acyclic} if ${\rm H}^n(X^\ydt)=0$ for all $n\in \Z$. Moreover, the {\it twist} $\mathfrak{t}(X^\cdt)$ of $X^\cdt$ is a complex whose $n$-component is $X^n$ and whose $n$-th differential is $-d_X^n$; see \cite[Section 4]{BaL2}. Clearly, $\mathfrak{t}(X^\cdt)\cong X^\cdt$. Finally, the mapping cone of a morphism $f^\ydt: X^\ydt\to Y^\ydt$ in $C(\cA)$ will be written as $C_{f^\pdt}$.

\medskip

We shall now study double complexes. Let $(M^{\ydt\hspace{.5pt} \ydt},v_{\hspace{-1pt}_M}^{\ydt\hspace{.5pt} \ydt}, h_{\hspace{-1pt}_M}^{\ydt\hspace{.5pt} \ydt})$ \vspace{1pt} be a double complex over $\cA$. Given $i, j\in \Z$,
the {\it $i$-th column} of $M^{\ydt\hspace{.5pt} \ydt}$ is the complex $(M^{i, \ydt}, v_{\hspace{-1pt}_M}^{i,\,\ydt\,})$, and the {\it $j$-th row} is the complex $(M^{\cdt, \hspace{.5pt} j}, h_{\hspace{-1pt}_M}^{\cdt, \hspace{.5pt} j})$.
We define the {\it horizontal shift} of $M^{\ydt\hspace{.5pt}\ydt}$ to be the double complex $(X^{\ydt\hspace{.5pt}\ydt}, v_{\hspace{-1pt}_X}^{\ydt\hspace{.5pt}\ydt}, h_{\hspace{-1pt}_X}^{\ydt\hspace{.5pt}\ydt})$ such that
$X^{i,j}=M^{i+1,j}$, and $v_{\hspace{-1pt}_X}^{i, j}=-v_{\hspace{-1pt}_M}^{i+1,j}$ and $h_{\hspace{-1pt}_X}^{i,j}=-h_{\hspace{-1pt}_M}^{i+1, j}$. In the sequel,
we shall write $M^{\ydt\hspace{.5pt} \ydt}[1]$ for the horizontal shift of $M^{\ydt\hspace{.5pt}\ydt}$. A {\it morphism} $f^{\ydt \pdt}: M^{\ydt\,\ydt}\to N^{\ydt\,\ydt}$ of double complexes over $\cA$ consists of morphisms $f^{i,j}: M^{i,j}\to N^{i,j}$, with $i,j\in \Z$, in $\cA$ such that 
$$\xymatrixcolsep{14pt}\xymatrixrowsep{14pt}\xymatrix{
&N^{i,j+1}\\
M^{i,j+1} \ar[ur]^{f^{i,j+1}}&  \\
& N^{i,j} \ar[rr]^{h_{\hspace{-1pt}_N}^{i,j}}\ar[uu]_{v_{\hspace{-1pt}_N}^{i,j}}&& N^{i+1,j} \\
M^{i,j} \ar[uu]^{v_{\hspace{-1pt}_M}^{i,j}} \ar[ur]^-{f^{i,j}} \ar[rr]^-{h_{\hspace{-1pt}_M}^{i,j}} && M^{i+1,j}\ar[ur]_-{f^{i+1,j}}
}\vspace{2pt}$$
commutes, for all $i, j\in \Z$; or equivalently, $f^{i,\ydt}: M^{i,\ydt}\to N^{i,\ydt}$ and $f^{\ydt,j}: M^{\cdt,j}\to N^{\cdt,j}$ are complex morphisms, for all $i,j\in \Z$. We shall say that $f^{\ydt \pdt}$ is {\it horizontally null-homotopic} if
there exist $u^{i, j}: M^{i,j}\to N^{i-1, j}$ such that $u^{i+1, j} h_{\hspace{-1pt}_M}^{i,\hspace{.5pt} j}+ h_{\hspace{-1pt}_N}^{i-1,\hspace{.5pt} j}u^{i,j}=f^{i,j}$ and
$v_{\hspace{-1pt}_N}^{i-1, j} u^{i,j}+u^{i,j+1} v_{\hspace{-1pt}_M}^{i, j} =0$,
for all $i, j\in \Z$. The double complexes over $\cA$ form an abelian category $DC(\cA)$. 
Recall that the {\it total complex} $\mathbb{T}(M^{\ydt\hspace{.6pt} \ydt})$ of $M^{\ydt\hspace{.6pt} \ydt}$ is a complex over $\cA$ defined $\T(M^{\ydt\hspace{.6pt} \ydt})^n=\oplus_{i\in \Z} \, M^{i, n-i}$ and
$$d_{\T(M^{\ydt\hspace{.6pt} \ydt})}^{\,n}=(d_{\T(M^{\ydt\hspace{.6pt} \ydt})}^{\,n}(j,i))_{(j,i)\in \Z\times \Z}:
\oplus_{i\in \mathbb{Z}} \, M^{i, n-i} \to \oplus_{j\in \mathbb{Z}} \, M^{j, n+1-j},$$
where $d_{\T(M^{\ydt\hspace{.6pt} \ydt})}^{\,n}(j,i):  M^{i, n-i} \to M^{j, n+1-j}$ is given by
$$d_{\T(M^{\ydt\hspace{.6pt}\ydt})}^{\,n}(j,i) = \left\{\begin{array}{ll}
v_{\hspace{-1pt}_M}^{i, n-i}, & j=i;  \vspace{3pt} \\
h_{\hspace{-1pt}_M}^{i, n-i}, & j=i+1; \vspace{1pt} \\
0,            & j\ne i, i+1.
\end{array}\right.
$$

\medskip

\begin{Prop}\label{Tot}

Let $\cA$ be a concrete abelian category with countable direct sums. Taking total complexes yields a functor $\mathbb{T}: DC(\cA) \to C(\cA)$.

\begin{enumerate}[$(1)$]

\item If $M^{\ydt\hspace{.5pt} \ydt}\in DC(\cA)$, then $\T(M^{\ydt\hspace{.5pt} \ydt}[1])=\T(M^{\ydt\hspace{.5pt} \ydt})[1];$

\vspace{1.5pt}

\item If $f^{\ydt \pdt}: M^{\ydt\,\ydt}\to N^{\ydt\,\ydt}$ is horizontally null-homotopic, then $\T(f^{\ydt \pdt})$ is null-homotopic.

\end{enumerate}\end{Prop}

\noindent{\it Proof.} Let $f^{\pdt\hspace{.4pt} \pdt}: M^{\ydt\hspace{.6pt}\ydt}\to N^{\ydt\hspace{.6pt}\ydt}$ be a morphism in $DC(\cA)$. Given an integer $n$, we put
$$\T(f^{\ydt\hspace{.4pt} \ydt})^n=\left(\T\left(f^{\ydt \hspace{.5pt} \ydt}\right)^n\hspace{-3pt}(j, i)\right)_{(j,i)\in \mathbb{Z}\times \mathbb{Z}}:
\oplus_{i\in \Z} \, M^{i, n-i} \to \oplus_{j\in \Z} \,  N^{j, n-j},$$
where $\T(f^{\ydt\hspace{.4pt}\ydt})^n(i, i)=f^{i, n-i}$ and $\T(f^{\ydt\hspace{.5pt}\ydt})^n(j, i)=0$ in case $j\ne i$, or equivalently,  $\T(f^{\ydt\hspace{.4pt} \ydt})^n=\oplus_{i\in \Z}\, f^{i, n-i}.$
It is easy to see that $\T(f^{\ydt \hspace{.4pt} \ydt})^{n+1}\circ d_{\T(M^{\ydt \hspace{.4pt} \ydt})}^{\,n}= d^n_{\T(N^{\ydt\hspace{.4pt}\ydt})}\circ \T(f^{\ydt\hspace{.5pt} \ydt})^n$. This yields a morphism $\T(f^{\ydt \hspace{.4pt} \ydt})=(\T(f^{\ydt\hspace{.4pt} \ydt})^n)_{n\in \mathbb{Z}}:\T(M^{\ydt\hspace{.5pt} \ydt})\to \T(N^{\ydt\hspace{.5pt} \ydt})$ in $C(\cA)$. That is, $\T$ is indeed a functor. Statement (1) is evident.

Assume that $f^{\ydt \pdt}$ is horizontally null-homotopic. Let $u^{i, j}: M^{i,j}\to N^{i-1, j}$ be such that \vspace{1pt}$f^{i,j}=u^{i+1, j} \circ h_{\hspace{-1pt}_M}^{i,\hspace{.5pt} j}+ h_{\hspace{-1pt}_N}^{i-1,\hspace{.5pt} j} \circ u^{i,j}$ and $v_{\hspace{-1pt}_N}^{i-1, j} u^{i,j}+ u^{i,j+1} v_{\hspace{-1pt}_M}^{i, j} =0$, for all $i, j\in \Z$. Set $h^n=(h^n(j,i))_{(j,i)\in \Z\times \Z}: \oplus_{i\in \Z} M^{i, n-i}\to \oplus_{j\in \Z} N^{j, n-j}$, \vspace{1pt} where $h^n(i-1,i)=u^{i,n-i}: M^{i, n-i}\to N^{i-1, n-i}$ and $h^n(j,i)=0$ in case $j\ne i-1$. Given any $n, i, j\in \Z$, we obtain
$$\begin{array}{rcl}
\sum_{p\in \Z} h^{n+1}(j,p) \circ d^n_{\T(M^{\ydt \hspace{.5pt} \ydt})}(p,i) & = & h^{n+1}(j,j+1) \circ d^n_{\T(M^{\ydt \hspace{.5pt} \ydt})}(j+1,i)\vspace{5pt} \\
&=& \left\{\begin{array}{ll}
u^{i+1, n-i}\circ h_{\hspace{-1pt}_M}^{i,\hspace{.5pt} n-i},   & j=i; \\ \vspace{-10pt} \\
u^{i, n+1-i} \circ v_{\hspace{-1pt}_M}^{i, n-i},               & j=i-1; \\ \vspace{-10pt} \\
0,                                                             & j\ne i, i-1,
\end{array}\right.$$
\end{array}$$
and
$$\begin{array}{rcl}
\sum_{q\in \Z} d^{n-1}_{\T(N^{\ydt \hspace{.5pt} \ydt\hspace{.5pt}})}(j,q) \circ h^n(q,i) &=& d^{n-1}_{\T(N^{\ydt \hspace{.5pt} \ydt})}(j,i-1) \circ h^n(i-1,i) \vspace{5pt} \\
&=& \left\{\begin{array}{ll}
 h_{\hspace{-1pt}_N}^{i-1,\hspace{.5pt} n-i}\circ u^{i, n-i}, & j=i; \\ \vspace{-10pt} \\
 v_{\hspace{-1pt}_N}^{i-1, n-i} \circ u^{i, n-i},               & j=i-1; \\ \vspace{-10pt} \\
0,                                                            & j\ne i, i-1,
\end{array}\right.$$
\end{array}$$
from which we deduce that $\T(f^{\ydt \hspace{.5pt} \ydt})^n=h^{n+1}\circ d^n_{\T(M^{\ydt \hspace{.5pt} \ydt})}+ d^{n-1}_{\T(N^{\ydt \hspace{.5pt} \ydt\hspace{.5pt}})} \circ h^n.$ That is, $\T(f^{\ydt \hspace{.5pt} \ydt})$ is null-homotopic.
The proof of the proposition is completed.

\medskip

Next, we shall study when the total complex of a double complex is acyclic. Let $M^{\pdt \hspace{.4pt} \pdt}$ be a double complex over $\cA$. Given $n\in \Z$, the $n$-{\it diagonal} of $M^{\pdt \hspace{.4pt} \pdt}$ is the set of objects $M^{i, n-i}$ with $i\in \Z$. We shall say that $M^{\pdt \hspace{.4pt} \pdt}$ is {\it $n$-diagonally bounded} (respectively, {\it bounded-above}, {\it bounded-below}) if $M^{i, n-i}=0$ for all but finitely many (respectively, positive, negative) integers $i$. Moreover, $M^{\pdt \hspace{.4pt} \pdt}$ is called {\it diagonally bounded} (respectively, {\it bounded-above},  {\it bounded-below}) if it is $n$-diagonally bounded (respectively, bounded-above, bounded-below) for all $n\in \Z$.

\medskip

\begin{Lemma}\label{Homology-zero}

Let $\cA$ be a concrete abelian category with countable direct sums, and let $M^{\ydt \hspace{.4pt} \ydt} \in DC(\cA).$ If $n\in \Z$, then ${\rm H}^n(\T(M^{\ydt \hspace{.4pt} \ydt}))=0$ in case

\vspace{-2pt}

\begin{enumerate}[$(1)$]

\item $M^{\ydt \hspace{.4pt} \ydt}$ is $n$-diagonally bounded-below with ${\rm H}^{n-j}(M^{\ydt\hspace{.7pt}, \hspace{.5pt} j})=0$ for all $j\in \Z;$ or

\vspace{2pt}

\item $M^{\ydt \hspace{.4pt} \ydt}$ is $n$-diagonally bounded-above with ${\rm H}^{n-i}(M^{i, \ydt})=0$ for all $i\in \Z.$

\end{enumerate} \end{Lemma}

\noindent{\it Proof.} We shall consider only the case where Statement (1) holds for some $n$. In particular, there exists an integer $t<0$ such that $M^{i, n-i}=0$ for all $i<t$. Let $v^{\ydt \hspace{.5pt} \ydt}$ be the vertical differential, and $h^{\ydt \hspace{.4pt} \ydt}$ the horizontal differential, of $M^{\ydt \hspace{.4pt} \ydt}$. Write $(X^\ydt, d^\pdt)$ for the total complex of $M^{\ydt \hspace{.4pt} \ydt}$. Consider an element $c=(c_{i, n-i})_{i\in \mathbb{Z}}\in {\rm Ker}(d^n).$ Then, $v^{i, n-i}(c_{i, n-i})+h^{i-1, n-i+1}(c_{i-1, n-i+1})=0$, for all $i\in \Z$. Since $c$ has at most finitely many non-zero entries, we may assume that $c_{i, n-i}=0$ for all $i> 0$. Then, $h^{0,n}(c_{0,n})=-v^{1, n-1}(c_{1, n-1})=0$. Since ${\rm H}^0(M^{\ydt\hspace{.5pt}, n})=0$, there exists some $x_{-1,n}\in M^{-1, n}$ such that $c_{0,n}=h^{-1,n}(x_{-1,n})$. This yields that
$$h^{-1, n+1}(c_{-1, n+1} - v^{-1, n}(x_{-1,n})) = h^{-1, n+1}(c_{-1, n+1}) + v^{0, n}(c_{0,n})=0.$$

Since ${\rm H}^{-1}(M^{\ydt\hspace{.6pt},n+1})=0$, we see that $c_{-1, n+1}-v^{-1, n}(x_{-1,n})=h^{-2,n+1}(x_{-2,n+1})$, with
$x_{-2,n+1}\in M^{-2, n+1}$. Continuing this process, we find $x_{-i, n-1+i}\in M^{i, n-1-i}$ such that
$c_{i, n-i}=v^{i, n-1-i}(x_{i,n-1-i})+h^{i-1,n-i}(x_{i-1,n-i})$, for $i=-1, \ldots, t$. Since $M^{t-1, n-t}=0$, we see that $h^{t, n-1-t}(x_{t, n-1-t})=0$. Setting $x=(x_{i, n-1-i})_{i\in \mathbb{Z}}$, where $x_{i, n-1-i}=0$ for $i\geq 0$ or $i<t$, we obtain $c=d^{n-1}(x)$. The proof of the lemma is completed.


\medskip

As an immediate consequence of Lemma \ref{Homology-zero}, we obtain the promised generalization of the Acyclic Assembly Lemma stated in \cite[(2.7.1)]{Wei}.

\medskip

\begin{Cor}\label{AAL}

Let $\cA$ be a concrete abelian category with countable direct sums. If $M^{\ydt \hspace{.6pt} \ydt}\in DC(\cA)$, then $\T(M^{\ydt \hspace{.4pt} \ydt})$ is acyclic in case $M^{\ydt \hspace{.6pt} \ydt}$ is diago\-nally bounded-below with acyclic rows or diagonally bounded-above with acyclic columns.

\end{Cor}

\medskip

Given a double complex morphism $f^{\ydt \pdt}: M^{\ydt\hspace{.6pt}\ydt}\to N^{\ydt\,\ydt},$ we shall define a {\it horizontal cone} $H_{\hspace{-1.3pt}f^{\ydt \pdt}}$ and a {\it vertical cone} $V_{\hspace{-1pt}f^{\ydt \pdt}}$. Indeed, $H_{\hspace{-1.3pt}f^{\ydt \pdt}}$ is the double complex $(H^{\ydt\,\ydt},v^{\ydt\,\ydt},h^{\ydt\,\ydt})$ defined by
$H^{i,j}=M^{i+1,j}\oplus N^{i,j}$ and
$$
v^{i,j}=\left(\begin{array}{cc}
-v_{\hspace{-1pt}_M}^{i+1,j} &  0  \vspace{2pt} \\
         0   & v_{\hspace{-1pt}_N}^{i,j}
\end{array}\right), \;
h^{i,j}=\left(
\begin{array}{cc}
-h_{\hspace{-1pt}_M}^{i+1,j} &  0 \vspace{2.5pt} \\
f^{\,i+1,j} &  h_{\hspace{-1pt}_N}^{i,j}
\end{array}\right),
\vspace{2pt}$$
which can be visualized as follows:
$$\xymatrixcolsep{25pt}\xymatrixrowsep{20pt}
\xymatrix{
             &  \vdots                                                    &&  &           \vdots \\
\cdots \ar[r] & M^{i+1,j+1}\oplus N^{i,j+1} \ar[u]
\ar[rrr]^---{\left(\hspace{-5pt}\begin{array}{cc} -h_{\hspace{-1pt}_M}^{i+1,j+1} & 0 \vspace{1pt} \\ f^{\,i+1,j+1}  & h_{\hspace{-1pt}_N}^{i,j+1} \end{array}\hspace{-5pt}\right)}
&&& M^{i+2,j+1}\oplus N^{i+1,j+1} \ar[u] \ar[r] &                                            \cdots \\ \\
\cdots \ar[r] &  M^{i+1,j}\oplus N^{i,j} \hspace{10pt} \ar[rrr]^---{\left(\hspace{-5pt}\begin{array}{cc} -h_{\hspace{-1pt}_M}^{i+1,j} &  0 \vspace{1pt} \\ f^{\,i+1,j}  & h_{\hspace{-1pt}_N}^{i,j} \end{array}\hspace{-5pt}\right)}
\ar[uu]^{\left(\hspace{-5pt}\begin{array}{cc} -v_{\hspace{-1pt}_M}^{i+1,j} &  0 \vspace{2pt} \\ 0  &  v_{\hspace{-1pt}_N}^{i,j} \end{array}\hspace{-5pt}\right)}
&&& M^{i+2,j} \oplus N^{i+1,j}\ar[uu]_{\left(\hspace{-5pt}\begin{array}{cc} -v_{\hspace{-1pt}_M}^{i+2,j} &  0 \vspace{2pt} \\ 0  &  v_{\hspace{-1pt}_N}^{i+1,j} \end{array}\hspace{-5pt}\right)}
\ar[r] & \cdots \\ &\vdots\ar[u] &&& \vdots\ar[u]
}\vspace{4pt}$$
In a similar fashion, one can define the vertical cone $V_{\hspace{-1pt}f^{\ydt \pdt}}$ of $f^{\ydt\,\ydt}$. The following observation will be useful.

\medskip

\begin{Lemma}\label{V-cone}

Let $\cA$ be a concrete abelian category with countable direct sums, and let $f^{\ydt \pdt}: M^{\ydt\hspace{.6pt}\ydt}\to N^{\ydt\,\ydt}$ be a morphism in $DC(\cA)$.

\begin{enumerate}[$(1)$]

\item If $f^{i,\ydt}: M^{i,\ydt}\to N^{i,\pdt}$ is a quasi-isomorphism for some $i\in \Z$, then the $i$-th column of $V_{\hspace{-1pt}f^{\ydt \pdt}}$ is acyclic.

\vspace{1pt}

\item If $f^{\ydt, j}: M^{\ydt, j}\to N^{\ydt, j}$ is a quasi-isomorphism for some $j\in \Z$, then the $j$-th row of $H_{\hspace{-1.3pt}f^{\ydt \pdt}}$ is acyclic.

\end{enumerate}

\end{Lemma}

\noindent{\it Proof.} In view of the definition, we see that the $i$-th column of $V_{f^{\ydt\,\ydt}}$ is the mapping cone $C_{f^{i,\ydt}}$ of $f^{i,\ydt}: M^{i,\ydt}\to N^{i, \ydt},$ while the $j$-th row of $H_{\hspace{-1pt}f^{\ydt \pdt}}$ is the mapping cone $C_{f^{\ydt, j}}$ of $f^{\ydt, j}: M^{\ydt, j}\to N^{\ydt, j}.$ The proof of the lemma is completed.

\medskip

The following results says that $\T: DC(\cA)\to C(\cA)$ sends vertical cones and horizontal cones to mapping cones.

\medskip

\begin{Lemma}\label{Cone}

Let $\cA$ be concrete abelian category with countable direct sums. If $f^{\ydt\hspace{.5pt}\ydt}: M^{\ydt\,\ydt}\to N^{\ydt\,\ydt}$ in $DC(\cA)$,
then $\T(H_{\hspace{-1.3pt}f^{\ydt\hspace{.5pt}\ydt}})=C_{\mathbb{T}(f^{\ydt\hspace{.5pt}\ydt})}=\T(V_{{f^{\ydt\hspace{.5pt}\ydt}}}).$

\end{Lemma}

\noindent{\it Proof.} Let $f^{\ydt\hspace{.5pt}\ydt}: M^{\ydt\hspace{.5pt}\ydt}\to N^{\ydt\hspace{.5pt}\ydt}$ be a morphism in $DC(\cA)$.
Given any $n\in \Z$, we obtain $\T(H_{\hspace{-1.3pt}f^{\ydt\hspace{.5pt}\ydt}})^n=\oplus_{i\in \Z} (M^{i+1,n-i}\oplus N^{i,n-i})$ and \vspace{1pt}
$d^n_{\mathbb{T}(H_{\hspace{-1pt}f^{\ydt\hspace{.5pt}\ydt}})}=(d^n_{\mathbb{T}(H_{\hspace{-1.3pt}f^{\ydt\hspace{.5pt}\ydt}})}(j,i))_{(j,i)\in \mathbb Z\times \mathbb Z},$
where $d^n_{\mathbb{T}(H_{\hspace{-1.3pt}f^{\ydt\hspace{.5pt}\ydt}})}(j,i): M^{i+1,n-i}\oplus N^{i,n-i}  \to M^{j+1,n+1-j}\oplus N^{j,n+1-j}$ is defined by
$$
\begin{array}{rcl}
d^n_{\mathbb{T}(H_{\hspace{-1.5pt}f^{\ydt\hspace{.5pt}\ydt}})}(j,i) &=& \left\{\begin{array}{ll}
\hspace{-4pt}\left(\begin{array}{cc}
-v_{\hspace{-1pt}_M}^{i+1,n-i}      &  0              \vspace{3pt}\\
0                                   & v_{\hspace{-1pt}_N}^{i,n-i}
\end{array}\right) \vspace{3pt},                                              & j=i;   \vspace{3pt} \\
\hspace{-4pt}\left(\begin{array}{cc}
-h_{\hspace{-1pt}_M}^{i+1,n-i} &  0              \vspace{3pt}\\
f^{\,i+1,n-i}                  &  h_{\hspace{-1pt}_N}^{i,n-i}
\end{array}\right) \vspace{3pt},                                              & j=i+1; \vspace{3pt} \\
0,                                                                            & j\neq i,i+1.
\end{array}\right.
\end{array}
\vspace{3pt}
$$
On the other hand, $\mathbb{T}(f^{\ydt\hspace{.5pt}\ydt}): \T(M^{\cdt \cdt})\to \T(N^{\cdt \cdt})$ is a morphism in $C(\cA)$, whose mapping cone $C_{\mathbb{T}(f^{\ydt\hspace{.5pt}\ydt})}$ is defined by
$$
C_{\mathbb{T}(f^{\ydt\hspace{.5pt}\ydt})}^n = \T(M^{\ydt\hspace{.5pt}\ydt})^{n+1} \oplus \T(N^{\ydt\hspace{.5pt}\ydt})^n = \oplus_{i\in \Z} (M^{i+1,n-i}\oplus N^{i,n-i})=\T(H_{f^{\ydt\hspace{.5pt}\ydt}})^n,
$$
and
$$d^n_{C_{_{\mathbb{T}(f^{\ydt\hspace{.5pt}\ydt})}}}=\left(\hspace{-2pt}\begin{array}{cc}
-d^{n+1}_{\T(M^{\ydt\hspace{.5pt}\ydt})}     &\hspace{-5pt}  0     \vspace{3pt}\\
\mathbb{T}(f^{\ydt\hspace{.5pt}\ydt})^{n+1}  &\hspace{-5pt} d^n_{\T(N^{\ydt\hspace{.5pt}\ydt})}
\end{array}\hspace{-2pt}\right)=(d^n_{_{C_{\mathbb{T}(f^{\ydt\hspace{.5pt}\ydt})}}}\hspace{-1pt}(j,i))_{(j,i)\in \mathbb Z\times \mathbb Z},$$

\vspace{5pt}

\noindent where $d^n_{_{C_{\mathbb{T}(f^{\ydt\hspace{.5pt}\ydt})}}}\hspace{-1pt}(j,i): M^{i+1,n-i}\oplus N^{i,n-i}\to M^{j+1,n+1-j}\oplus N^{j,n+1-j}$ is defined by
$$
\left(\hspace{-5pt}\begin{array}{cc}
-d^{n+1}_{\T(M^{\ydt\hspace{.5pt}\ydt})}(j,i)    &  0 \vspace{5pt} \\
\mathbb{T}(f^{\ydt\hspace{.5pt}\ydt})^{n+1}(j,i) &  d^n_{\T(N^{\ydt\hspace{.5pt}\ydt})}(j,i)
\end{array}\hspace{-5pt}\right)
=
\left\{
\begin{array}{ll}
\hspace{-3pt}
\left(\hspace{-2pt}\begin{array}{cc}
-v_{\hspace{-1pt}_M}^{i+1,n-i} & 0 \vspace{3pt}\\
0                              &  v_{\hspace{-1pt}_N}^{i,n-i}
\end{array}\hspace{-2pt}\right),                         &   j=i; \vspace{5pt}\\
\hspace{-3pt}
\left(\hspace{-2pt}\begin{array}{cc}
-h_{\hspace{-1pt}_M}^{i+1,n-i} &  0  \vspace{5pt}\\
f^{\,i+1,n-i}                             &  h_{\hspace{-1pt}_N}^{i,n-i}
\end{array}\hspace{-2pt}\right) \vspace{3pt},            &   j=i+1;   \vspace{1pt}\\
0,                                                       &   j\neq i,i+1.
\end{array}\right.
$$
Thus, $d^n_{_{C_{\mathbb{T}(f^{\ydt\hspace{.5pt}\ydt})}}}\hspace{-1pt}(j,i)= d^n_{\mathbb{T}(H_{\hspace{-1.5pt}f^{\ydt\hspace{.5pt}\ydt}})}(j,i),$ for $i,j\in \Z$. This establishes the first part of the lemma, and the second part can be shown in a similar fashion. The proof of the lemma is completed.

\medskip

Let $\cA, \cB$ are concrete abelian categories with countable direct sums. Consider a functor $\mathfrak{F}: \cA\to C(\cB): M\to \mathfrak{F}(M)^\ydt; f\mapsto \mathfrak{F}(f)^\ydt.$ We shall extend it to the complexes over $\cA$. Given $M^\ydt\in C(\cA)$, applying $\mathfrak{F}$ to each of its components yields a double complex $\mathfrak{F}(M^\pdt)^\cdt$ over $\cB$ as follows: \vspace{-5pt}

$$\xymatrixrowsep{20pt}\xymatrixcolsep{22pt}\xymatrix{
             &  \vdots                                                    &&             \vdots\\
\cdots\ar[r] & \mathfrak{F}(M^i)^{j+1} \ar[u] \ar[rr]^{\mathfrak{F}(d_{\hspace{-.5pt}M}^i)^{j+1}}  && \mathfrak{F}(M^{i+1})^{j+1} \ar[u] \ar[r]
& \cdots\\
\cdots\ar[r] & \mathfrak{F}(M^i)^j \ar[rr]^{\mathfrak{F}(d_{\hspace{-.5pt}M}^{\hspace{.5pt}i})^j} \ar[u]^{(-1)^i d^j_{\mathfrak{F}(M^i)}} && \mathfrak{F}(M^{i+1})^j \ar[u]_{(-1)^{i+1} d_{\mathfrak{F}(M^{i+1})}^j} \ar[r] & \cdots \\ & \vdots\ar[u] && \vdots\ar[u] \\
}\vspace{0pt}$$

Note that the $i$-th column of $\mathfrak{F}(M^\cdt)^\cdt$ is $\frak{t}^i(\mathfrak{F}(M^i)^\cdt)$, the $i$-th twist of $\mathfrak{F}(M^i)^\cdt$. Given a morphism $f^\ydt: M^\ydt\to N^\ydt$ in $C(\cA)$, the diagram
$$\xymatrixcolsep{14pt}\xymatrixrowsep{14pt}\xymatrix{
& \fF(N^i)^{j+1}\\
\fF(M^i)^{j+1} \ar[ur]^{\fF(f^i)^{j+1}}&  \\
& \fF(N^i)^j \ar[rr]^{\fF(d_N^i)^j}\ar[uu]_{(-1)^id^j_{\fF(N^i)}}&& \fF(N^{i+1})^j, \\
\fF(M^i)^j. \ar[uu]^{(-1)^id^j_{\fF(M^i)}}\ar[ur]^-{\fF(f^i)^j} \ar[rr]^-{\fF(d_M^i)^j} && \fF(M^{i+1})^j\ar[ur]_-{\fF(f^{i+1)^j}}
}\vspace{5pt}$$
is commutative, for all $i,j\in \Z$. Therefore, $\fF(f^i)^j$ with $i,j\in \Z$ form a morphism $\fF(f^\ydt)^\ydt: \fF(M^\ydt)^\ydt \to \fF(N^\ydt)^\ydt$ in $DC(\cB)$.

\medskip

\begin{Prop}\label{cal-F}

Let $\cA, \cB$ be concrete abelian categories with countable direct sums. A functor $\fF: \cA\to C(\cB)$ induces a functor \vspace{0pt}
$$\fF^{DC}: C(\cA) \to DC(\cB): M^\ydt\mapsto \fF(M^\ydt)^\ydt\hspace{.5pt}; \hspace{.5pt} f^\ydt\mapsto \fF(f^\ydt)^\ydt.
\vspace{-1pt}$$


\begin{enumerate}[$(1)$]

\item If $M^{\ydt}$ is an object in $C(\cA)$, then $\fF^{DC}(M^{\ydt}[1]) = \fF^{DC}(M^{\ydt})[1].$

\vspace{1.5pt}

\item If $f^{\ydt}$ is a morphism in $C(\cA)$, then $\fF^{DC}(C_{\hspace{-1pt}f^{\ydt}})=H_{\mathfrak{F}^{DC}(f^{\pdt})}$
and $\fF^{DC}(f^{\ydt})$ is horizontally null-homotopic whenever $f^{\ydt}$ is null-homotopic.

\end{enumerate}

\end{Prop}

\noindent{\it Proof.} Let $\fF: \cA\to C(\cB)$ be a functor. Statement (1) is evident. Let $f^\ydt: M^\ydt\to N^\cdt$ be a morphism in $C(\cA)$. Write $(C^\pdt, d_{_C}^{\hspace{.5pt}\ydt})$ for its mapping cone. For any $n\in \Z$, we obtain $\fF(C^n)^\cdt=\fF(M^{n+1})^\cdt \oplus \fF(N^n)^\cdt$ with \vspace{0pt}
$$
d_{\fF(C^n)}^{\pdt}=\left(\begin{array}{cc}
d_{\fF(M^{n+1})}^\pdt & \hspace{-8pt} 0  \\ \vspace{-8pt}\\
                    0 & \hspace{-8pt} d_{\fF(N^n)}^{\pdt}
\end{array}\right); \quad \fF(d_{_C}^{\hspace{.6pt}n})^\cdt=\left(\begin{array}{cc}
- \fF(d_{\hspace{-.8pt}_M}^{n+1})^\cdt & \hspace{-8pt} 0\\ \vspace{-8pt}\\
\fF(f^{n+1})^\cdt & \hspace{-8pt} \fF(d_{\hspace{-.8pt}_N}^n)^{\hspace{.6pt}\cdt}
\end{array}\right).
$$

\vspace{2pt}

Write $(H^{\pdt\hspace{.5pt}\pdt}, v^{\ydt\hspace{.5pt}\ydt}_{\hspace{-1pt}_H}, h^{\ydt\hspace{.5pt}\ydt}_{\hspace{-1pt}_H})$ for the horizontal cone of $\mathfrak{F}^{DC}(f^{\pdt}):
\mathfrak{F}^{DC}(M^{\pdt})\to \mathfrak{F}^{DC}(N^{\pdt}),$ where $\fF^{DC}(M^\ydt)=F(M^\cdt)^\cdt$. Given $i,j\in \Z$, by definition, we obtain
$$H^{i,j} = \fF(M^{i+1})^{j} \oplus \fF(N^i)^{j}= \fF(C^i)^j  = \fF^{DC}(C^\ydt)^{i, j}$$
with
$$
\begin{array}{rcl}
h_{\hspace{-1pt}_H}^{i,j}
=
\left(\hspace{-5pt}\begin{array}{cc}
-\fF(d_{M^{i+1}})^j           & \hspace{-3pt} 0 \vspace{5pt} \\
\fF(f^{i+1})^j      & \hspace{-3pt} \fF(d_{N^{i}})^j
\end{array}\hspace{-5pt}\right) =
\fF(d_C^{i})^j=h_{\hspace{-2pt}_{\fF^{DC}(C^\ydt)}}^{i,j}
\end{array}
$$
and
$$
v^{i,j}_{\hspace{-1pt}_H} = \left(\hspace{-5pt}\begin{array}{cc}
(-1)^id_{\fF(M^{i+1})}^j & \hspace{-3pt} 0 \vspace{5pt} \\
0                        & \hspace{-3pt} (-1)^id_{\fF(N^i)}^j
\end{array}\hspace{-5pt}\right) = (-1)^id^j_{\fF(C^i)}=v_{\hspace{-2pt}_{\fF^{DC}(C^\ydt)}}^{i,j}. \vspace{2pt}
$$
That is, $C^\pdt=H^\cdt$. Next, suppose that $f^\cdt$ is null-homotopic. Let $u^i: M^i\to N^{i-1}$\vspace{1pt} be morphisms such that $f^i=u^{i+1}\circ d_{\hspace{-1pt}_M}^i+d_{\hspace{-1pt}_N}^{i-1}\circ u^i$, for all $i\in \Z$. In particular, $\fF(f^i)^j=\fF(u^{i+1})^j\circ \fF(d_{\hspace{-1pt}_M}^i)^j+\fF(d_{\hspace{-1pt}_N}^{i-1})^j\circ \fF(u^i)^j,$ for all $j\in \Z$. Moreover, since $\fF(u^i)^\cdt: \fF(M^i)\to \fF(N^{i-1})^\cdt$ is a complex morphism, we obtain
$$(-1)^i \fF(u^i)^{j+1}\circ d_{\fF(M^i)}^j  + (-1)^i d_{\fF(N^{i-1})}^j \circ \fF(u^i)^j=0,$$ for all $j\in \Z$.
Considering $F(u^i)^j: \fF(M^i)^j\to \fF(N^{i-1})^j$ with $i, j\in \Z$, we see that $\fF^{DC}(f^{\ydt})$ is horizontally null-homotopic. The proof of the proposition is completed.

\medskip

The following statement follows immediately from Propositions \ref{Tot} and \ref{cal-F}, which is a general version of Lemma 3.7 stated in \cite{BaL2}.

\medskip

\begin{Prop}\label{F-extension}

Let $\cA, \cB$ be concrete abelian categories with countable direct sums. A functor $\fF: \cA\to C(\cB)$ extends to a functor
$\fF^C=\T\circ \fF^{DC}: C(\cA) \to C(\cB).$

\begin{enumerate}[$(1)$]

\item If $M$ is an object in $\cA$, then $\fF^{C}(M) = \fF(M)^\ydt.$

\vspace{1.5pt}

\item If $M^\ydt$ is a complex in $C(\cA)$, then $\fF^{C}(M^\ydt[1]) = \fF^{C}(M^\ydt)[1].$

\vspace{1.5pt}

\item If $f^\ydt$ is a morphism in $C(\cA)$, then $\fF^{C}(C_{\hspace{-1pt}f^\pdt})=C_{\mathfrak{F}^C(f^\pdt)}$ and  $\fF^{C}(f^\ydt)$  is null-homotopic whenever $f^\ydt$ is null-homotopic.
    
\end{enumerate}

\end{Prop}

\medskip

The following statement says that such extension of functors is compatible with the composition of functors.

\medskip

\begin{Lemma}\label{F-composition}

Let $\cA, \cB, \cC$ be concrete abelian categories with countable direct sums. If $\fF: \cA\to C(\cB)$ and $\mathfrak{G}: \cB\to C(\cC)$ are functors, then $(\mathfrak{G}^C\circ \mathfrak{F})^C=\mathfrak{G}^C\circ \mathfrak{F}^C.$

\end{Lemma}

\noindent{\it Proof.} Let $\fF: \cA\to C(\cB)$ and $\mathfrak{G}: \cB\to C(\cC)$ be functors. Fix some complex $M^\ydt\in C(\cA)$. 
Given an integer $n$, we obtain
$(\mathfrak{G}^C\circ \mathfrak{F})^C(M^\pdt)^n=\oplus_{i\in \Z}\, \mathfrak{G}^C(\mathfrak{F}(M^i)^\ydt)^{n-i}
$
and $d^n_{(\mathfrak G^C\circ \mathfrak F)^C\left(M^\ydt\right)}=
(d^n_{(\mathfrak G^C\circ \mathfrak F)^C(M^\ydt)}(j,i))_{(j,i)\in \mathbb Z^2},$ where
$$d^n_{(\mathfrak G^C\circ \mathfrak F)^C(M^\ydt)}(j,i): \mathfrak G^C\hspace{-1.5pt}(\mathfrak F(M^i)^\ydt)^{n-i}\to \mathfrak G^C\hspace{-1.5pt}(\mathfrak F(M^j)^\ydt)^{n+1-j}$$ is given by
%
%
$$d^n_{(\mathfrak G^C\circ \mathfrak F)^C(M^\ydt)}(j,i)=\left\{\begin{array}{ll}
(-1)^i d_{\mathfrak G^C\hspace{-1.5pt}(\mathfrak F(M^i)^\ydt)}^{n-i}, & j=i; \vspace{2.5pt}\\
\mathfrak G^{\hspace{-1pt}C\hspace{-1pt}}(\mathfrak F(d_M^i)^\ydt)^{n-i}, & j=i+1 \\
0, & j\ne i, i+1.
\end{array}\right.\vspace{3pt}$$
Furthermore, by definition, we obtain a diagram
$$\xymatrixcolsep{18pt}\xymatrixrowsep{14pt}\xymatrix{
\mathfrak G^C(\mathfrak F(M^i)^\ydt)^{n-i} \ar[rrrrr]^-{\mathfrak G^{\hspace{-.6pt}C\hspace{-1pt}}(\mathfrak F(d_M^i)^{\hspace{-.5pt}\cdt\hspace{.5pt}})^{n-i}} \ar@{=}[d]&&&&& \mathfrak G^C(\mathfrak F(M^{i+1})^\ydt)^{n-i} \ar@{=}[d] \\
\oplus_{p\in \Z\,} \mathfrak G(\mathfrak F(M^i)^p)^{n-i-p}
\ar[rrrrr]^-{\left(\mathfrak G^{\hspace{-.6pt}C\hspace{-1pt}}(\mathfrak F(d_M^i)^{\hspace{-.5pt}\cdt\hspace{.5pt}})^{n-i}(q,p)\right)_{(q,p)\in \mathbb Z^2}}
&&&&& \oplus_{q\in \Z\,} \mathfrak G(\mathfrak F(M^{i+1})^q)^{n-i-q},}$$
where
$$\mathfrak G^{\hspace{-.6pt}C\hspace{-1pt}}(\mathfrak F(d_M^i)^{\hspace{-.5pt}\cdt\hspace{.5pt}})^{n-i}(q,p)
=\left\{\begin{array}{ll} \mathfrak G (\mathfrak F(d_M^i)^p)^{n-i-p}, & q=p;\\
0, & q\ne p, \end{array}\right.
\vspace{2pt}$$
and a diagram
$$\xymatrixcolsep{16pt}\xymatrixrowsep{14pt}\xymatrix{
\mathfrak G^C(\mathfrak F(M^i)^\ydt)^{n-i} \ar[rrrrr]^-{d^{n-i}_{\mathfrak G^{\hspace{-.4pt}C\hspace{-1pt}}(\mathfrak F(M^i)^\ydt)}} \ar@{=}[d]&&&&& \mathfrak G^C(\mathfrak F(M^i)^\ydt)^{n+1-i} \ar@{=}[d] \\
\oplus_{p\in \Z\,} \mathfrak G(\mathfrak F(M^i)^p)^{n-i-p}
\ar[rrrrr]^-{\left(d^{n-i}_{\mathfrak G^{\hspace{-.4pt}C\hspace{-1pt}}(\mathfrak F(M^i)^\ydt)}(q,p)\right)_{(q,p)\in \mathbb Z^2}}
&&&&& \oplus_{q\in \Z\,} \mathfrak G(\mathfrak F(M^i)^q)^{n+1-i-q},}$$
where

$$
d^{n-i}_{\mathfrak G^{\hspace{-.4pt}C\hspace{-1pt}}(\mathfrak F(M^i)^\ydt)}(q,p)=\left\{\begin{array}{ll}
(-1)^p d_{\mathfrak{G}(\mathfrak{F}(M^i)^p)}^{\hspace{.4pt}n-i-p}, & q=p;   \\ \vspace{-9pt} \\
\mathfrak G(d_{\mathfrak{F}(M^i)}^p)^{n-i-p},             & q=p+1; \\ \vspace{-9pt} \\
0,                          & q\ne p, p+1.
\end{array} \right.
\vspace{5pt} $$
Therefore, $(\mathfrak{G}^C\circ \mathfrak{F})^C(M^\ydt)$ is the complex described by the diagram
$$\xymatrixcolsep{18pt}\xymatrixrowsep{14pt}\xymatrix{
(\mathfrak G^C\circ \mathfrak F)^C(M^\ydt)^n \ar[rrrrr]^-{d^n_{(\mathfrak G^{\hspace{-.5pt}C\hspace{-.6pt}}\circ \mathfrak F)^{\hspace{-.5pt}C\hspace{-.6pt}}(M^\ydt)}} \ar@{=}[d] &&&&&  (\mathfrak G^C\circ \mathfrak F)^C(M^\ydt)^{n+1} \ar@{=}[d] \\
\oplus_{(i,p)\in \mathbb{Z}^2} \mathfrak{G}(\mathfrak{F}(M^i)^p)^{n-i-p}
\ar[rrrrr]^-{(d^n_{(\mathfrak G^{\hspace{-.5pt}C\hspace{-.6pt}}\circ \mathfrak F)^{\hspace{-.5pt}C\hspace{-.6pt}}(M^\ydt)}(j,q;i,p))_{(j,q;i,p)\in \Z^4}} &&&&& \oplus_{(j,q)\in \mathbb{Z}^2} \mathfrak{G}(\mathfrak{F}(M^j)^q)^{n+1-j-q}}$$
where
$$
d^n_{(\mathfrak G^{\hspace{-.5pt}C\hspace{-.6pt}}\circ \mathfrak F)^{\hspace{-.5pt}C\hspace{-.6pt}}(M^\ydt)}(j,q;i,p)
= \left\{\begin{array}{ll}
(-1)^{i+p} d_{\mathfrak G(\mathfrak F(M^i)^p)}^{n-i-p}, & j=i; q=p;  \vspace{3pt}  \\
(-1)^i \mathfrak G(d_{\mathfrak F(M^i)}^p)^{n-i-p},     & j=i; q=p+1; \vspace{3pt} \\
\mathfrak G(\mathfrak F(d_M^i)^p)^{n-i-p}               & j=i+1, q=p; \vspace{3pt} \\
0,                                                      & \text{otherwise}.
\end{array}\right.
$$

\smallskip

Next, given any integer $n$, we obtain
$\mathfrak{G}^C(\mathfrak{F}^C(M^\ydt))^n=\oplus_{s\in \Z}\mathfrak{G}(\mathfrak{F}^C(M^\ydt)^s)^{n-s}
$
and $d^n_{\mathfrak{G}^{\hspace{-.5pt}C\hspace{-.6pt}}(\mathfrak{F}^{\hspace{-.5pt}C\hspace{-.6pt}}(M^\pdt))}=
(d^n_{\mathfrak{G}^{\hspace{-.5pt}C\hspace{-.6pt}}(\mathfrak{F}^{\hspace{-.5pt}C\hspace{-.6pt}}(M^\pdt))}(t,s))_{(t,s)\in \mathbb Z^2}$, where
$$d^n_{\mathfrak{G}^{\hspace{-.5pt}C\hspace{-.6pt}}(\mathfrak{F}^{\hspace{-.5pt}C\hspace{-.6pt}}(M^\pdt))}(t,s):
\mathfrak{G}(\mathfrak{F}^C(M^\ydt)^s)^{n-s} \to \mathfrak{G}(\mathfrak{F}^C(M^\ydt)^t)^{n+1-t}$$
is given by
$$d^n_{\mathfrak{G}^C(\mathfrak{F}^C(M^\pdt))}(t,s) = \left\{\begin{array}{ll}
(-1)^s d_{\mathfrak G (\mathfrak F^C(M^\pdt)^s)}^{\hspace{.4pt} n-s}, & t=s;   \\\vspace{-9pt} \\
\mathfrak G(d_{\mathfrak F^C(M^\pdt)}^s)^{n-s},             & t=s+1; \\\vspace{-9pt} \\
0,                          & t\ne s, s+1.
\end{array}\right.$$
Furthermore, by definition, we obtain a diagram
$$\xymatrixcolsep{16pt}\xymatrixrowsep{14pt}\xymatrix{
\mathfrak G(\mathfrak F^C(M^\ydt)^s)^{n-s} \ar@{=}[d] \ar[rrrrr]^-{d_{\mathfrak G (\mathfrak F^C(M^\ydt)^s)}^{\hspace{.4pt} n-s}} &&&&& \mathfrak G(\mathfrak F^C(M^\ydt)^s)^{n-s+1}\ar@{=}[d] \\
\oplus_{i\in \mathbb Z} \,\mathfrak{G}(\mathfrak{F}(M^i)^{s-i})^{n-s}
\ar[rrrrr]^-{(d_{\mathfrak G (\mathfrak F^C(M^\ydt)^s)}^{\hspace{.4pt} n-s}(j,i))_{(j,i)\in \mathbb Z\times \mathbb \Z}}
&&&&& \oplus_{j\in \mathbb Z} \,\mathfrak{G}(\mathfrak{F}(M^j)^{s-j})^{n-s+1},}$$
where

$$d_{\mathfrak G (\mathfrak F^C(M^\ydt)^s)}^{\hspace{.4pt} n-s}(j,i)=
\left\{\begin{array}{ll}
d_{\mathfrak G(\mathfrak F(M^i)^{s-i})}^{n-s}, & j=i;\\
0, & j\ne i
\end{array} \vspace{5pt} \right. $$
and a diagram
$$\xymatrixcolsep{16pt}\xymatrixrowsep{14pt}\xymatrix{
\mathfrak G\left(\mathfrak F^C(M^\ydt)^s\right)^{n-s} \ar@{=}[d] \ar[rrrrr]^{\mathfrak G (d_{\mathfrak F^C(M^\pdt)}^s )^{n-s}} &&&&& \mathfrak G\left(\mathfrak F^C(M^\ydt)^{s+1}\right)^{n-s} \ar@{=}[d] \\
\oplus_{i\in \Z} \mathfrak G(\mathfrak F(M^i)^{s-i})^{n-s}  \ar[rrrrr]^{(\mathfrak G(d_{\mathfrak F^C(M^\ydt)}^s)^{n-s}(j,i))_{(j,i)\in \mathbb Z\times \mathbb Z} } &&&&& \oplus_{j\in \Z} \mathfrak G(\mathfrak F(M^j)^{s+1-j})^{n-s}, }$$
where

$$\mathfrak G(d_{\mathfrak F^C(M^\ydt)}^s)^{n-s}(j,i)=
\left\{\begin{array}{ll}
(-1)^i \, \mathfrak G (d_{\mathfrak F(M^i)}^{s-i})^{n-s}, & j=i;\vspace{2pt}\\
\mathfrak G( \mathfrak F(d_M^i)^{s-i})^{n-s}, & j=i+1;\vspace{1.5pt}\\
0, & j\ne i, i+1.
\end{array}\right.\vspace{5pt}$$
Thus, $\mathfrak{G}^C(\mathfrak{F}^C(M^\ydt))$ is the complex described by the diagram
$$\xymatrixcolsep{16pt}\xymatrixrowsep{14pt}\xymatrix{
\mathfrak{G}^{\hspace{-.5pt}C\hspace{-.6pt}}(\mathfrak{F}^{\hspace{-.5pt}C\hspace{-.6pt}}(M^\ydt))^n \ar@{=}[d] \ar[rrrrr]^{d^n_{\mathfrak{G}^{\hspace{-.5pt}C\hspace{-.6pt}}(\mathfrak{F}^{\hspace{-.5pt}C\hspace{-.6pt}}(M^\ydt))}} &&&&& \mathfrak{G}^{\hspace{-.5pt}C\hspace{-.6pt}}(\mathfrak{F}^{\hspace{-.5pt}C\hspace{-.6pt}}(M^\ydt))^{n+1} \ar@{=}[d] \\
\oplus_{(i,s)\in \mathbb Z^2\,} \mathfrak{G}(\mathfrak{F}(M^i)^{s-i})^{n-s}
\ar[rrrrr]^-{(d^n_{\mathfrak{G}^{\hspace{-.5pt}C\hspace{-.6pt}}(\mathfrak{F}^{\hspace{-.5pt}C\hspace{-.6pt}}(M^\ydt))}(j,t;i,s))_{(j,t;i,s)\in \mathbb Z^4}}&&&&& \oplus_{(j,t)\in \mathbb Z^2\,} \mathfrak{G}(\mathfrak{F}(M^j)^{t-j})^{n+1-t}
}$$
where
$$
d^n_{\mathfrak{G}^C(\mathfrak{F}^C(M^\pdt))}(j,t;i,s)=
%
%
\left\{\begin{array}{ll}
(-1)^{s} d_{\mathfrak G(\mathfrak F(M^i)^{s-i})}^{\hspace{.4pt}n-s}, & t=s, j=i;   \\ \vspace{-9pt} \\
(-1)^i \mathfrak G(d^{s-i}_{\mathfrak F(M^i)})^{n-s},             & t=s+1, j=i; \\ \vspace{-9pt} \\
\mathfrak G( \mathfrak F(d_M^i)^{s-i})^{n-s}, & t=s+1, j=i+1; \\ \vspace{-9pt} \\
0,                          & \text{otherwise}.
\end{array}\right.
$$
Setting $p=s-i$ and $q=t-j$, we see that $\mathfrak{G}^C(\mathfrak{F}^C(M^\ydt))$ is also described by
$$\xymatrixcolsep{16pt}\xymatrixrowsep{18pt}\xymatrix{
\mathfrak{G}^{\hspace{-.5pt}C\hspace{-.6pt}}(\mathfrak{F}^{\hspace{-.5pt}C\hspace{-.6pt}}(M^\ydt))^n \ar@{=}[d] \ar[rrrr]^{d^n_{\mathfrak{G}^{\hspace{-.5pt}C\hspace{-.6pt}}(\mathfrak{F}^{\hspace{-.5pt}C\hspace{-.6pt}}(M^\ydt))}} &&&& \mathfrak{G}^{\hspace{-.5pt}C\hspace{-.6pt}}(\mathfrak{F}^{\hspace{-.5pt}C\hspace{-.6pt}}(M^\ydt))^{n+1} \ar@{=}[d] \\
\oplus_{(i,p)\in \mathbb Z^2} \mathfrak{G}(\mathfrak{F}(M^i)^p)^{n-i-p}
\ar[rrrr]^-{\left(d^n(j,q;i,p)\right)_{(j,q;i,p)\in \mathbb Z^4}}&&&& \oplus_{(j,q)\in \mathbb Z^2} \mathfrak{G}(\mathfrak{F}(M^j)^q)^{n+1-j-q}
}$$
where
$$
\begin{array}{rcl}
d^n(j,q;i,p) & = & d^n_{\mathfrak{G}^C\hspace{-1.5pt}(\mathfrak{F}^C(M^\pdt))}(j,q + j;i,p + i) \vspace{6pt} \\
&= & \left\{ \begin{array}{ll}
(-1)^{p+i} d_{\mathfrak G(\mathfrak F(M^i)^p)}^{\hspace{.4pt}n-i-p}, &   q= p, j = i; \vspace{4pt} \\
(-1)^i \mathfrak G(d^{s-i}_{\mathfrak F(M^i)})^{n-i-p} ,             &  q=p + 1,   j= i; \vspace{4pt} \\
\mathfrak G( \mathfrak F(d_M^i)^p)^{n-i-p},                          &   q=p, j=i + 1; \vspace{4pt} \\
0,                                                                   &  \text{otherwise}.
\end{array}\right.
\vspace{2pt}
\end{array}$$

Thus, we conclude that $(\mathfrak{G}^C\circ \mathfrak{F})^C(M^\pdt)=(\mathfrak{G}^C\circ \mathfrak{F}^C)(M^\pdt)$, for all $M^\pdt\in C(\cA)$. Similarly, we may show that $(\mathfrak{G}^C\circ \mathfrak{F})^C(f^\ydt)=(\mathfrak{G}^C\circ \mathfrak{F}^C)(f^\ydt)$, for every morphism $f^\ydt: M^\ydt\to N^\ydt$ in $C(\cA)$.
The proof of the proposition is completed.

%

\medskip

To conclude this section, we shall study how to extend morphisms of functors.

\medskip

\begin{Lemma}\label{main-lemma1}

Let $\fF, \mathfrak G: \cA\to C(\cB)$ be functors, where $\cA, \cB$ are concrete abelian categories with countable direct sums. \vspace{1pt} A functorial morphism $\eta: \mathfrak F\to \mathfrak G$ induces functorial morphisms $\eta^{DC}: \fF^{DC}\to \mathfrak G^{DC}$ and $\eta^C: \fF^C\to \mathfrak G^C.$

\end{Lemma}

\noindent{\it Proof.} Let $\eta=(\eta_{_{\hspace{-.6pt}M}}^\ydt)_{M\in \mathcal A}: \mathfrak F\to \mathfrak G$ be a functorial morphism. Fix a complex $M^{\ydt}\in C(\cA)$. Given $i, j\in \Z$, since $\eta_{_{\hspace{-.6pt}M}}^\ydt: \fF(M)\to \mathfrak G(M)$ is natural in $M$, we obtain a commutative diagram

$$\xymatrixcolsep{14pt}\xymatrixrowsep{10pt}\xymatrix{
& \mathfrak G(M^i)^{j+1}\\
\fF(M^i)^{j+1} \ar[ur]^{\eta_{M^i}^{\hspace{.5pt}j+1}}&  \\
& \mathfrak G(M^i)^j \ar[rr]^{\mathfrak G(d_M^i)^j}\ar[uu]_{(-1)^id^j_{\mathfrak G(M^i)}}&& \mathfrak G(M^{i+1})^j. \\
\fF(M^i)^j \ar[uu]^{(-1)^id^j_{\fF(M^i)}}\ar[ur]^-{\eta_{M^i} ^{\hspace{.5pt}j} } \ar[rr]^-{\fF(d_M^i)^j} && \fF(M^{i+1})^j\ar[ur]_-{\eta_{M^{i+1}}^{\hspace{.5pt}j}}
}\vspace{5pt}$$
Setting $\eta^{i,j}_{_{\hspace{-1pt}M^\ydt}}=\eta_{_{\hspace{-1pt}M^i}}^{\hspace{.5pt}j},\vspace{1pt}$ we obtain a morphism $\eta_{_{\hspace{-.6pt}M^\ydt}}^{\ydt\hspace{.5pt}\ydt}: \fF^{DC}(M^\ydt)\to \mathfrak G^{DC}(M^\ydt)$ in $DC(\cB)$,
and a morphism \vspace{1pt} $\eta^\ydt_{_{\hspace{-1pt}M^\ydt}}=\T(\eta_{_{\hspace{-1pt}M^\ydt}}^{\ydt\hspace{.5pt}\ydt}): \fF^C(M^\ydt)\to \mathfrak G^C(M^\ydt)$ in $C(\cB)$. Consider a morphism $f^\cdt: M^\ydt\to N^\ydt$ in $C(\cA)$. Given  $i, j\in \Z$, we obtain a commutative diagram
$$\xymatrix{
\fF(M^i)^j \ar[r]^{\eta_{M^i}^j} \ar[d]_{\fF(f^i)^j} & \mathfrak G (M^i)^j \ar[d]^{\mathfrak G(f^i)^j}\\
\fF(N^i)^j \ar[r]^{\eta_{N^i}^j} & \mathfrak G (N^i)^j.
}$$
Hence, $\mathfrak G^{DC}(f^\ydt)\circ \eta_{_{\hspace{-.6pt}M^\ydt}}^{\ydt \hspace{.6pt}\ydt}=\eta_{_{\hspace{-.6pt}N^\ydt}}^{\ydt \hspace{.6pt}\ydt}\circ \fF^{DC}(f^\ydt).$\vspace{2pt} Applying $\T: DC(\cB)\to C(\cB)$ yields $\mathfrak G^C(f^\cdt)\circ \eta_{_{\hspace{-.6pt}M^\ydt}}^\ydt=\eta_{_{\hspace{-.6pt}N^\ydt}}^\ydt\circ \fF^C(f^\ydt).$ Thus, $\eta_{_{\hspace{-.6pt}M^\ydt}}^{\ydt \hspace{.6pt}\ydt}$ \vspace{1.5pt} and $\eta^\ydt_{_{\hspace{-.6pt}M^\ydt}}$ are natural in $M^\ydt$. That is, we have a functorial morphism $\eta^{DC}=(\eta_{_{\hspace{-1pt}M^\ydt}}^{\ydt\hspace{.5pt}\ydt})_{M^\ydt\,\in C(\cA)}: \fF^{DC}\to \mathfrak G^{DC}$ \vspace{1pt} and a functorial morphism $\eta^C=(\eta_{_{\hspace{-.8pt}M^\ydt}}^\ydt)_{M^\ydt \,\in C(\cA)}: \fF^C\to \mathfrak G^C$. The proof of the lemma is complete.

\medskip

In general, the extension of a functorial quasi-isomorphism is not necessarily a functorial quasi-isomorphism. The following statement gives some sufficient conditions for this to be the case.

\medskip

\begin{Lemma}\label{main-lemma2}

Let $\fF, \mathfrak G: \cA\to C(\cB)$ be functors, where $\cA, \cB$ are abelian categories with countable direct sums. Let $\eta: \mathfrak F\to \mathfrak G$ be a functorial morphism inducing functorial morphisms $\eta^{DC}\hspace{-3pt}: \hspace{-1pt} \mathfrak F^{DC}\to \mathfrak G^{DC}\hspace{-3pt}$ and $\eta^{\hspace{.5pt}C}: \mathfrak F^C\to \mathfrak G^C\hspace{-2pt}.$ Suppose that $\fF(M^\cdt)^\cdt$ and $\mathfrak G(M^\cdt)^\cdt$ are diagonally bounded-above for some $M^{\ydt}\in C(\cA)$. If $\eta^{\ydt}_{\hspace{-.8pt}_{M^i}}$ is a quasi-isomorphism for every $i\in \Z$, then $\eta^C_{_{\hspace{-.8pt}M^\ydt}}$ is a quasi-isomorphism.

\end{Lemma}

\noindent{\it Proof.} Write $\eta=(\eta_{_{\hspace{-.6pt}M}}^\ydt)_{M\in \mathcal A}$, \vspace{.5pt} where $\eta_{_M^\ydt}: \mathfrak F(M)^\ydt\to \mathfrak G(M)^\ydt$ is a morphism in $C(\cB)$. By definition, $\eta^{DC}=(\eta^{\cdt\,\cdt}_{_{\hspace{-1pt}M^\ydt}})_{M^\cdt\in C(\mathcal{A})}$, where $\eta_{_{\hspace{-.6pt}M^\ydt}}^{\ydt\hspace{.5pt}\ydt}: \fF(M^\cdt)^\cdt\to \mathfrak{G}(M^\cdt)^\cdt$ is defined so that $\eta^{i,j}_{_{\hspace{-1pt}M^\ydt}}=\eta_{\hspace{-1pt}_{M^i}}^{\hspace{.5pt}j}: \fF(M^i)^j\to \mathfrak G(M^i)^j.\vspace{1pt} $ Moreover, $\eta^C=(\eta_{_{\hspace{-.8pt}M^\ydt}}^\ydt)_{M^\ydt \,\in C(\cA)}$, where $\eta^\ydt_{_{\hspace{-1pt}M^\ydt}}=\T(\eta_{_{\hspace{-1pt}M^\ydt}}^{\ydt\hspace{.5pt}\ydt}): \T(\fF(M^\cdt)^\cdt)\to \T(\mathfrak{G}(M^\cdt)^\cdt)$.

\vspace{1pt}

Since $\fF(M^\ydt)^\ydt$ and $\mathfrak G(M^\ydt)^\ydt$ are diagonally bounded-above, the vertical cone $V_{\eta_{\hspace{-.6pt}M^\ydt}^{\ydt\hspace{.5pt}\ydt}}$ of $\eta_{_{\hspace{-.6pt}M^\ydt}}^{\ydt\hspace{.5pt}\ydt}$ \vspace{1pt} is diagonally bounded-above. Assume that $\eta^\ydt_{\hspace{-.5pt}_{M^i}}: \fF(M^i)^\ydt\to \mathfrak G(M^i)^\ydt$ is a quasi-isomorphism, for every $i\in \Z$. \vspace{1pt} Then $\eta^{i, \ydt}: \frak{t}^i(\fF(M^i)^\ydt\to \frak{t}^i(\mathfrak G(M^i)^\ydt)$, that is the morphism from the $i$-th column of $\fF(M^\ydt)^\ydt$ to the $i$-th column of $\mathfrak G(M^\ydt)^\ydt$, is a quasi-isomorphism, for every $i\in \Z$. By Lemma \ref{V-cone}(1), the columns of $V_{\eta_{\hspace{-.6pt}M^\ydt}^{\ydt\hspace{.5pt}\ydt}}$\vspace{1pt} are exact, and by Corollary \ref{AAL}, $\T(V_{\eta_{\hspace{-.6pt}M^\ydt}^{\ydt\hspace{.5pt}\ydt}})$ is acyclic. Since $C_{\hspace{.5pt}\eta^\ydt_{_{\hspace{-1pt}M^\ydt}}}= C_{\hspace{.5pt}\T(\eta_{M^\ydt}^{\ydt\hspace{.5pt}\ydt})} = \T(V_{\eta_{\hspace{-.6pt}M^\ydt}^{\ydt\hspace{.5pt}\ydt}})$; see (\ref{Cone}), we see that $C_{\hspace{.5pt}\eta^\ydt_{M^{\ydt}}}$ \vspace{.4pt} is acyclic. Thus, $\eta^C_{_{\hspace{-1pt}M^\ydt}}=\eta^\ydt_{_{\hspace{-1pt}M^\ydt}}$ is a quasi-isomorphism. The proof of the lemma is completed.

\smallskip

\section{Modules and representations}

\medskip

\noindent The objective of this section is to study modules over an algebra defined by a locally finite quiver with relations. Most of the results will be generalizations of some classical results for modules over a locally bounded category; see \cite{BoG,GSZ} or for representations of a strongly locally finite quiver; see \cite{BLP}. We shall first study projective $J$-covers in the most general case, and then investigate injective envelopes in the strongly locally finite dimensional case, and finally discuss projective $n$-resolutions in the graded case. In particular, we shall show that a graded algebra is quadratic if and only if every simple module admits a linear projective $2$-resolution. This generalizes a well known result that a Koszul algebra is quadratic; see \cite{BGS}.

\medskip

Throughout this section, we shall consider a $k$-algebra $\La=kQ/R,$ where $Q$ is a locally finite quiver and $R$ is a weakly admissible ideal in $kQ$. Let $M$ be a left $\La$-module. The module $M$ is called {\it unitary} if $M=\sum_{x\in Q_0} e_xM$ and an element $u\in M$ is called {\it normalized} if $u\in e_xM$ for some $x\in Q_0.$ Moreover, we shall say that $M$ is {\it finitely supported} if $e_xM=0$ for all but finitely many $x\in Q_0$. 
We shall denote by $\ModLa$ the category of all unitary left $\Lambda$-modules, and by $\ModbLa$ 
and ${\rm mod}^{\hspace{.4pt}b\hspace{-2.5pt}}\La$ the full subcategories of $\ModLa$ of finitely supported modules
and of finite dimensional modules, respectively.

\medskip

On the other hand, a {\it representation} $M$ of the bound quiver $(Q, R)$ consists of a family of $k$-vector spaces $M(x)$ with $x\in Q_0$ and a family of $k$-linear maps $M(\alpha): M(x)\to M(y)$ with $\alpha: x\to y\in Q_1$, such that $M(\rho)=0$ for every relation $\rho\in R(x,y)$ with $x, y\in Q_0$. Here, $M(\gamma)=\sum_i\, \lambda_i M(\alpha_{i, m_i}) \circ \cdots \circ M(\alpha_{i,1})$ for any $\gamma=\sum_i\, \lambda_i \alpha_{i, m_i}\cdots \alpha_{i,1}\in kQ(x,y)$ with $\lambda_i\in k$ and $\alpha_{ij}\in Q_1$. In particular, we may write $M(\bar{\gamma})=M(\gamma)$, where $\bar{\gamma\hspace{1.6pt}}\hspace{-1pt}=\gamma+R\in \La$. A morphism $f: M\to N$ of representations of $(Q, R)$ consists of a family of $k$-linear maps $f(x)$ with $x\in Q_0$ such that $f(y)\circ M(\alpha)= N(\alpha) \circ f(x),$ for every $\alpha: x\to y$ in $Q_1$. We shall denote by ${\rm Rep}(Q, R)$ the category of all representations of $(Q, R)$. It is well known that we may regard a module $M\in \Mod\La$ as a representation $M\in {\rm Rep}(Q, R)$ in such a way that $M(x)=e_xM$ for $x\in Q_0$, and $M(\alpha): M(x)\to M(y)$ with $\alpha\in Q_1(x, y)$ is given by the left multiplication by $\bar{\alpha}$. And we may regard a morphism $f: M\to N$ in $\ModLa$ as a morphism $(f(x))_{x\in Q_0}: M\to N$ in ${\rm Rep}(Q, R)$, where $f(x): M(x)\to N(x)$ is obtained by restricting $f$.

\medskip

Given $a\in Q_0$, we obtain a projective module $P_a=\La e_a$ in $\ModLa$ and a simple module $S_a=\La e_a/ Je_a$. Regarded as a representation, $P_a(x)=e_x\La e_a$ for $x\in Q_0$, and $P_a(\alpha): P_a(x)\to P_a(y)$ with $\alpha \in Q_1(x, y)$ is the left multiplication by $\bar{\alpha}$. We shall denote by ${\rm Proj}\,\La$ the full additive subcategory of ${\rm Mod}\La$ generated by the modules isomorphic to $P_x\otimes V$ with $x\in Q_0$ and $V\in {\rm Mod}\,k$, and by ${\rm proj}\,\La$ the one generated by the modules isomorphic to $P_x$ with $x\in Q_0$. The following easy statement is well known in the locally bounded case.

\medskip

\begin{Lemma}\label{proj-rad}

Let $\La=kQ/R$ be a strongly locally finite dimensional algebra.

\vspace{-2pt}

\begin{enumerate}[$(1)$]

\item If $a\in Q_0$, then $JP_a$ is the largest proper submodule of $P_a.$

\vspace{.6pt}

\item If $P\in {\rm Proj}\hspace{.5pt}\La$, then $JP$ is the Jacobson radical of $P$.

\end{enumerate}

\end{Lemma}

\noindent{\it Proof.} Let $M$ be a submodule of $P_a$ but not contained in $JP_a$, for some $a\in Q_0$. Then, $e_a-u\in M$ for some $u\in JP_a$. By Lemma \ref{Alg-dec}(2), $u^n=0$ for some $n\ge 2$. Since $ue_a=u$, we have
$(e_a+e_au+\cdots + e_au^{n-1})(e_a-u)=e_a\in M$. Thus, $M=P_a$. Statement (2) follows from Statement (1).
%
%
The proof of the lemma is completed.

\medskip

\noindent{\sc Remark.} By Lemma \ref{proj-rad}(1), $P_a$ is indecomposable in case $\La$ is strongly locally finite dimensional. However, this is not true even if $\La$ is locally finite dimensional, for instance, if $\La$ is given by a single loop $\alpha$ with a relation $\alpha^2=\alpha^3$.

\medskip

We shall study morphisms involving modules in ${\rm Proj}\hspace{.5pt}\La$. Given $\gamma\in kQ(x,y)$ with $x, y\in Q_0$, the left multiplication by $\bar{\gamma\hspace{1.6pt}}\hspace{-1pt}$ yields a $k$-linear map $P_a(\bar{\gamma\hspace{1.6pt}}): P_a(y)\to P_a(x)$ for every $a\in Q_0$, while the right multiplication by $\bar{\gamma\hspace{1.6pt}}$ yields a $\La$-linear morphism
$P[\bar{\gamma\hspace{1.6pt}}\hspace{-1pt}]: P_y\to P_x$.
The following statement is evident.

\medskip

\begin{Lemma}\label{Mor}

Let $\La=kQ/R$, where $Q$ is a locally finite quiver and $R$ is weakly admissible ideal. Consider $M\in {\rm Mod}\La$ and $V\in {\rm Mod}\,k$. Given $a, b\in Q_0$, we obtain

\begin{enumerate}[$(1)$]

\item a $k$-linear isomorphism $\mathcal{P}_{a,b}: e_b\La e_a \to {\rm Hom}_{\it\Lambda}(P_b, P_a): u \mapsto P[u];$

\vspace{.5pt}

\item a $k$-linear isomorphism ${\mathcal M}_a:  {\rm Hom}_{\it\Lambda}(P_a, M) \to e_aM : f\mapsto f(e_a);$

\vspace{.5pt}

\item a $k$-linear isomorphism $\psi_M: \Hom_{\it\Lambda}(P_a\otimes V, M)\rightarrow \Hom_k(V,e_aM);$

\vspace{1pt}

\item a $k$-linear map ${\mathcal M}_{a,b}: e_b \La  e_a \to \Hom_k(e_aM, e_bM): u \mapsto M(u)$, where $M(u)$ denotes the left multiplication by $u$.

\end{enumerate}

\end{Lemma}

\noindent{\it Proof.} Statement (1), (2) and (4) are evident. Observing that $P_a$ is a $\La$-$k$-bimodule, we deduce Statement (3) easily from the adjoint isomorphism and Statement (2). The proof of the lemma is completed.

\medskip

In case $\La$ is locally finite dimensional, the following statement describes all the morphisms in ${\rm Proj}\hspace{.4pt}\La$; compare \cite[(7.6)]{BaL}.

\medskip

\begin{Lemma}\label{rqz-pm}

Let $\La=kQ/R$ be a locally finite dimensional algebra. Given $a, b\in Q_0$ and $V, W\in {\rm Mod}\,k$, every $\La$-linear morphism $f: P_a\otimes V \to P_b\otimes W$ is uniquely written as $f=\sum \, P[u]\otimes f_u,$ where $u$ runs over a basis of $e_a\La e_b$, and $f_u\in {\rm Hom}_k(V, W).$

\end{Lemma}

\noindent{\it Proof.} Let $f: P_a\otimes V\to P_b\otimes W$ be $\La$-linear. Then, $f(e_a\otimes V)\subset e_a \La e_b \,\otimes W$. Let  $\{u_1, \ldots, u_n\}$ be a finite basis of $e_a\La e_b$. If $v\in V$, then $f(e_a\otimes v)=\textstyle{\sum}_{i=1}^n\, u_i\otimes w_i,$ for some unique $w_1, \ldots, w_n\in W$. This yields $k$-linear maps $f_i: V\to W: v\mapsto w_i$, for $i=1, \ldots, n$. We see easily that $f=\textstyle{\sum}_{i=1}^n P[u_i]\otimes
f_i,$ and this expression is unique. The proof of the lemma is completed.

\medskip

Given $M\in {\rm Mod}\La$, an epimorphism $d: P\to M$ with $P\in {\rm Proj}\hspace{.4pt}\La$ is called a {\it projective $J$-cover} of $M$ if ${\rm Ker}(d)\subseteq JP$. 
For instance, for $a\in Q_0$, the canonical projection $d_a: P_a\to S_a$ is a projective $J$-cover of $S_a$.
In general, we shall put $T(M)=M/JM$, called {\it $J$-top} of $M$. A generating set $\{u_1, \ldots, u_s\}$ of $M$ is called a {\it $J$-top basis} if 
$\{u_1+JM, \ldots, u_s+JM\}$ is a $k$-basis of $T(M)$, and such a $J$-top basis is {\it normalized} if $u_1, \ldots, u_s$ are normalized. The following statement is well known in the finite dimensional case; see \cite[(1.1)]{LiM}, whose proof is left to the reader.

\medskip

\begin{Lemma}\label{p-cover}

Let $\La=kQ/R$, where $Q$ is locally finite and $R$ is weakly admissible. A module $M\in \ModLa$ has a $J$-top basis $\{u_1, \ldots, u_s\}$ with $u_i\in e_{a\hspace{-.8pt}_i}M$  if and only if it has a projective $J$-cover $d: P_{a_1} \oplus \cdots \oplus P_{a\hspace{-.8pt}_s} \to M$ with $d(e_{a\hspace{-.8pt}_i})=u_i$, where $a_1, \ldots, a_s\in Q_0$.

\end{Lemma}

%

\medskip

Given a module $M\in {\rm Mod}\hspace{.5pt}\La$, we shall call an exact sequence

$$\xymatrix{P^{-n} \ar[r]^{d^{-n}} & P^{-n+1} \ar[r] & \cdots \ar[r] & P^{-1} \ar[r]^{d^{-1}} &  P^0 \ar[r]^{d^0} & M \ar[r] & 0}$$
a {\it $J$-minimal projective $n$-resolution} of $M$ over ${\rm proj}\,\La$ if $P_i\in {\rm proj}\hspace{.5pt}\La$ and $d^{-i}$ co-restricts to a projective $J$-cover of ${\rm Im}(d^{-i})$, for $i=0, \ldots, n$. The following statement is well known in case $Q$ is finite; compare \cite[(2.5)]{GMV}.

\medskip

\begin{Cor} \label{proj-pres}

Let $\La=kQ/R$, where $Q$ is a locally finite quiver and $R$ is weakly admissible ideal in $kQ$, and let $a\in Q_0$ with $Q_1(a, -)=\{\alpha_i: a \to b_i \mid i=1, \ldots, r\}$. Then $S_a$ admits a $J$-minimal projective $1$-resolution
$$\xymatrixcolsep{27pt}
\xymatrix{P_{b_1} \oplus \cdots \oplus P_{b_r} \ar[rr]^--{(P[\bar{\alpha}_1], \cdots, P[\bar{\alpha}_r]\hspace{.5pt})} && P_a \ar[r]^{d_a} & S_a \ar[r] & 0.}$$

\end{Cor}

\noindent{\it Proof.} Clearly, ${\rm Ker}(d_a)=JP_a$, which has a $J$-top basis $\{\bar{\alpha}_1, \cdots, \bar{\alpha}_r\}$. Considering the inclusion map $j: JP_a\to P_a$, by Lemma \ref{p-cover}, we obtain a projective $J$-cover $d: P_{b_1} \oplus \cdots \oplus P_{b_r} \to JP_a$ such that $(P[\bar{\alpha}_1], \cdots, P[\bar{\alpha}_r])=j\circ d$. The proof of the corollary is completed.

\medskip

Next, we shall define an exact functor $D: {\rm Mod}\La\to {\rm Mod}\La^{\rm o}$ in order to define injective-like modules. It is more convenient to identify unitary modules as representations. Given a module $M\in {\rm Mod}\La$, we define a module $DM\in {\rm Mod}\La^{\rm o}$ by $(DM)(x)=\Hom_k(M(x), k)$ for $x\in Q_0$, and $(DM)(\alpha^{\rm o})=\Hom_k(M(\alpha), k)$ for $\alpha\in Q_1$. Given a morphism $f: M\to N$, we define a morphism $Df: DN\to DM$ by $(Df)(x)=\Hom_k(f(x),k): (DN)(x)\to (DM)(x)$ for $x\in Q_0$.

\medskip

Given $a\in Q_0$, we shall write $P_a^{\rm o}=\La^{\rm o} e_a\in {\rm Mod}\hspace{.4pt}\La^{\rm o}$ and $I_a=DP^{\rm o}_a\in {\rm Mod}\La$. More explicitly, $I_a(x)=\Hom_k(e_x \La^{\rm o}e_a, k)$ for all $x\in Q_0$; and if $u\in e_y\La e_x$ and $f\in I_a(x)$, then $u f \in I_a(y)$ such that $(u f)(v^{\rm o})=f(u^{\rm o}v^{\rm o}),$ for all $v\in e_a\La e_y$. In general, $I_a$ is neither indecomposable nor injective in ${\rm Mod}\La$. By abuse of notation, we shall denote by ${\rm Inj}\hspace{.4pt}\La$ the full additive subcategory of ${\rm Mod}\La$ generated by the modules isomorphic to $I_x\otimes V$, where $x\in Q_0$ and $V\in {\rm Mod}\,k$, and by ${\rm inj}\hspace{.4pt}\La$ the one generated by the modules isomorphic to $I_x$ with $x\in Q_0$. In case $\La$ is locally finite dimensional, however, the following statement says that $I_a\otimes V$ is indeed injective in $\ModLa$ for any $k$-vector space $V$; compare \cite[(1.3)]{BLP}.

\medskip

\begin{Prop}\label{proj-inj}

Let $\La=kQ/R$ be a locally finite dimensional algebra. Consider $M\in {\rm Mod}\hspace{.5pt}\La$ and $V\in {\rm Mod}\hspace{.5pt}k$. Given $a\in Q_0$, we obtain a $k$-linear isomorphism
$$\phi_{\hspace{-1pt}_M}: \Hom_{\it\Lambda}(M, I_a\otimes V)\rightarrow \Hom_k(e_aM,V).$$

\end{Prop}

\noindent{\it Proof.} Fix $a\in Q_0$. For each $x\in Q_0$, since
$e_x\La^{\rm o} e_a$ is finite dimensional, we deduce from Corollary \ref{Cor 1.2} a $k$-linear isomorphism
$$\sigma_x: I_a(x)\otimes V = \Hom_k(e_x\La^{\rm o} e_a, k) \otimes V \to \Hom_k(e_x\La^{\rm o} e_a,V)$$ such that $\sigma_x(h\otimes v)(u^{\rm o})=h(u^{\rm o})v$, for $h\in I_a(x)$, $v\in V$ and $u\in e_a \La e_x$. Moreover, we have a $k$-linear map $\psi_a: \Hom_k(e_a\La^{\rm o} e_a, V)\rightarrow V: g\mapsto g(e_a).$ Observing that every $\La$-linear morphism $f: M\to I_a\otimes V$ consists of a family of $k$-linear maps $f_x: e_xM \to I_a(x)\otimes V$ with $x\in Q_0$, we obtain a $k$-linear map
$$\phi_{\hspace{-1pt}_M}:\Hom_{\it\Lambda}(M, I_a\otimes V)\rightarrow \Hom_k(e_aM, V): f\rightarrow \psi_a \circ \sigma_a\circ f_a.$$

Suppose that $\phi_{\hspace{-1pt}_M}(f)=0$. We claim that $f=0$. Let $x\in Q_0$. For any $m\in e_xM$, we write $f_x(m)=\sum_{i=1}^s h_i\otimes v_i$, where $h_i\in \Hom_k(e_x \La^{\rm o}e_a,  k)$ and $v_i\in V$. We may assume that the $v_i$ are $k$-linearly independent. Given any $u\in e_a\La e_x$, we obtain $u m \in e_a M$. Since $f$ is $\La$-linear, $f_a(um)=u f_x(m) = \sum_{i=1}^s (u h_i)\otimes v_i.$ Thus,
$$0=\phi_{\hspace{-1pt}_M}(f)(um)={\textstyle\sum}_{i=1}^s\sigma_a(uh_i\otimes v_i)(e_a)={\textstyle\sum}_{i=1}^s(uh_i)(e_a) v_i={\textstyle\sum}_{i=1}^sh_i(u^{\rm o})v_i.$$
Since the $v_i$ are assumed to be $k$-linearly independent, $h_i(u^{\rm o})=0$, and hence, $h_i=0$, for $i=1, \ldots, s$. Therefore, $f_x(m)=0$. This establishes our claim.

Conversely, consider a $k$-linear map $g_a: e_aM \rightarrow V$. Let $x\in Q_0$ and $m\in e_xM$. We obtain a $k$-linear map $g_x(m): e_x\La^{\rm o}e_a \to V: u^{\rm o} \rightarrow  g_a(u m),$ and hence, a $k$-linear map $f_x: e_xM \to I_a(x)\otimes V: m\mapsto \sigma_x^{-1}(g_x(m)).$ We claim that this yields a $\La$-linear morphism $f=(f_x)_{x\in Q_0}: M\to I_a\otimes V$. Indeed, given $w\in e_y\La e_x$ and $m\in e_xM$, we obtain $\sigma_y(f_y(w m))(u^{\rm o})=g_y(wm)(u^{\rm o})=g_a((uw)m)$. On the other hand, we may write \vspace{1pt} $g_x(m)=\sum_{i=1}^s \sigma_x (h_i\otimes v_i),$ for some $h_i\in I_a(x)$ and $v_i\in V$. Now, $w f_x(m)= w \sigma_x^{-1}(g_x(m)) = \sum_{i=1}^s (w h_i)\otimes v_i$. For any $u\in e_a\La e_y$, we see that
$$\textstyle\sigma_y(w f_x(m))(u^{\rm o}) = \sum_{i=1}^s h_i(w^{\rm o}u^{\rm o}) v_i = \sum_{i=1}^s \sigma_x(h_i\otimes v_i)((uw)^{\rm o})
                           = g_a((uw)m).$$
%
Thus, $\sigma_y(w f_x(m))=\sigma_y(f_y(w m))$, and hence, $f_y(w m))=w f_x(m).$ This establishes our second claim. Clearly, $\phi_{\hspace{-1pt}_M}(f)=g_a$. The proof of the proposition is completed.

\medskip

We shall now describe the morphisms in ${\rm Inj}\La$. Given $u\in e_b\La e_a$ with $a, b\in Q_0$, we obtain a $\La^{\rm o}$-linear morphism $P[u^{\rm o}]: P^{\hspace{.5pt}\rm o}_a\to P^{\hspace{.4pt}\rm o}_b$. Applying the exact functor $D: {\rm Mod}\,\La^{\rm o}\to {\rm Mod}\,\La$, we obtain a $\La$-linear morphism $I[u]=DP[u^{\rm o}]: I_b\to I_a$ such that
$I[u](f)(v^{\rm o})=f(v^{\rm o}u^{\rm o})$, for all $f\in I_a(x)$ and $v \in e_a\La e_x$ with $x\in Q_0$. 

\medskip

\begin{Lemma}\label{Inj-Mor-2}

Let $\La=kQ/R$ be a locally finite dimensional algebra. Given $a, b\in Q_0$ and $V, W\in {\rm Mod}\,k,$ every $\La$-linear morphism $f: I_a\otimes V \to I_b\otimes W$ is uniquely written as $f=\sum \, I[u]\otimes f_u,$ where $u$ runs over a $k$-basis of $e_a\La e_b$ and $f_u\in {\rm Hom}_k(V, W).$

\end{Lemma}

\noindent{\it Proof.} Fix $a, b\in Q_0$. Since $e_a\La e_b$ is finite dimensional, we obtain a $k$-linear isomorphism
$\theta_{a, b}: e_a\La e_b \to \Hom_k(\Hom_k(e_b\La^{\rm o} e_a, k), k): u\mapsto \theta_{a,b}(u),$ where
$\theta_{a,b}(u)$ sends $f\in I_a(u)$ to $f(u^{\rm o})$. Given $V, W\in {\rm Mod}\,k$, we have $k$-linear isomorphisms
$$\xymatrixcolsep{20pt}\xymatrixrowsep{16pt}\xymatrix{e_a\La e_b \otimes \Hom_k(V, W) \ar[rr]^-{\theta_{a, b} \otimes 1} && \Hom_k(I_a(b), k)\otimes \Hom_k(V, W) \ar[d]^-{\varphi} \\
\Hom_{\it\Lambda}(I_a\otimes V, I_b\otimes W)\ar[rr]^--\phi &&
{\Hom_k(I_a(b)\otimes V, W),\hspace{40pt}}}$$ where $\varphi$ is as defined in Lemma \ref{Tensor}, and $\phi$ is as defined in Proposition \ref{proj-inj}. We claim, for $u\in e_a\La e_b$ and $h\in \Hom_k(V, W)$, that $\phi(I[u]\otimes h)=(\varphi\circ (\theta_{a,b}\otimes \id)) (u\otimes h)$. Indeed, $\phi(I[u]\otimes h)$ is the composite of the maps in the sequence
$$\xymatrixcolsep{32pt}\xymatrix{I_a(b)\otimes V \ar[r]^-{I[u]\otimes h} & I_b(b)\otimes W \ar[r]^-{\sigma_b} & \Hom_k(e_b\La^{\rm o} e_b, W) \ar[r]^-{\psi_b} & W,}$$ where $\sigma_b$ and $\psi_b$ are as defined in the proof of Proposition \ref{proj-inj}. Now,
given any $g\in I_a(b)$ and $v\in V$, we obtain $\left(\varphi(\theta_{a,b}(u)\otimes h)\right)(g\otimes v)=\theta_{a,b}(u)(g)h(v)=g(u^{\rm o})h(v)$ and $\phi(I[u]\otimes h)(g\otimes v)=\sigma_b(I[u](g)\otimes h(v))(e_b)=I[u](g)(e_b)\, h(v)=g(u^{\rm o})\, h(v).$
This establishes our claim. As a consequence, we obtain a $k$-linear isomorphism
$$\phi^{-1}\circ \varphi \circ (\theta_{a,b}\otimes 1): e_a\La e_b \otimes \Hom_k(V, W) \to \Hom_{\it\Lambda}(I_a\otimes V, I_b\otimes W):
u\otimes h\to I[u]\otimes h.$$ The proof of the lemma is completed.

\medskip

Given a module $M\in \ModLa$, we shall write $S(M)=\{m\in M \mid J m=0\}$, called the {\it $J$-socle} of $M$.

\medskip

\begin{Lemma}\label{Special-func}

Let $\La=kQ/R,$ where $Q$ is a locally finite quiver and $R$ is a weakly admissible ideal. If $a\in Q_0$, then

\vspace{-1pt}

\begin{enumerate}[$(1)$]

\item $S(I_a)$ has a basis $\{e_a^\star\}$, where $e_{\hspace{-1pt}a}^\star\in I_a(a)$ with $e_{\hspace{-1pt}a}^\star(e_a)=1$ and $e_{\hspace{-1pt}a}^\star(e_a J^{\rm o} e_a)=0;$ 

\vspace{1.5pt}

\item  $S(I_a/S(I_a))$ has a basis $\{\alpha^\star+S(I_a) \mid \alpha: x\to a \in Q_1(-, a)\},$ where $\alpha^\star\in I_a(x)$ such that
$\alpha^\star(\bar\alpha^{\rm o})=1$, and $\alpha^\star(\bar{\gamma\hspace{1.6pt}}\hspace{-1pt}^{\rm o})=0$ for all $\gamma\in Q(x, a)$ with $\gamma \ne \alpha.$

\end{enumerate}\end{Lemma}

\noindent{\it Proof.} Fix $a\in Q_0$. 
Clearly, $e_{\hspace{-1pt}a}^\star\in S(I_a)$. If $f\in I_a(x)$ for some $x\in Q_0$, which is neither zero nor a multiple of $e_a^\star$, then $f(u^{\rm o})\ne 0$ for some $u\in e_a J e_x$, that is, $(u\cdot f)(e_a^\star)\ne 0$. Hence, $f\notin S(I_a)$. Thus, $S(I_a)=k\hspace{.5pt}e_a^\star$.

Let $\alpha \in Q_1(x, a)$. The existence of $\alpha^\star$ follows from Lemma \ref{Alg-dec}(1). It is evident that $\bar \alpha \cdot \alpha^\star=e_a^\star$. Consider $\beta\in Q_1(x, y)$ with $\beta \ne \alpha$. For $\delta \in Q(y, a)$, since $\delta \beta \ne \alpha$, we obtain $(\bar\beta \cdot \alpha^\star)(\bar\delta ^{\rm o})=\alpha^\star(\bar\beta^{\rm o}\bar\delta ^{\rm o})=0$. Thus, $\alpha^\star + S(I_a)\in S(I_a/S(I_a))$. Let $\alpha_1, \ldots, \alpha_r\in Q_1(x, a)$ such that $\sum_{i=1}^r\lambda_i \alpha^\star_i\in S(I_a)$, where $\lambda_1$, $\ldots,$ $\lambda_r\in k$ are non-zero. Since $(\sum_{i=1}^r\lambda_i \alpha^\star_i)(\alpha_1)=\lambda_1$, by Statement (1), $\sum_{i=1}^r\lambda_i \alpha^\star_i=\lambda \hspace{.5pt} e_a^\star$, where $\lambda\in k$ is non-zero. Then, $x=a$ and $e_a^\star(\alpha_1)=\lambda^{-1}\lambda_1\ne 0$, a contradiction. Thus the classes $\alpha^\star + S(I_a)$, with $\alpha\in Q_1(-, a)$, are $k$-linearly independent in $S(I_a/S(I_a))$.

Consider $g+S(I_a)\in S(I_a/S(I_a)),$ where $g\in I_a(x)$ for some $x\in Q_0$. Let $\rho\in Q(x, a)$ be of length at least two. Write $\rho=\delta \alpha$, where $\alpha: x\to y$ is an arrow and $\delta: y\to a$ is a non-trivial path. Since $\bar\alpha g \in S(I_a)$ and $\delta: y\to a$ is  non-trivial, we obtain $g(\bar\rho^{\rm o})=(\bar\alpha g)(\bar\delta^{\rm o})= 0$. Hence $g(e_x(J^{\rm o})^2 e_a)=0$. By Lemma \ref{Alg-dec}(1), $g=\sum_{\gamma\in Q_{\le 1}(x, a)} g(\bar{\gamma\hspace{1.6pt}}\hspace{-1pt}^{\rm o}) \gamma^\star,\vspace{1pt}$ and hence, $g+S(I_a)=\sum_{\alpha\in Q_1(x, a)} g(\bar\alpha^{\rm o}) (\alpha^\star+S(I_a))$. The proof of the lemma is completed.

\medskip

The following statement is well known in the locally bounded case.

\smallskip

\begin{Cor}\label{inj-soc}

Let $\La=kQ/R$ be strongly locally finite dimensional. If $a\in Q_0$, then $S(I_a)$ and
$S(I_a/S(I_a))$ are essential socles of $I_a$ and $I_a/S(I_a)$, respectively.

\end{Cor}

\noindent{\it Proof.} By Proposition \ref{Alg-dec}(2), $J$ is the Jacobson radical of $\La$. Thus, the $J$-socle of a module is the socle. Let $h\in I_a(x)\backslash S(I_a)$, for some $x\in Q_0$. Then, $h(e_x J^{\rm o} e_a)\ne 0$. Since $J^{\rm o}$ is locally nilpotent, there exists a maximal positive integer $s$ such that $h(e_x(J^{\rm o})^s e_a)\ne 0$. Then, $h(\bar\zeta^{\,\rm o})=\lambda \ne 0$, for some $\zeta\in Q_s(x, a)$.

First, $\bar\zeta  h\in I_a(a)$ with $(\bar\zeta  h)(e_a)=h(\bar\zeta^{\rm o})=\lambda$. By the maximality of $s$, we see that $(\bar\zeta h)(e_a J^{\rm o} e_a)=0$, and hence, $\bar\zeta h=\lambda \hspace{.5pt} e_a^\star$. Thus, $S(I_a)$ is essential in $I_a$.

Next, write $\zeta=\beta \xi$, where $\beta\in Q_1(b, a)$ and $\xi\in Q_{s-1}(x, b)$ for some $b\in Q_0$. Then, $\bar\xi f \in I_a(b)$ with $(\bar\xi h)(\bar\beta^{\rm o})=h(\bar\zeta^{\,\rm o})\ne 0$. In particular, $\bar\xi (h+ S(I_a))\ne 0$. By the maximality of $s$, however, $(\bar\xi h)(e_b (J^{\rm o})^2 e_a)=0$.  In view of Lemma \ref{Alg-dec}(1), we see that $\bar\xi (h+ S(I_a))=\bar\xi h + S(I_a)=\sum_{\alpha\in Q_1(b, a)} (\bar\xi h)(\bar\alpha^{\rm o})\cdot (\alpha^\star+S(I_a))$. \vspace{1pt} This shows that $S(I_a/S(I_a))$ is essential in $I_a/S(I_a)$. The proof of the corollary is completed.

\medskip

A set $\{u_1, \ldots, u_s\}$ of normalized elements of a module $M \in {\rm Mod}\,\La$ is called an {\it essential socle basis} of $M$ if $\{u_1, \ldots, u_s\}$ is a basis $S(M)$, while $S(M)$ is essential in $M$. The following result is well known in case $\La$ is finite dimensional, and its proof is left to the reader.

\medskip

\begin{Lemma}\label{i-envelope}

Let $\La=kQ/R$ be a strongly locally finite dimensional algebra. A module $M\in {\rm Mod}\hspace{.5pt}\La$ has an essential socle basis $\{u_1, \ldots, u_s\}$ with $u_i\in e_{\hspace{-.6pt}a_{\hspace{-.5pt}i}}M$ if and only if $M$ has an injective envelope $j: M\to I_{a_1} \oplus \cdots \oplus I_{a_s}$ with $j(u_i)=e_{\hspace{-.6pt}a_i}^\star$, where $a_1, \ldots, a_s\in Q_0$. 

\end{Lemma}

%

\medskip

The following statement is well known in case $Q$ is a finite.

\medskip

\begin{Cor}\label{inj-pres}

Let $\La=kQ/R$ be a strongly locally finite dimensional algebra. If $a\in Q_0$ with $Q_1(-, a)=\{\beta_i: b_i \to a \mid i=1, \ldots, s\}$, then \vspace{-2pt}
$$\xymatrixcolsep{25pt}\xymatrix{0 \ar[r] & S_a \ar[r]^{j_a} &
I_a \ar[rr]^--{\left(I[\bar{\beta}_1], \ldots, I[\bar{\beta}_s]\right)^t} && I_{b_1} \oplus \cdots \oplus I_{b_s},}
$$ is a mi\-nimal injective co-presentation of $S_a$, where $j_a$ sends $e_a+Je_a$ to $e_a^\star$.

\end{Cor}

\noindent{\it Proof.} Let $a\in Q_0$ with $Q_1(-, a)=\{\beta_i: b_i \to a \mid i=1, \ldots, s\}$. By Corollary \ref{inj-soc} and Lemma \ref{i-envelope}, $j_a$ is an injective envelope of $S_a$ with ${\rm Im}(j_a)=S(I_a)$. By Lemma \ref{Special-func} and Corollary \ref{inj-soc}, $I_a$ has an essential socle basis $\{\beta_1^\star+S(I_a), \ldots, \beta_s^\star + S(I_a)\}$. By Lemma \ref{i-envelope}, \vspace{1pt} we obtain an injective envelope $j: I_a/S(I_a)\to I_{b_1} \oplus \cdots \oplus I_{b_s},$ sending $\beta_i^\star+S(I_a)$ to $e^\star_{b_i}$, for $i=1, \ldots, s$. It is easy to see that $I[\bar\beta_i](\beta_i^\star)=e_{b_i}^\star$. As a consequence, $\left(I[\bar{\beta}_1], \ldots, I[\bar{\beta}_s]\right)^t$ is the composite of the canonical projection $I_a\to I_a/S(I_a)$ and the injective envelope $j$. The proof of the corollary is completed.

\medskip

The following statement will be needed later.

\medskip

\begin{Lemma}\label{ess-socle}

Let $\La=kQ/R$ be a strong finite dimensional algebra. Suppose that $M\in {\rm Mod}\La$ has a finitely supported essential socle. If $N$ is a submodule of $M$, then $M/N$ has an essential socle.

\end{Lemma}

\noindent{\it Proof.} By Proposition \ref{Alg-dec}(2), $J$ is the Jacobson radical of $\La$. Suppose that $S(M)$ is supported by $a_1, \ldots, a_r\in Q_0$. Let $N$ be a submodule of $M$ containing $S(M)$. Consider a non-zero element $w+N \in M/N$, where $w\in e_{b_1}M+\cdots + e_{b_s}M$, for some $b_1, \ldots, b_s\in Q_0.$ Since $J$ is locally nilpotent, there exists some $t>0$ such that $e_{a_j}J^te_{b_i}=0$ for all $1\le i\le s$ and $1\le j\le r$. Suppose that $v(w+N)\ne 0$ for some $v\in J^t$. Since $S(M)$ is essential in $M$, there exists some $u\in \La$ such that $0\ne (uv)w\in S(M)$. In particular, $e_{b_j}(uv)e_{a_i}\ne 0$ for some $1\le i\le r$ and $1\le j\le s$, which is a contradiction. Thus, $J^t(w+N)=0$. As a consequence, there exists some maximal $n$ with $0\le n<t$ for which
$J^n(w+N)\ne 0$. That is, $0\ne J^n(w+N)\subseteq S(M/N)$. The proof of the lemma is completed.

\medskip

For the rest of this section, we shall assume that $\La=kQ/R$, where $R$ is a homogeneous ideal in $kQ$. Then, $\La$ is positively graded with  $\La=\oplus_{i\ge 0}\La_i$, called {\it $J$-radical graduation}, where $\La_i$ is the $k$-vector subspace of $\La$ generated by $\bar{\gamma\hspace{1.6pt}}$ with $\gamma\in Q_i$. One says that $\La$ is {\it quadratic} if $R$ is a quadratic ideal. Let $M=\oplus_{i\in \Z}M_i$ and $N=\oplus_{i\in \Z}N_i$ be graded modules in ${\rm Mod}\,\La$. One says that $M$ is {\it generated} in degree $n$ if $M=\La M_n$. Given $a\in Q_0$, we see that $P_a$ and $S_a$ are graded modules generated in degree $0$, but $I_a$ is not necessarily graded. A $\La$-linear morphism $f: M\to N$ is called {\it homogeneous of degree} $n$ if $f(M_i)\subseteq N_{i+n}$ for all $i\in \Z$; and in this case, we write $f=(f_i)_{i\in \Z}$, where $f_i: M_i\to N_{i+n}$ is obtained by restricting $f$. Observe that $f$ is {\it graded} if it is homogeneous of degree $0$. The following statement is evident.

\medskip

\begin{Lemma}\label{grad-exact}

Let $\La=kQ/R$ be a grade algebra. If $L, M, N\in {\rm Mod}\,\La$ are graded, then a sequence
$\xymatrixrowsep{15pt}\xymatrixcolsep{20pt}\xymatrix{L\ar[r]^f & M\ar[r]^g & N}$ of homogeneous morphisms of degree $n$ is exact if and only if
$\xymatrixrowsep{15pt}\xymatrixcolsep{22pt}\xymatrix{L_{i-n}\ar[r]^{f_{i-n}} & M_i\ar[r]^{g_{i+n}} & N_{i+n}}$ is exact, for every $i\in \Z$.

\end{Lemma}

\medskip

The following statement is well known in case $Q$ is finite; see, for example, \cite[(2.4)]{GMV}, and its proof is left to the reader.

\medskip

\begin{Lemma}\label{gr-proj-c}

Let $\La=kQ/R$, where $Q$ is locally finite and $R$ is homogeneous. Let $M\in{\rm Mod}\La$ be graded and finitely generated with homogeneous projective $J$-covers $f: P\to M$ and $f': P'\to M$. Then $f'=f\circ g$, for some graded isomorphism $g$.

\end{Lemma}


\medskip

The following result describes a $J$-minimal projective $2$-resolution of a simple module in the graded case; compare \cite[(2.4)]{GSZ}.

\medskip

\begin{Lemma}\label{ex-p-pres}

Let $\La=kQ/R$, where $Q$ is locally finite and $R$ is homogeneous with a minimal generating set $\Oa$. Let $a\in Q_0$ with $Q_1(a, -)=\{\alpha_i: a \to b_i \mid i=1, \ldots, r\}$ and $\Oa(a, -)=\{\rho_1, \ldots, \rho_s\}$. If $\rho_j=\sum_{j=1}^r \gamma_{ij}\alpha_i$\vspace{1pt} with $\gamma_{ij}\in kQ(b_i, c_j)$, then $S_a$ has $J$-minimal projective $2$-resolution
$$\xymatrixcolsep{25pt}\xymatrix{P_{c_1}\oplus \cdots \oplus P_{c_s} \ar[rr]^{(P[\bar{\gamma\hspace{1pt}}\hspace{-1pt}_{ij}])_{r\times s}} && P_{b_1} \oplus \cdots \oplus P_{b_r} \ar[rr]^-{\left(P[\bar{\alpha}_1], \cdots, P[\bar{\alpha}_r]\right)} && P_a \ar[r]^{d_a} & S_a \ar[r] & 0.}$$

\end{Lemma}

\noindent{\it Proof.} Let $\rho_j=\sum_{j=1}^r \gamma_{ij}\alpha_i,\vspace{1pt}$ where $\gamma_{ij}\in kQ(b_i, c_j)$.  Write $d_1=\left(P[\bar{\alpha}_1], \cdots, P[\bar{\alpha}_r]\right)$ and $d_2=(P[\bar{\gamma\hspace{1pt}}_{\hspace{-1pt}ij}])_{r\times s}$.
By Corollary \ref{proj-pres}, it suffices to show that $d_2$ co-restricts to a projective $J$-cover of ${\rm Ker}(d_1)$. Since $u_j=(\bar{\gamma}_{1j}, \ldots, \bar{\gamma}_{rj}) \in {\rm Ker}(d_1)$, by Lemma \ref{p-cover}, it amounts to show that $\{u_1, \ldots, u_s\}$ is a $J$-top basis of ${\rm Ker}(d_1)$.

Let $v=(\bar{\delta}_1, \ldots, \bar{\delta}_r)\in {\rm Ker}(d_1),$ where $\delta_i\in kQ(b_i, -)$. We may assume that $\delta_i\in kQ(b_i, c)$, for some $c\in Q_0$. Since $d_1(v)=0$, we obtain $\sum_{i=1}^r \delta_i \alpha_i \in R(a, c)$, and hence,
$\textstyle{\sum}_{i=1}^r \delta_i \alpha_i={\sum}_{j=1}^s \omega_j \rho_j + {\sum}_{i=1}^r \eta_i \alpha_i ={\sum}_{i=1}^r ({\sum}_{j=1}^s \omega_j\gamma_{ij} + \eta_i)\alpha_i,$ where $\omega_j\in kQ(c_j, c)$ and $\eta_i\in R(b_i, c)$. This yields $\delta_i = {\sum}_{j=1}^s \omega_j\gamma_{ij} + \eta_i$, and hence,
$\bar\delta_i = {\sum}_{j=1}^s \bar\omega_j\bar{\gamma\hspace{1.6pt}}\hspace{-1pt}_{ij}$, for $i=1, \ldots, r$. As a consequence, $v={\sum}_{j=1}^s \bar\omega_j u_j.$

\vspace{1pt}

Assume next that $\sum_{j=1}^s \lambda_j \hspace{.5pt}u_j \in J \hspace{.5pt}{\rm Ker}(d_1)$, where $\lambda_j \in k$. As shown above, $\sum_{j=1}^s\lambda_ju_j = \sum_{j=1}^s \bar\nu_j u_j$, with $\nu_j \in kQ^+$. \vspace{1pt} Thus, $\sum_{j=1}^s \lambda_j \gamma_{ij}= \sum_{j=1}^s (\nu_j \gamma_{ij} +\eta_{ij}),$\vspace{.5pt} where $\eta_{ij}\in R(b_i, c_j)$, for $i=1, \ldots, r$. Calculating $\sum_{i=1}^r\sum_{j=1}^s \lambda_j \gamma_{ij} \alpha_i$, we obtain \vspace{.5pt} $\sum_{j=1}^s \lambda_j \rho_j = \sum_{j=1}^s \nu_j(\rho_j + \zeta_j),$ where $\zeta_j\in R(a, c_j)$. By Lemma \ref{min-rel}, $\lambda_j=0,$ for $i=1, \ldots, s.$ The proof of the lemma is completed.

\medskip

Let $M$ be a graded module in ${\rm Mod}\hspace{.5pt}\La$. A $J$-minimal projective $n$-resolution of $M$ over ${\rm proj}\,\La$ is called a {\it linear projective $n$-resolution} if the morphisms between the projective modules are homogenous of degree one.
The following statement extends a well known result that a Koszul algebra is quadratic; see \cite[(2.3.3)]{BGS}.

\medskip

\begin{Theo}\label{qua-alg}

Let $\La=kQ/R$, where $Q$ is locally finite and $R$ is homogeneous. Then, $\La$ is quadratic if and only if every simple $\La$-module admits a linear projective $2$-resolution over ${\rm proj}\hspace{.5pt}\La$.

\end{Theo}

\noindent{\it Proof.} Let $\Oa$ be a minimal generating set of $R$. Fix $a\in Q_0$. Since $\Oa(a, -)$ contains only finitely many quadratic relations, the necessity follows immediately from Lemma \ref{ex-p-pres}. Assume that $S_a$ admits a linear projective $2$-resolution over ${\rm proj}\,\La$. Letting $Q_1(a, -)=\{\alpha_i: a\to b_i \mid i=1, \ldots, r\}$, we deduce from Lemmas \ref{Inj-Mor-2}, \ref{gr-proj-c} and \ref{ex-p-pres} a commutative diagram with exact rows
$$\xymatrixcolsep{25pt}\xymatrix{P_2 \ar@{=}[d]\ar[rr]^-{d_2} && P_1 \ar[rr]^-{d_1} \ar[d]^{f_1} && P_a \ar[r]^{d_0} \ar[d]^{f_0} & S_a \ar[r] \ar@{=}[d] & 0\\
P_{c_1}\oplus \cdots \oplus P_{c_s} \ar[rr]^-{\left(P[\bar{\gamma}_{ij}]\right)_{r\times s}} && P_{b_1} \oplus \cdots \oplus P_{b_r} \ar[rr]^-{(P[\bar{\alpha}_1],\cdots, P[\bar{\alpha}_r])} && P_a \ar[r]^{d_a} & S_a \ar[r]  & 0,}$$
where the upper row is a linear projective $2$-resolution, $f_0, f_1$ are graded isomorphisms, and $\gamma_{ij}\in kQ(b_i, c_j)$. Since $f_1 \circ d_2$ is homogeneous of degree one, $\gamma_{ij}\in kQ_1(b_i, c_j)$ and $\eta_j=\sum_{i=1}^r\gamma_{ij}\alpha_i\in R_2(a, c_j),$ for $j=1, \ldots, s.$ By Lemma \ref{p-cover}, $\{u_j=(\bar{\gamma}_{1j}, \ldots, \bar{\gamma}_{rj}) \mid j=1, \ldots, s\}$ is a $J$-top basis of ${\rm Ker}(P[\bar{\alpha}_1],\cdots, P[\bar{\alpha}_r])$.

\vspace{1pt}

Let $\rho\in \Oa(a, c)$ be a relation of degree $n>2$. Write $\rho=\sum_{i=1}^r\, \gamma_i \alpha_i$, for some $\gamma_i\in kQ_{n-1}(b_i,c)$. Since $(\bar{\gamma}_1, \ldots, \bar{\gamma}_r)\in {\rm Ker}(P[\bar{\alpha}_1],\cdots, P[\bar{\alpha}_r])$, we see that $(\bar{\gamma}_1, \ldots, \bar{\gamma}_r)=\sum_{j=1}^s \bar\delta_{\hspace{-1pt}j} \hspace{1pt} u_j$, \vspace{1pt} for some $\delta_j\in kQ_{n-2}(c_j, c).$ Then, $\gamma_i=\sigma_i+\sum_{j=1}^s \delta_j \gamma_{ij}$, where $\sigma_i\in R(b_i, c),$ for $i=1, \ldots, r$. This yields
$\rho= 
\sum_{i=1}^r \sigma_i\alpha_i + \sum_{j=1}^s \delta_j \eta_j.$
Since $n>2$, we see that $\rho \in R(kQ^+)+(kQ^+)R,$ a contradiction to Lemma \ref{min-rel}. The proof of the theorem is completed.

\medskip

Let $M$ be a graded module in ${\rm Mod}\,\La$. Consider an exact sequence
$$\xymatrix{\cdots \ar[r] & P^{-n} \ar[r]^{d^{-n}} & P^{-n+1} \ar[r] & \cdots \ar[r] & P^{-1} \ar[r]^{d^{-1}} &  P^0 \ar[r]^{d^0} & M \ar[r] & 0,}$$
where $d^{-n}$ co-restricts to a projective $J$-cover of ${\rm Im}(d^{-n})$ for every $n\ge 0$. If all $d^{-n}$ with $n>0$ are homogeneous of degree one, then the complex
$$\xymatrix{\cdots \ar[r] & P^{-n} \ar[r]^{d^{-n}} & P^{-n+1} \ar[r] & \cdots \ar[r] & P^{-1} \ar[r]^{d^{-1}} &  P^0 \ar[r] & 0\ar[r] & \cdots }$$
is called a {\it linear projective resolution} of $M$ over ${\rm proj}\,\La$.

\medskip

\begin{Defn}\label{Kos-alg}

Let $\La=kQ/R$, where $Q$ is a locally finite quiver and $R$ is a homogeneous ideal in $kQ$. We shall say that $\La$ is {\it Koszul} if every simple $\La$-module admits a linear projective resolution over ${\rm proj}\hspace{.5pt}\La$.

\end{Defn}

\medskip

\noindent{\sc Remark.} By Theorem \ref{qua-alg}, a Koszul algebra is quadratic; compare \cite[(2.3.3)]{BGS}.

\smallskip

\section{Koszul complexes and Koszul algebras}

\medskip

\noindent The objective of this section is to present a self-contained combinatorial account of the Koszul theory for quadratic algebras defined by a locally finite quiver and for all (not only graded) modules. Although our main results will be similar to those in \cite{BGS}, we shall take an elementary approach with a local viewpoint and provide more detailed arguments. Indeed, localizing the Koszul complex stated in \cite{BGS} yields a local Koszul complex associated with each simple module $S$, which contains always a projective 2-resolution of $S$, and is a projective resolution if and only if $S$ has a linear projective resolution. Using an alternative interpretation of a local Koszul complex, we shall show that a quadratic algebra is Koszul if and only if its quadratic dual is Koszul. In the locally finite dimensional case, we shall prove that a quadratic algebra is Koszul if and only if every simple module has a linear injective co-resolution, or equivalently, the opposite algebra is Koszul.

\medskip

Throughout this section, we put $\La=kQ/R,$ where $Q$ is a locally finite quiver and $R$ is a quadratic ideal in $kQ$. We start with some notation. Let $a, x\in Q_0$ and $n\ge 0$. Recall that $R_n(a, x)=R(a, x)\cap kQ_n(a, x)$. For $n=0, 1$, we set $R^{(n)}(a, x)=kQ_n(a,x)$. For $n \ge 2$, let $R^{(n)}(a,x)$ be the subspace of $R_n(a,x)$ of the elements $\gamma$ which, for any $0\le j \le n-2$, can be written as $\gamma=\sum_i \hspace{.4pt} \zeta_i \rho_i \delta_i,$ for some $\zeta_i\in kQ_{n-2-j}(c_i, x);$ $\rho_i\in R_2(b_i, c_i);$ $\delta_i\in kQ_j(a, b_i)$, where $b_i, c_i\in Q_0$. In particular, $R^{(2)}(a,x)=R_2(a, x)$. Furthermore, set $R^{(n)}(a, -)=\oplus_{x\in Q_0}\, R^{(n)}(a, x)$. Given an arrow $\alpha$ in $Q$, there exists a unique $k$-linear map $\partial_\alpha: kQ\to kQ$ such, for every path $\rho$, that $\partial_\alpha(\rho)=\delta$ if $\rho=\alpha \delta$; and $\partial_\alpha(\rho)=0$ if $\alpha$ is not a terminal arrow of $\rho$. In particular, $\partial_\alpha$ vanishes on all trivial paths in $Q$, while $\partial_\alpha(\alpha)=\varepsilon_{\hspace{-1.7pt}s\hspace{.4pt}(\hspace{-1.3pt}\alpha)}$.

%
%
%
%
%
%
%
%
%
%
%

\medskip

\begin{Lemma}\label{diff}

Let $\La=kQ/R$, where $Q$ is a locally finite quiver and $R$ is a quadratic ideal. Given $a, x, y\in Q_0$ and $n>0$, we obtain a $\La$-linear morphism
$$\partial_a^{-n}(y,x)={\textstyle\sum}_{\alpha\in Q_1(y,x)}P[\bar{\alpha}]\otimes \partial_\alpha: P_x \otimes R^{(n)}(a, x)\to P_y \otimes R^{(n-1)}(a, y)$$
such that $\partial_a^{-n}(y,x)(u\otimes \zeta\hspace{.5pt}\delta)=u \bar{\zeta} \otimes \delta$,
for all $u\in P_x$ and $\zeta\hspace{.5pt}\delta \in R^{(n)}(a, x)$ with $\zeta\in kQ_1(y, x)$ and $\delta\in kQ_{n-1}(a, y)$.

\end{Lemma}

\noindent{\it Proof.} Fix $a, x, y\in Q_0$ and $n>0$. We first claim, for any $\alpha\in Q_1(y, x)$, that $\partial_\alpha$ maps $R^{(n)}(a, x)$ into $R^{(n-1)}(a, y).$ Since $\partial_\alpha$ maps $kQ_n(a, x)$ into $kQ_{n-1}(a, y)$, we may assume that $n\ge 3$. Let $\gamma\in R^{(n)}(a, x)$ and $0\le j\le (n-1)-2.$ By definition, $\gamma={\sum}_{i=1}^r \alpha_i \zeta_i \rho_i \delta_i,$ where $\alpha_i\in Q_1,$ $\zeta_i\in kQ_{n-3-j},$ $\rho_i\in R_2,$ and $\delta_i\in kQ_j.$ We may assume that there exists some $1\le s\le r$ such that $\alpha_i=\alpha$ if and only if $1\le i\le s$. Then, $\partial_\alpha(\gamma)=\sum_{i=1}^s \zeta_i \rho_i \delta_i,$ where $\zeta_i\in kQ_{n-3-j};$ $\rho_i\in R_2;$ $\delta_i\in kQ_j.$ By definition, $\partial_\alpha(\gamma)\in R^{(n-1)}(a, y).$ This establishes our claim. As a consequence, we obtain a morphism
$\partial_a^{-n}(y,x)$ as stated in the lemma. Let $u\in P_x$ and $\zeta\hspace{.5pt}\delta \in R^{(n)}(a, x)$ with $\zeta\in kQ_1(y, x)$ and $\delta\in kQ_{n-1}(a, y)$. Writing $\zeta=\sum_{\beta\in Q_1(y, x)} \lambda_\beta \beta$ with $\lambda_\beta\in k$ yields
$$\partial_a^{-n}(x,y)(u\otimes \zeta\delta)={\textstyle\sum}_{\alpha, \beta\in Q_1(y,x)} P[\bar{\alpha}](u)\otimes \lambda_\beta \partial_\alpha( \beta \delta) = u \bar{\zeta} \otimes \delta.$$ The proof of the lemma is completed.

\medskip

Given $a\in Q_0$, in view of Lemma \ref{diff}, we may define a sequence $K_a^\cdt$ as follows:
$$\xymatrixrowsep{16pt}\xymatrix{\cdots \ar[r] & K_a^{-n} \ar[r]^-{\partial_a^{-n}} & K_a^{-n+1} \ar[r] &
\cdots \ar[r] & K_a^{-1} \ar[r]^{\partial_a^{-1}} & K_a^0  \ar[r] & 0 \ar[r] & \cdots,}$$
where $K_a^{-n}=\oplus_{x\in Q_0} P_x \otimes R^{(n)}(a, x)$ for $n\ge 0$ and $\partial_a^{-n}=(\partial_a^{-n}(y,x))_{(y,x)\in Q_0\times Q_0}\vspace{1.5pt}$, for $n> 0$. By Lemma \ref{diff}, $\partial_a^{-n}$ is homogeneous of degree one. Since $K_a^0=P_a\otimes k\varepsilon_a$, we have an additional $\La$-linear morphism $\partial_a^0: K_a^0 \to S_a: e_a\otimes \varepsilon_a\mapsto e_a+JP_a$.

\medskip

\begin{Lemma}\label{K-cplx}

Let $\La=kQ/R$, where $Q$ is a locally finite quiver and $R$ is a quadratic ideal. Given $a\in Q_0$, the sequence $K_a^\cdt$ is a complex over ${\rm proj}\hspace{.5pt}\La$
with
$$\xymatrixrowsep{18pt}\xymatrix{K_a^{-2} \ar[r]^-{\partial_a^{-2}} & K_a^{-1} \ar[r]^{\partial_a^{-1}} & K_a^0 \ar[r]^{\partial_a^0} & S_a \ar[r] & 0}$$ being a linear projective $2$-resolution of $S_a$ over ${\rm proj}\hspace{.6pt}\La$.

\end{Lemma}

\noindent{\it Proof.} Fix $a\in Q_0$. Observe that $R^{(n)}(a,x)$ is finite dimensional and vanishes for all but finitely many $x\in Q_0$. Thus, $K^{-n}_a\in {\rm proj}\hspace{.5pt}\La$ for all $n\ge 0$. Letting $\alpha_i: a \to b_i$, $i=1, \ldots, r$, be the arrows in $Q_1(a, -)$, we see that $K^{-1}_a=\otimes_{i=1}^r P_{b_i} \otimes k \alpha_i$. Let $\Oa$ be a minimal generating set for $R$, and let $\rho_j\in \Oa(a, c_j)$ with $c_j\in Q_0$, $j=1, \ldots, s$, be the relations in $\Oa(a, -)$. Then, $K^{-2}_a= \oplus_{j=1}^s P_{c_j}\otimes k\rho_j$. Writing $\rho_j=\sum_{j=1}^r \gamma_{ij}\hspace{.5pt}\alpha_i$ for some $\gamma_{ij}\in kQ_1(b_i, c_j)$, in view of Lemma \ref{diff}, we obtain a commutative diagram
$$\xymatrixcolsep{32pt}\xymatrix{
P_{c_1}\oplus \cdots \oplus P_{c_s} \ar[r]^-{\left(P[\bar{\gamma}_{ij}]\right)_{r\times s}} \ar[d]^{f_2} & P_{b_1} \oplus \cdots \oplus P_{b_r} \ar[rr]^-{(P[\bar{\alpha}_1],\cdots, P[\bar{\alpha}_r])} \ar[d]^{f_1} && P_a \ar[r]^{d_a} \ar[d]^{f_a} & S_a \ar[r] \ar@{=}[d] & 0\\
K^{-2}_a\ar[r]^{\partial_a^{-2}} & K^{-1}_a \ar[rr]^{\partial_a^{-1}} && K^0_a\ar[r]^{\partial_a^{\hspace{.4pt}0}} & S_a \ar[r]  & 0,}$$ where $f_a, f_1, f_2$ are graded isomorphisms defined in such a way that $f_a(e_a)=e_a\otimes \varepsilon_a;$ $f_1(e_{b_i})=e_{b_i}\otimes \alpha_i$, $i=1, \ldots, r;$ and $f_2(e_{c_j})=e_{c_j}\otimes \rho_j$, $j=1, \ldots, s$. By Lemma \ref{ex-p-pres}, we conclude that the lower row is a linear projective $2$-resolution of $S_a$.

Let now $n > 2$. Consider $v\in P_x$ and $\gamma \in R^{(n)}(a, x)$, where $x\in Q_0$. To show that $(\partial_a^{1-n}\circ \partial_a^{-n})(v\otimes \gamma)=0$, we may assume that $\gamma=\rho \delta,$ for some $\rho\in R_2(z, x)$ and $\delta\in kQ_{n-2}(a, z)$ with $z\in Q_0$. Write $\rho=\sum_{i=1}^s\lambda_i\beta_i\alpha_i$, for some $\lambda_i\in k$, $\alpha_i\in Q_1(z, y_i)$ and $\beta_i\in Q_1(y_i, x)$ with $y_i\in Q_0$.
By Lemma \ref{diff}, we obtain

$$\textstyle(\partial_a^{1-n}\circ \partial_a^{-n})(v\otimes \gamma)={\sum}_{i=1}^s (\partial_a^{1-n}\circ \partial_a^{-n})\left(v\otimes \lambda_i\beta_i\alpha_i \delta\right)=v \left({\sum}_{i=1}^s \lambda_i\bar{\beta}_i\bar{\alpha}_i\right)\otimes \delta=0.$$
The proof of the lemma is completed.

\medskip

In the sequel, $K_a^\cdt$ will be called the {\it local Koszul complex} of $\La$ at $a$. The following result is a local version of a well known result stated in \cite{BGS}.

\medskip

\begin{Theo}\label{Koz-proj-rls}

Let $\La=kQ/R$, where $Q$ is a locally finite quiver and $R$ is a quadratic ideal in $kQ$. If $a\in Q_0$, then $S_a$ has a linear projective resolution if and only if $K_a^\cdt$ is a projective resolution of $S_a$, that is, $K_a^\cdt$ is exact at all negative degrees.

\end{Theo}

\noindent{\it Proof.} The sufficiency is evident. Suppose that $S_a$ has a linear projective resolution over ${\rm proj}\hspace{.5pt}\La$. By Lemmas \ref{gr-proj-c} and \ref{K-cplx}, there exists a commutative diagram
$$\xymatrixrowsep{18pt}\xymatrix{
\cdots \ar[r] & P^{-p-1} \ar[r]^-{d^{-p-1}} & P^{-p} \ar[r]^{d^{-p}} \ar[d]^{f^{-p}} & P^{1-p} \ar[r] \ar[d]^{f^{1-p}} & \cdots \ar[r] & P^{-1} \ar[r]^{d^{-1}} \ar[d]^{f^{-1}} & P^0 \ar[r] \ar[d]^{f^{0}} & 0\\
\cdots \ar[r] & K_a^{-p-1}\ar[r] \ar[r]^-{\partial_a^{-p-1}} & K_a^{-p} \ar[r]^-{\partial_a^{-p}} & K_a^{1-p} \ar[r] & \cdots \ar[r] & K_a^{-1} \ar[r]^{\partial_a^{-1}} & K_a^0 \ar[r]  & 0,}$$
where $p\ge 2$, the upper row is a linear projective resolution of $S_a$, and $f^{-p}, \cdots, f^0$ are graded isomorphisms. In particular, $\partial_a^i$ co-restricts to a projective $J$-cover of ${\rm Ker}(\partial_a^{i-1})$, for $i=1, \ldots, p$.
We claim that $\partial^{-p-1}_a$ co-restricts to projective $J$-cover of ${\rm Ker}(\partial_a^{-p})$. We may assume that  $K_a^{-p}$ is non-zero. Then, $K_a^{1-p}=\oplus_{i=1}^m P_{x_i}\otimes k\zeta_i$, where
$\zeta_i\in R^{(p-1)}(a, x_i)$, $i=1, \ldots, m$, form a basis of $R^{(p-1)}(a, -)$. Moreover, $K_a^{-p}=\oplus_{j=1}^n P_{y_j}\otimes k\rho_j$, where
$\rho_j\in R^{(p)}(a, y_j)$, $j=1, \ldots, n$, form a basis of $R^{(p)}(a, -)$. Observe that $f^{-p}\circ d^{-p-1}$ is a projective $J$-cover of ${\rm Ker}(\partial_a^{-p})$, which is homogeneous of degree one. By Lemma \ref{p-cover}, ${\rm Ker}(\partial_a^{-p})$ admits a normalized $J$-top basis $T^p$, consisting of homogeneous elements of degree one.

Consider $u=(u_1, \ldots, u_n)\in T^p \cap \hspace{.5pt} e_z K^{-p}_a$, where $z\in Q_0$ and $u_j\in P_{y_j}\otimes k\rho_j$. Then $u_j=\bar{\gamma\hspace{1.6pt}}_{\hspace{-1pt}j} \otimes \rho_j,$ where $\gamma_j \in kQ_1(y_j, z); j=1, \ldots, n.$
By the definition of $R^{(p)}(a, y_j)$, we may write $\rho_j={\sum}_{i=1}^m \delta_{ij} \, \zeta_i,$ where $\delta_{ij} \in kQ_1(x_i, y_j).$ \vspace{1pt} Since $\partial_a^p(u)=0$, by Lemma \ref{diff},
${\sum}_{i=1}^m  ({\sum}_{j=1}^n \bar{\gamma}_j \bar{\delta}_{ij}) \otimes \zeta_i={\sum}_{j=1}^n {\sum}_{i=1}^m \partial_a^p (\bar{\gamma}_j \otimes \delta_{ij} \zeta_i )=  0 .$ Since the $\zeta_i$ are linearly independent, ${\sum}_{j=1}^n \bar{\gamma}_j \bar{\delta}_{ij} =0$, that is, $\eta_i={\sum}_{j=1}^n \gamma_j \delta_{ij}\in R(x_i, z)$, and consequently, $\eta_i\in R_2(x_i,z),$ for $i=1, \ldots, m$. Set $\omega ={\sum}_{j=1}^n \gamma_j \rho_j.$ Then $\omega=\sum_{i=1}^m\eta_i\zeta_i$, where $\eta_i\in R_2(x_i, z)$ and $\zeta_i\in kQ_{p-1}(a, x_i)$. Given $0\le s< p-1$, since $\rho_j\in R^{(p)}(a, y_j)$, we see that $\omega=\sum \eta_l \rho_l \delta_l$, where $\eta_l\in kQ_{p-2-s}(-, y_j)$, $\rho_l\in R_2$, $\delta_l\in kQ_s(a, -)$. This shows that $\omega \in R^{(p+1)}(a, z).$ In particular, $e_z\otimes \omega \in K_a^{p+1}$.

Let $f_i$ be the composite of $\partial^{-p-1}(y_i, z): P_z\otimes R^{(p+1)}(a, z)\to P_{y_i}\otimes R^{(p)}(a, y_i)$ and the canonical projection $p_i: P_{y_i}\otimes R^{(p)}(a, y_i) \to P_{y_i}\otimes k \rho_i,$ for $i=1, \ldots, n$. Since $\gamma_j \in kQ_1(y_j, z)$, given any $1\le i\le n$, we deduce from Lemma \ref{diff} that
$$\textstyle f_i(e_z\otimes \omega) \hspace{-2pt} = \hspace{-2pt} p_i ({\sum}_{j=1}^n \partial^{p+1}_a(y_i, z)(e_z\otimes \gamma_j \rho_j)) \hspace{-2pt} = \hspace{-2pt} p_i({\sum}_{y_j=y_i} \bar{\gamma}_j \otimes \rho_j) \hspace{-2pt} = \hspace{-2pt} \bar{\gamma}_i \otimes \rho_i.$$
Thus, $\partial^{-p-1}_a(e_z\otimes \omega)=(f_1(e_z\otimes \omega), \ldots, f_n(e_z\otimes \omega))=(u_1, \ldots, u_n)=u.$ This establishes our claim. In view of Lemma \ref{gr-proj-c}, we obtain a graded isomorphism $f^{-p-1}: P^{-p-1}\to K_a^{-p-1}$ such that $f^{-p}\circ d^{-p-1}= \partial^{-p-1}_a \circ f^{-p-1}.$ By induction, $K_a^\cdt$ is a linear projective resolution of $S_a$. The proof of the theorem is completed.

\medskip

Next, we shall define the quadratic dual of $\La$. Given $x, y\in Q_0$ and $n\ge 0$, since $Q_n(x,y)$ is finite, $D(kQ_n(x,y))$ has a basis $\{\xi^* \mid \xi\in Q_n(x, y)\}$, that is the dual basis of $Q_n(x, y)$. If $\gamma=\sum \lambda_i \xi_i$, with $\lambda_i\in k$ and $\xi_i\in Q_n(x, y)$, we shall write $\gamma^*=\sum \lambda_i \xi_i^*\in D(kQ_n(x, y))$. The following statement will be needed later.

\medskip

\begin{Lemma}\label{dual-comp}

Let $\zeta\in kQ_s(x, y)$ and $\gamma\in Q_t(y, z)$, where $x, y\in Q_0$ and $s, t\ge 0$.

\begin{enumerate}[$(1)$]

\item If $\delta\in kQ_s$ and $\xi\in kQ_t$, then $(\gamma\zeta)^*(\xi\delta)=\gamma^*\hspace{-1pt}(\xi) \hspace{.5pt} \zeta^*\hspace{-1pt}(\delta).$

\vspace{1pt}

\item If $\gamma\in Q_1(y, z),$ then $(\gamma\zeta)^*\hspace{-1.5pt}(\eta)=\zeta^*\hspace{-1.5pt}(\partial_\alpha(\eta))$ for $\eta \hspace{-1.5pt} \in kQ_{s+1}(x, z).$

\end{enumerate} \end{Lemma}

\noindent{\it Proof.} Statement (1) is evident. Let $Q_1(-, z)=\{\alpha_i: y_i\to z \mid i=1, \ldots, r\}$, where $\alpha_1=\gamma \in Q_1(y, z)$. If $\eta\in kQ_{s+1}(x, z),$ then $\eta=\sum_{i=1}^r \alpha_i \delta_i,$ where $\delta_i\in kQ_s(x, y_i)$ and $\partial_\gamma(\eta)= \delta_1$. Hence, $(\gamma\zeta)^*(\eta)={\sum}_{i=1}^r \gamma^*(\alpha_i) \zeta^*(\delta_i) =
\zeta^*(\delta_1)=\zeta^*(\partial_\gamma(\eta)).$
The proof of the lemma is completed.

\medskip

Consider now $kQ^{\rm o}$, where $Q^{\rm o}$ is of the opposite quiver of $Q$. Given $x, y\in Q_0$ and $n\ge 0$, we have a $k$-linear isomorphism $\psi^n_{x, y}: kQ_n^{\rm o}(y, x)\to D(kQ_n(x,y)): \gamma^{\rm o}\mapsto \gamma^*.$
Taking $n=2$, we denote by $R^!_2(y, x)$ the subspace of $kQ^{\rm o}_2(y, x)$ of elements $\rho^{\rm o}$, with $\rho\in kQ_2(x, y)$, such that $\rho^*\in R_2(x, y)^\perp$, that is, $\rho^*$ vanishes on $R_2(x, y)$. The {\it quadratic dual} of $\La$ is then defined by
$\La^!=kQ^{\rm o} / R^!,$ where $R^!$ is the ideal in $kQ^{\rm o}$ generated by all $R^!_2(y, x)$ with $x, y\in Q_0$; compare \cite[page 69]{Mar} and \cite[(2.8.1)]{BGS}.

\medskip

\begin{Prop}\label{q-dual}

Let $\La=kQ/R$, where $Q$ is a locally finite quiver and $R$ is a quadratic ideal. Then $\La^!$ and $\La^{\rm o}$ are quadratic with $(\La^!)^!= \La$ and $(\La^{\rm o})^!=(\La^!)^{\rm o}$.

\end{Prop}

\noindent{\it Proof.} Fix $x, y\in Q_0$. Consider $\psi_{y, x}^{2, {\rm o}}: kQ_2(x, y) \to D(kQ_2^{\rm o}(y, x)): \gamma \mapsto (\gamma^{\rm o})^*,$ where $(\gamma^{\rm o})^*(\rho^{\rm o})=\rho^*(\gamma)=\gamma^*(\rho)$, for all $\rho\in kQ_2(x, y)$. Now, $\gamma\in (R^!)^!_2(x, y)$ if, and only if, $(\gamma^{\rm o})^*(\rho^{\rm o})=0$, that is $\rho^*(\gamma)=0$ for all $\rho^{\rm o}\in R^!(y, x)$ if, and only if, $\gamma\in R_2(x, y)$. Moreover, $\gamma\in (R^{\rm o})^!_2(x, y)$ if and only if $(\gamma^{\rm o})^*(\rho^{\rm o})=0$, that is $\gamma^*(\rho)=0$ for all $\rho\in R_2^{\rm o}(y, x)$ if, and only if, $\gamma^{\rm o}\in R^!(y, x)$  if, and only if, $\gamma\in (R^!)^{\rm o}(x, y)$. Thus, $(R^!)^!=R$ and $(R^{\rm o})^!=(R^!)^{\rm o}.$ The proof of the proposition is completed.

\medskip

By definition, we obtain $R^!_2(y, x)\cong R_2(x, y)^\perp=R^{(2)}(x,y)^\perp$, for $x, y\in Q_0$. The following statement describes $R^!_n(y, x)=R^!(y, x)\cap kQ_n^{\rm o}(y, x)$ for all $n \ge 0$.

\medskip

\begin{Lemma}\label{perp}

Let $\La=kQ/R$, where $Q$ is locally finite and $R$ is quadratic. Given $n\ge 0$, we obtain a $k$-linear isomorphism
$\psi_{x,y}^n: R^!_n\hspace{-1pt}(y,x)\to R^{(n)}(x,y)^\perp: \gamma^{\rm o} \mapsto \gamma^*.$

\end{Lemma}

\noindent{\it Proof.} We need only to consider the case where $n\ge 3$. Fix $x, y\in Q_0$. We have a $k$-linear isomorphism
$\psi^n_{x, y}: kQ_n^{\rm o}(y, x)\to D(kQ_n(x,y)): \gamma^{\rm o}\mapsto \gamma^*.$ By definition, $R^!_n(y, x) = {\sum}_{j=0}^{n-2} R^!_{n,j}(y, x)$ and $R^{(n)}(x, y) = \cap_{j=0}^{n-2}\, R^{(n,j)}(x,y),$ where $$\textstyle R^!_{n,j}(y, x)={\sum}_{a,b\in Q_0} kQ^{\rm o}_j(a, x) \cdot R^!_2(b,a) \cdot kQ^{\rm o}_{n-j-2}(y,b)$$ and $\textstyle R^{(n,j)}(x,y)={\sum}_{a, b\in Q_0}\, kQ_{n-j-2}(b,y) \cdot R_2(a,b) \cdot kQ_j(x,a).$

\vspace{1pt}

First, we claim that $\psi_{x, y}^n(R^!_n(y,x))\subseteq R^{(n)}(x,y)^\perp.$ Indeed, consider $\sigma^{\rm o}\in R^!_n(y,x)$, where $\sigma\in kQ_n(x, y)$. We may assume that $\sigma^{\rm o}=\gamma^{\rm o} \eta^{\rm o} \delta^{\rm o}$, where $\delta \in kQ_{n-2-j}(b, y)$, and $\eta \in kQ_2(a, b)$ such that $\eta^{\rm o}\in R^!_2(b, a)$, and $\gamma \in kQ_j (x, a)$, for some $j$ with $0\le j\le n-2$.
Given any $w \in R^{(n)}(x, y)$, we may write $w = \sum_{i} \delta_i \cdot \rho_i \cdot \gamma_i,$ for some $\gamma_i\in kQ(x,a_i)$, and $\rho_i\in R_2(a_i, b_i)$, and $\delta_i\in kQ(b_i, y)$. By Lemma \ref{dual-comp}(1), $\sigma^*(w)={\sum}_{i} \delta^*\hspace{-1pt}(\delta_i)  \, \eta^*\hspace{-1pt}(\rho_i)  \,\gamma^*\hspace{-1pt}(\gamma_i) =0.$ This establishes our claim.

On the other hand, by Lemma \ref{orthogonal}(1), $R^{(n)}(x,y)^\perp = {\sum}_{j=0}^{n-2}\, R^{(n,j)}(x,y)^\perp,$ where $R^{(n,j)}(x,y)^\perp=\cap_{\,a, b\in Q_0}\;( kQ_{n-j-2}(b,y) \cdot R_2(a,b) \cdot kQ_j(x,a))^\perp.$
Let $f\in R^{(n)}(x,y)^\perp$ be non-zero. We may assume that $f\in R^{(n,j)}(x,y)^\perp,$ for some $0\le j\le n-2.$ Then, $f=\omega^*$, with $\omega\in kQ_n(x,y).$ Assume further that $\omega= \delta \zeta  \gamma$, where $\gamma\in kQ_j(x, a)$, $\zeta\in kQ_2(a, b)$, and $\delta \in kQ_{n-j-2}(b, y)$, for some $a, b\in Q_0$. Since $f\ne 0$, there exist $\gamma_0\in kQ_j(x, a)$ and $\delta_0 \in kQ_{n-j-2}(b, y)$ such that $\gamma^*(\gamma_0)=\delta^*(\delta_0)=1$.

Consider a basis $\{\rho_1, \ldots, \rho_r; \rho_{r+1}, \ldots, \rho_s\}$ of
$kQ_2(a, b)$, where $\{\rho_1, \ldots, \rho_r\}$ is a basis of $R_2(a, b)$. Then, $kQ_2(a, b)$ has a basis $\{\eta_1, \ldots, \eta_r; \eta_{r+1}, \ldots, \eta_s\}$ such that $\{\eta_1^*, \ldots, \eta_r^*; \eta_{r+1}^*, \ldots, \eta_s^*\}$ is the dual basis of $\{\rho_1, \ldots, \rho_r; \rho_{r+1}, \ldots, \rho_s\}$. Then, $\{\eta_{r+1}^{\rm o}, \ldots, \eta_s^{\rm o}\}$ is a basis of $R^!_2(b, a)$. Now, we may write $\zeta=\sum_{i=1}^s \lambda_i \eta_i$, with $\lambda_i\in k\vspace{1pt}$. Then, $f={\sum}_{i=1}^s \lambda_i (\delta \eta_i \gamma)^*.$ Since $f \in \left( kQ_{n-j-2}(b,y) \cdot R_2(a,b) \cdot kQ_j(x,a)\right)^\perp,$
\vspace{1.5pt} by Lemma \ref{dual-comp}(1),  $0=f(\delta_0\rho_j\gamma_0)={\sum}_{i=1}^s \lambda_i \hspace{.5pt} \delta^*\hspace{-1pt}(\delta_0) \hspace{.5pt} \eta_i^*\hspace{-1pt}(\rho_j) \hspace{.5pt} \gamma^*\hspace{-1pt}(\gamma_0)=\lambda_j,$ $j=1, \ldots, r.$ \vspace{1.5pt} Thus, $f={\sum}_{i>r} \lambda_i (\delta \eta_i \gamma )^*=\psi_{x,y}^n\left({\textstyle\sum}_{i>r} \lambda_i \gamma^{\rm o} \eta_i^{\rm o} \delta^{\rm o}\right).$ The proof of the lemma is completed.

\medskip

We shall fix some notation for the quadratic dual $\La^!$. For $\gamma\in kQ$, we write $\gamma^!=\gamma^{\rm o}+ R^! \in \La^!$. Let $x, a\in Q_0$. We write again $e_x=\varepsilon_x+R^!$ and put $P^!_x=\La^!e_x$. In particular, $e_a\La_n^!e_x=\{\gamma^! \mid \gamma\in kQ_n(a, x)\}$ for all $n\ge 0$. Given $\alpha\in Q_1(y, x)$ and $n>0$, the $\La^!$-linear morphism $P[\alpha^{!}]: P^!_y\to P^!_x$ restricts to a $k$-linear map $e_a\La^!_{n-1}e_y\to e_a\La^!_{n}e_x$ which, for the simplicity of notation, will be denoted again by $P[\alpha^!]$. Now, we define a sequence $L_a^\cdt$ of morphisms in ${\rm proj}\hspace{.5pt}\La$ as follows:
$$\xymatrixrowsep{18pt}\xymatrix{\cdots \ar[r] & L_a^{-n} \ar[r]^-{d^{-n}} & L_a^{1-n} \ar[r] &
\cdots \ar[r] & L_a^{-1} \ar[r]^{d^{-1}} & L_a^0 \ar[r] & 0 \ar[r] & \cdots ,}$$ where $L_a^{-n}=\oplus_{x\in Q_0} P_x \otimes D(e_a\La^!_ne_x)$ for $n\ge 0$, and $d^{-n}=(d^{-n}(y,x))_{(y,x)\in Q_0\times Q_0}$ with
$d^{-n}(y,x)={\textstyle\sum}_{\alpha\in Q_1(y,x)} P[\bar{\alpha}]\otimes DP[\alpha^!]: P_x \otimes  D(e_a\La^!_ne_x)\to P_y \otimes  D(e_a\La^!_{n-1}e_y)$ for all $n>0$.

\smallskip

\begin{Lemma}\label{k-cplx-iso}

Let $\La=kQ/R$, where $Q$ is a locally finite quiver and $R$ is a quadratic ideal. If $a\in Q_0$, then $L_a^\cdt$ is isomorphic to the local Koszul complex of $\La$ at $a$.

\end{Lemma}

\noindent{\it Proof.} Given $a, x\in Q_0$ and $n\in \Z$, we obtain $D(kQ_n(a, x))=\{\gamma^* \mid \gamma\in kQ_n(a, x) \}$ and $e_a\La^!_ne_x=\{\gamma^! \mid \gamma\in kQ_n(x, a)\} \cong kQ_n^{\rm o}(x, a)/R^!_n(x, a)$. It follows from Lemma \ref{perp} that $\gamma^*\in  R^{(n)}(a,x)^\perp$ if and only if $\gamma^{\rm o}\in R^!_n(x, a)$. This leads to a bilinear form
$$<\hspace{-2.5pt}-, -\hspace{-2.5pt}>: R^{(n)}(a,x)\times e_a\La^!_ne_x \to k: (\delta, \gamma^!) \mapsto \gamma^*(\delta),$$ which, by definition, is right non-degenerate. If $\delta\in R^{(n)}(a,x)$ is non-zero, then $\gamma^*(\delta)\ne 0$ for some $\gamma\in kQ_n(a, x)$, that is, $<\hspace{-2.5pt}\delta, \gamma^!\hspace{-2.5pt}>\ne 0$. Hence, $<\hspace{-2.5pt}-, -\hspace{-2.5pt}>$ is non-degenerate. Since $R^{(n)}(a,x)$ and $e_a\La^!_ne_x$ are finite dimensional, we obtain a $k$-linear isomorphism $\phi^n_{x}: R^{(n)}(a,x)\to D(e_a\La^!_ne_x): \delta \to {\hspace{-4pt}}<\hspace{-2.5pt}\delta, -\hspace{-3pt}>.$ We claim, for $x, y\in Q_0$ and $n>0$, that the following diagram commutes:
$$
\xymatrixrowsep{20pt}\xymatrix{
R^{(n)}(a,x) \ar[rrr]^-{\textstyle\sum_{\alpha\in Q_1(y,x)}\partial_\alpha} \ar[d]_{\phi^n_{x}} &&& R^{(n-1)}(a, y)\ar[d]^{\phi^{n-1}_y}\\
D(e_a\La^{!}_ne_x)  \ar[rrr]^-{\sum_{\alpha\in Q_1(y,x)}DP[\alpha^!]} &&& D(e_a\La^{!}_{n-1}e_y).
}$$ Indeed, given $\delta\in R^{(n)}(a, x)$ and $\zeta\in kQ_{n-1}(a, y)$, by Lemma \ref{dual-comp}(2),
we obtain
$$\begin{array}{rcl}
(\textstyle\sum_{\alpha\in Q_1(y, x)}D P[\alpha^!])(\phi^n_{x}(\delta))(\zeta^!)
&=&\textstyle\sum_{\alpha\in Q_1(y, x)} \phi^n_{x}(\delta)(\zeta^!\alpha^!)\vspace{1pt}\\
&=&\textstyle\sum_{\alpha\in Q_1(y, x)} (\alpha\zeta)^*(\delta)\vspace{1pt}\\
&=&\textstyle\sum_{\alpha\in Q_1(y, x)}  \zeta^*(\partial_\alpha(\delta)) \vspace{1pt}\\
&=& \phi^{n-1}_{x}({\textstyle\sum}_{\alpha\in Q_1(x,y)} \partial_\alpha(\delta))(\zeta^!).
\end{array}$$
Our claim is established. As a consequence, we obtain a commutative diagram
$$
\xymatrix{
\textstyle\oplus_{x\in Q_0}P_x\otimes R^{(n)}(a, x) \ar[r]^-{\partial^{-n}_a} \ar[d]_{\oplus \hspace{.5pt} 1_{\hspace{-1.5pt}P_x}\otimes \phi^n_{x}} & \oplus_{y\in Q_0} P_y\otimes R^{(n-1)}(a,y)
\ar[d]^{\oplus\hspace{.5pt} 1_{\hspace{-1.5pt}P_y}\otimes \phi^{n-1}_{y}}\\
\oplus_{x\in Q_0} P_x\otimes D(e_a\La^{!}_ne_x) \ar[r]^-{d^{-n}} & \oplus_{y\in Q_0} P_y\otimes D(e_a\La^{!}_{n-1}e_y)
}
$$ with vertical isomorphisms, for all $n>0$. The proof of the lemma is completed.

\medskip

The following result is a generalization of Proposition 2.9.1 stated in \cite{BGS}, where $\La$ is assumed to be locally finitely generated; see also \cite[Theorem 30]{MOS}.

\medskip

\begin{Theo}\label{Dual-Koszul}

Let $\La=kQ/R$, where $Q$ is a locally finite quiver and $R$ is a quadratic ideal. Then $\La$ is Koszul if and only if $\La^{!}$ is Koszul.

\end{Theo}

\noindent{\it Proof.} \vspace{1pt} Since $(\La^!)^!=\La$, it suffices to prove the necessity. Suppose that $\La$ is Koszul. Fix $a\in Q_0$. By Proposition \ref{q-dual}, $e_a(\La^!)^!_ne_x=e_a\La_ne_x$, for all $x\in Q_0$ and $n\ge 0.$ By Proposition \ref{k-cplx-iso}, the local Koszul complex of $\La^!$ at $a$ is isomorphic to
$$L^\cdt: \quad \xymatrix{\cdots\ar[r] & L^{-n}\ar[r]^{d^{-n}} & L^{1-n}\ar[r] & \cdots \ar[r] & L^{-1}\ar[r]^{d^{-1}} & L^0
\ar[r] & 0 \ar[r] & \cdots }$$ where $L^{-n}=\oplus_{x\in Q_0}P^!_x\otimes D(e_a\La_ne_x)$ for $n\ge 0$, and $d^{-n}=(d^{-n}(y,x))_{(y,x)\in Q_0\times Q_0}$ with $d^{-n}(y,x)={\sum}_{\alpha\in Q_1(x, y)}P[\alpha^!]\otimes DP[\bar \alpha]: P^!_x\otimes D(e_a\La_ne_x)\to P^!_y\otimes D(e_a\La_{n-1}e_y)$ for $n>0$. By Theorem \ref{Koz-proj-rls}, it suffices to show that $L^\cdt$ is exact at every negative degree $n$. Indeed, the $\La^!$-linear morphism $d^{-n}: L^{-n}\to L^{1-n}$ consists of a family of $k$-linear maps $d^{-n}(b): L^{-n}(b)\to L^{1-n}(b)$, with $b\in Q_0$.  Since $L^\cdt$ is a linear complex, by Lemma \ref{grad-exact}, it amounts to establish the exactness of the sequence
$$(*) \quad \xymatrixrowsep{8pt}\xymatrixcolsep{35pt}\xymatrix{
\oplus_{x\in Q_0}e_b\La^!_{s-1} e_x\otimes D(e_a\La_{n+1}e_x)
\ar[r]^-{d^{-n-1}(b)} & \oplus_{y\in Q_0}e_b\La^!_s e_y\otimes D(e_a\La_n e_y) \hspace{20pt} \\  \hspace{135pt}\ar[r]^-{\hspace{5pt} d^{-n}(b)} &
\oplus_{z\in Q_0} e_b \La^!_{s+1} e_z\otimes D(e_a\La_{n-1}e_z),}$$ for all $b\in Q_0$ and $s\in \Z$. If $s<0$, then $e_b\La^!_s e_y=0$, and hence, $(*)$ is exact. In case $s=0$, the sequence $(*)$ becomes
$$\xymatrix{0 \ar[r] & e_b\La^!_0 e_b \otimes D(e_a \La_n e_b) \ar[r]^-{d_b^{-n}} & \oplus_{z\in Q_0} e_b \La^!_1 e_z\otimes D(e_a\La_{n-1}e_z),}$$ where $d_b^{-n}=(d_b^{-n}(z,b))_{z\in Q_0}$ with $d_b^{-n}(z,b)=\sum_{\alpha\in Q_1(b, z)} P[\alpha^!] \otimes DP[\bar\alpha]$. Consider a non-zero linear form $f\in D(e_a \La_n e_b)$. Then $f(u \bar\beta) \ne 0$, for some $\beta\in Q_1(b, z)$ and $u\in e_a\La_{n-1} e_z$ with $z\in Q_0$. That is, $(DP[\bar\beta])(f)(u)\ne 0,$ and consequently, $(DP[\bar\beta])(f)\ne 0.$ Now, $d_b^{-n}(z,b)(e_b\otimes f)={\textstyle\sum}_{\alpha\in Q_1(b, z)}\, \alpha^!\otimes (DP[\bar\alpha])(f),$ which is non-zero since the $\alpha^!$ with $\alpha\in Q_1(b, z)$ are $k$-linearly independent. This shows that $d^{-n}(b)$ is a monomorphism. That is, $(*)$ is exact.

Finally, consider the case $s>0$. Since $\La$ is Koszul, by Theorem \ref{Koz-proj-rls} and Proposition \ref{k-cplx-iso}, $L_b^\cdt$ is exact at degree $s$. By Lemma \ref{grad-exact}, we obtain an exact sequence
$$\xymatrixrowsep{8pt}\xymatrix{\oplus_{z\in Q_0}e_a \La_{n-1} e_z \otimes D(e_b\La^!_{s+1} e_z) \ar[r]^-{\hspace{8pt}d_{b,a}^{\hspace{.5pt}-s-1}} & \oplus_{y\in Q_0}e_a \La_n e_y\otimes D(e_b\La^!_s e_y) \hspace{20pt} \\  \hspace{140pt}\ar[r]^{d_{b,a}^{-s}} &
\oplus_{x\in Q_0} e_a\La_{n+1} e_x \otimes D(e_b\La^!_{s-1}e_x),}$$
where $d_{b,a}^{-s}=(d_{b,a}^{-s}(x, y))_{(x, y)\in Q_0\times Q_0}$ \vspace{1pt} with $d_{b,a}^{-s}(x, y)=\textstyle{\sum}_{\alpha\in Q_1(x,y)} P[\bar\alpha]\otimes DP[\alpha^!].$ Applying the duality $D$ to the above sequence, by Lemma \ref{dual-mix}, we obtain an exact sequence isomorphic to the sequence $(*)$. The proof of the theorem is completed.

\medskip

\noindent{\sc Remark.} In case $\La$ is Koszul, one calls $\La^!$ the {\it Koszul dual} of $\La$.

\medskip

Given $M, N\in {\rm Mod}\hspace{.5pt}\La$, one can define in a canonical way the Yoneda Ext-group $\Ext_{_{\hspace{-1pt}\mathit\Lambda}}^n(M,N)$, for all $n\ge 0$; see \cite[Section III.5]{Mac}.

\medskip

\begin{Lemma}\label{Ext-Simple}

Let $\La=kQ/R$ be a locally finite dimensional Koszul algebra. Given $a,b\in Q_0$ and $n\ge 0$, we obtain $\Ext_{_{\hspace{-1pt}\mathit\Lambda}}^n(S_b,S_a)= e_b\La_n^{!}e_a$.

\end{Lemma}

\noindent{\it Proof.} Fix $a, b\in Q_0.$ By Proposition \ref{Alg-dec}(3), $\La$ is strong locally finite dimensional, and by Theorem \ref{Koz-proj-rls} and Lemma \ref{k-cplx-iso}, $L_b^\cdt$ is a minimal projective resolution of $S_b$. Given $n\ge 0$, it is well known that $\Ext_{_{\hspace{-1pt}\mathit\Lambda}}^n(S_b, S_a)\cong \Hom_{_{\hspace{-1pt}\mathit\Lambda}}(L_b^{-n}, S_a)$; see, for example, \cite[(III.6.4)]{Mac}. Since $e_b\La_n^{!}e_a$ is finite dimensional, by Proposition \ref{Mor}(3), we see that
$$\Ext_{_{\hspace{-1pt}\mathit\Lambda}}^n(S_b, S_a)\cong \Hom_{_{\hspace{-1pt}\mathit\Lambda}}(P_a\otimes D(e_b\La_n^{!}e_a), S_a) \cong {\rm Hom}_k(D(e_b\La_n^{!}e_a), k)\cong e_b\La_n^{!}e_a.$$
The proof of the lemma is completed.

\medskip

In order to study the opposite algebra of a Koszul algebra, we shall introduce some notation for $(\La^!)^{\rm o}=kQ/(R^!)^{\rm o}$. Given $\gamma\in kQ$, we shall write $\hat{\gamma\hspace{1.5pt}}=\gamma+ (R^!)^{\rm o}$. Given $x\in Q_0$, however, we shall write $e_x=\varepsilon_x+(R^!)^{\rm o}$ for the simplicity of notation. We shall associate with each $a\in Q_0$ a sequence $T_a^\pdt$ of morphisms in ${\rm inj}\hspace{.8pt}\La$ as follows:
$$\xymatrix{\cdots\ar[r] & 0 \ar[r] & T_a^0 \ar[r]^{d^0} & T_a^1 \ar[r] & \cdots \ar[r] & T_a^n \ar[r]^{d^n} & T_a^{n+1} \ar[r] &  \cdots}$$ with $T_a^{n}=\oplus_{x\in Q_0} I_x \otimes e_x\La^!_ne_a$ and $d^{n}=(d^{n}(y,x))_{(y,x)\in Q_0\times Q_0}$ for all $n\ge 0$, where
$d^{n}(y,x)=\textstyle\sum_{\alpha\in Q_1(y,x)} I[\bar{\alpha}]\otimes P_a^!(\alpha^{\rm o}): I_x \otimes  e_x\La^!_ne_a\to I_y \otimes  e_y\La^!_{n+1}e_a.$

\medskip

\begin{Lemma}\label{Koz-cplx-dual}

Let $\La=kQ/R$ be a locally finite dimensional quadratic algebra. If $a\in Q_0$, then $T_a^\pdt$
is isomorphic to the dual of the local Koszul complex of $\La^{\rm o}$ at $a$.

\end{Lemma}

\noindent{\it Proof.} Fix $a\in Q_0$. By Lemma \ref{q-dual}, $e_a(\La^{\rm o})^!_ne_x=e_a(\La^!)^{\rm o}_ne_x$, for $x\in Q_0$ and $n\ge 0$. By Lemma \ref{k-cplx-iso}, the local Koszul complex of $\La^{\rm o}$ at $a$ is isomorphic to 
$$\xymatrix{L^\cdt: \qquad \cdots \ar[r] & L^{-n} \ar[r]^{d^{-n}} & L^{1-n}\ar[r] & \cdots \ar[r] & L^{-1}\ar[r]^{d^{-1}} & L_0 \ar[r] & 0,}$$ with $L^{-n}=\oplus_{x\in Q_0} P_x^{\rm o} \otimes D(e_a(\La^!)^{\rm o}_ne_x)$ and $d^{-n}=(d^{-n}(x,y))_{(x,y)\in Q_0\times Q_0}$, where
$d^{-n}(x,y)=\textstyle\sum_{\alpha\in Q_1(y,x)} P[\bar{\alpha}^{\rm o}]\otimes DP[\hat\alpha]:
P_y^{\rm o} \otimes D(e_a(\La^!)^{\rm o}_{n+1}e_y) \to P_x^{\rm o}\otimes D(e_a(\La^!)^{\rm o}_ne_x).$

\vspace{1pt}

Since $e_a(\La^!)^{\rm o}_ne_x$ \vspace{1pt} is finite dimensional, we have a canonical $k$-linear isomorphism $\theta_{x,n}: D^2(e_a(\La^!)^{\rm o}_ne_x) \to e_a(\La^!)^{\rm o}_ne_x$. Composing this with the $k$-linear isomorphism
$e_a(\La^!)^{\rm o}_ne_x\to e_x\La^!_ne_a$, sending $\hat{\gamma\hspace{1.5pt}}$ to $\gamma^!$ for $\gamma\in kQ_n(x,a)$, we obtain a $k$-linear isomorphism
$\psi_{x,n}: D^2(e_a(\La^!)^{\rm o}_ne_x) \to e_x\La^!_ne_a$
such, for any $\alpha\in Q_1(y, x)$, that
$$\xymatrix{
I_x \otimes D^2( e_a (\La^!)^{\rm o}_n e_x)
\ar[d]_{\id \otimes \psi_{x,n}} \ar[rr]^-{I[\bar\alpha]\otimes DP[\hat\alpha]}   &&
I_y \otimes D^2(e_a (\La^!)^{\rm o}_{n+1} e_y)\ar[d]^{\id \otimes \psi_{y,n+1}} \\
I_x\otimes e_x \La^!_n e_a \ar[rr]^-{I[\bar\alpha]\otimes P_a^!(\alpha^{\rm o})} &&
I_y\otimes e_y \La^!_{n+1} e_a
}$$ commutes. Since the $L^{-n}$ are finite direct sums, in view of the above commutative diagram, we see that $D(L^\cdt)\cong T_a^\pdt $. The proof of the lemma is completed.

\medskip

In case $\La$ is locally finite dimensional, we obtain the following generalization of Proposition 2.2.1 stated in \cite{BGS}.

\medskip

\begin{Theo}\label{Opp-Koszul}

Let $\La=kQ/R$ be a locally finite dimensional qudratic algebra. The following statements are equivalent.

\begin{enumerate}[$(1)$]

\item The algebra $\La$ is Koszul.

\vspace{1pt}

\item The opposite algebra $\La^{\rm o}$ is Koszul.

\vspace{1pt}

\item For each $a\in Q_0$, the simple module $S_a$ has $T_a^\ydt$ as an injective co-resolution.

\end{enumerate}

\end{Theo}

\noindent{\it Proof.} By Proposition \ref{Alg-dec}(3), $\La$ is strongly locally finite dimensional, and by Proposition \ref{proj-inj}, $T_a^\pdt$ is a complex of injective $\La$-modules. Assume that $T_a^\pdt$ is an injective co-resolution of $S_a$ for every $a\in Q_0$. By Lemma \ref{Koz-cplx-dual}, every local Koszul complex of $\La^{\rm o}$ is exact at all non-zero degrees, and by Theorem \ref{Koz-proj-rls}, $\La^{\rm o}$ is Koszul.

Suppose now that $\La$ is Koszul. Fix  $a\in Q_0$. Recall that the complex $(T_a^\pdt, d^\pdt)$ is such that
$T_a^i=\oplus_{x\in Q_0} I_x \otimes e_x\La^!_ie_a$ and $d^i=(d^i(y,x))_{(y,x)\in Q_0\times Q_0}$ for $i\ge 0$, where
$d^i(y,x)=\textstyle\sum_{\alpha\in Q_1(y,x)} I[\bar{\alpha}]\otimes P_a^!(\alpha^{\rm o}): I_x \otimes  e_x\La^!_ie_a\to I_y \otimes  e_y\La^!_{i+1}e_a.$ In particular, $T_a^0=I_a\otimes ke_a$ and $T_a^1=\oplus_{j=1}^s I_{b_j}\otimes k \beta_j^!$, where $\{\beta_1, \ldots, \beta_s\}=Q_1(a, -)$. Consider the $\La$-linear morphism $d^{-1}: S_a\to T^0: e_a+Je_a\mapsto e_a^\star\otimes e_a$. In view of Corollary \ref{inj-pres}, we may assume that there exists an exact sequence
$$\xymatrixcolsep{18pt}\xymatrix{0\ar[r] & S_a \ar[r]^-{d^{-1}} & T_a^0 \ar[r]^-{d^0} & T_a^1 \ar[r] & \cdots \ar[r] & T_a^{n-1} \ar[r]^-{d^{n-1}} & T_a^n \ar[r]^{p^n} & C^{n+1}\ar[r] & 0,}$$ where $n\ge 1$, such that $S(T_a^i)\subseteq {\rm Im}(d^{i-1})$, for $i=0, 1, \ldots, n$. By Corollary \ref{inj-soc} and Lemma \ref{ess-socle}, $C^{n+1}$ has an essential socle. Given $y\in Q_0$, in view of the proof of \cite[(III.6.4)]{Mac}, we see
that ${\rm Ext}^{n+1}_{_{\hspace{-1pt}\it\Lambda\hspace{-1pt}}}(S_y, S_a)\cong {\rm Hom}_{_{\hspace{-1pt}\it\Lambda\hspace{-1.5pt}}}(S_y, C^{n+1})/{\rm Im}({\rm Hom}_{_{\hspace{-1pt}\it\Lambda\hspace{-1pt}}}(S_y, p^n));$ see also \cite[(III.8.2)]{Mac}. Setting $C^n={\rm Im}(d^{n-1})$, we obtain a short exact sequence
$$\xymatrix{0\ar[r] & C^n\ar[r]^-{j^n} & T_a^n \ar[r]^-{p^n} & C^{n+1} \ar[r] & 0,}$$ where $j^n$ is an injective envelope. Applying ${\rm Hom}_{_{\hspace{-.5pt}\it\Lambda\hspace{-1pt}}}(S_y, -)$ yields an exact sequence
$$\xymatrixcolsep{16pt}\xymatrix{\Hom_{_{\hspace{-.5pt}\it\Lambda\hspace{-1pt}}}(S_y,C^n) \ar[r]^-{j^n_*} & \Hom_{_{\hspace{-.5pt}\it\Lambda\hspace{-1pt}}}(S_y, T_a^n) \ar[r]^-{p^n_*} & \Hom_{_{\hspace{-.5pt}\it\Lambda\hspace{-1pt}}}(S_y, C^{n+1})\ar[r] & \Ext^{n+1}_{_{\hspace{-.5pt}\it\Lambda}}\hspace{-1pt}(S_y,S_a)\ar[r] & 0.}
$$
Since $S_y$ is simple, $j^n_*$ is surjective. Thus, $\Hom_{_{\hspace{-.5pt}\it\Lambda\hspace{-1pt}}}(S_y, C^{n+1})\cong \Ext^{n+1}_{_{\hspace{-.5pt}\it\Lambda\hspace{-1pt}}}(S_y,S_a).$ In view of Lemma \ref{Ext-Simple}, ${\rm dim}_k\Hom_{_{\hspace{-.5pt}\it\Lambda\hspace{-1pt}}}(S_y, C^{n+1})={\rm dim}_k\,e_y\La^{!}_{n+1}e_a,$ and consequently, $S(C^{n+1})\cong \oplus_{y\in Q_0}S_y\otimes e_y\La^{!}_{n+1}e_a$, which is finite dimensional. In view of Lemma \ref{i-envelope}, we obtain an injective envelope
$j^{n+1}: C^{n+1} \to \oplus_{y\in Q_0} I_y \otimes e_y\La^{!}_{n+1}e_a=T_a^{n+1}.$ Now, we claim that the sequence
$$\xymatrix{0\ar[r] & S_a \ar[r]^-{d} & T_a^0 \ar[r]^-{d^0} & \cdots \ar[r]^-{d^{n-2}}  & T_a^{n-1} \ar[r]^-{d^{n-1}} & T_a^n \ar[r]^-{d^n} & T_a^{n+1}}$$
is exact with $S(T_a^{n+1})\subseteq {\rm Im}(d^n)$, or equivalently, ${\rm Ker}(d^n)= {\rm Im}(d^{n-1})$. Indeed, set $g=j^{n+1}p_n: T_a^n\to T_a^{n+1}$. Since $d^{n}d^{n-1}=0$, it is not hard to see that $d^n= h g$ for some $\La$-linear morphism $h: T_a^{n+1} \to T_a^{n+1}$.
Write
$$\textstyle g=(g(z,x))_{(z,x)\in Q_0\times Q_0}: \oplus_{x\in Q_0} \,I_x\otimes  e_x \La^!_n e_a\to \oplus_{z\in Q_0} \, I_z \otimes e_z\La^!_{n+1}e_a,$$ where $g(z,x): I_x\otimes e_x \La^!_n e_a\to I_z \otimes e_z\La^!_{n+1}e_a$ is $\La$-linear, and
$$
\textstyle h=(h(y,z))_{(y,z)\in Q_0\times Q_0}: \oplus_{z\in Q_0} \, I_z\otimes e_z \La^!_{n+1} e_a \to \oplus_{y\in Q_0} \, I_y \otimes e_y\La^!_{n+1}e_a,$$
where $h(y,z):  I_z\otimes e_z \La^!_{n+1} e_a\to I_y \otimes e_y\La^!_{n+1}e_a$ is $\La$-linear.

\vspace{1pt}

Given $x, y, z\in Q_0$, consider a basis
$\{\bar\alpha \mid \alpha\in Q_1(z, x)\} \cup \hspace{.6pt}\mathcal{U}_{z,x}$ of $e_x J e_z$, where $\mathcal{U}_{z,x}$ consists of homogeneous elements of degrees $>1$, and a basis $\mathcal{V}_{y,z}$ of $e_z J e_y$ consisting of homogeneous elements. Since $g$ vanishes on $S(T_a)$, by Lemma \ref{Inj-Mor-2}, $g(z,x)=\sum_{\alpha\in Q_1(z,x)}\, I[\bar\alpha] \otimes g_{\alpha}+ \sum_{u\in \mathcal{U}_{z,x}}\, I[u]\otimes g_u$, where $g_{\alpha}, \, g_u$ are $k$-linear maps. Moreover, $h(y,y)=\id_{I_y} \otimes h_{e_y}+\sum_{v\in \mathcal{V}_{y,y}} I[v]\otimes h_v$, where $h_{e_y}, h_v$ are $k$-linear maps, and $h(y,z)=\sum_{v\in \mathcal{V}_{y,z}} I[v]\otimes h_v$ in case $z\ne y$.
Given $(y,x)\in Q_0\otimes Q_0$, by the uniqueness stated in Lemma \ref{Inj-Mor-2},
$d^n(y,x)=\sum_{\alpha\in Q_1(y, x)}\, I[\bar\alpha] \otimes (h_{e_y}\circ g_\alpha).$ Hence, $(h_{e_y}\circ g_\alpha)=P(\alpha^!)$, for all $\alpha\in Q_1(y,x)$. Therefore, we may assume that $h(y,z)=0$ for $z\ne y$, and $h(y,y)=\id_{I_y} \otimes h_{e_y}$. Fix some $y\in Q_0$. Let $w\in e_y \La^!_{n+1} e_a$, say $w=\xi^!$ for some path $\xi\in Q_{n+1}(y, a)$. Writing $\xi=\zeta \alpha$, where $\alpha\in Q_1(y, x)$ and $\zeta\in Q_n(x, a)$ for some $x\in Q_0$, we see that $w=\alpha^! \zeta^!=P(\alpha^!)(\zeta^!)=h_{e_y}(g_\alpha(\xi^!))$. Thus, $h_{e_y}$ is surjective, and since $e_y\La^!_{n+1} e_a$ is finite dimensional, $h_{e_y}$ is bijective. Thus, $h$ is a $\La$-linear isomorphism, and hence, ${\rm Ker}(d^n)={\rm Ker}(g)={\rm Im}(d^{n-1})$. Our claim is established. By induction, $T_a^\pdt$ is a minimal injective co-resolution of $S_a$. The proof of the theorem is completed.

\smallskip

\section{Koszul duality}

\medskip

\noindent The objective of this section is to apply the results obtained in Section 2 to establish the Koszul duality relating the non-graded derived categories of two dual Koszul algebras given by a locally finite gradable quiver. More precisely, we shall
construct a triangle-equivalence between a continuous family of pairs of triangulated subcategories, all but the classic one contain doubly unbounded complexes, of their respective derived categories of all module. In this context, our result extends
the classic Koszul duality stated in \cite[(1.2.6)]{BGS}; see also \cite[Theorem 30]{MOS}. In case the Koszul algebra is left (respectively, right) locally bounded and its Koszul dual is right (respectively, left) locally bounded, our Koszul duality restricts to an equivalence of the bounded derived categories of finitely supported modules (respectively, of finite dimensional modules).

\medskip

Throughout this section, we shall put $\La=kQ/R,$ where $Q$ is a locally finite gradable quiver and $R$ is a quadratic ideal in $kQ$. We fix a grading $Q_0=\cup_{n\in \Z}\, Q^n$ for $Q$. Observe that $Q(x, y)=Q_s(x, y)$, for $x\in Q^n$ and $y\in Q^{n+s}$, where $n, s\in \Z$ with $s\ge 0$; see \cite[(7.2)]{BaL}. In particular, $Q$ is strongly locally finite, and consequently, $\La$ is strongly locally finite dimensional. It will be more convenient for us to identify modules in ${\rm Mod}\La$ with representations in ${\rm Rep}(Q, I)$. Let $(M^\cdt, d^\ydt)$ be a complex over ${\rm Mod}\hspace{.5pt}\La$. For each $x\in Q_0$, we shall consider a complex $(M^\cdt(x), d^\pdt(x))$ over ${\rm Mod}\hspace{.5pt}k$, whose $n$-th component is $M^n(x)$ and whose $n$-th differential is $d^n(x)$. We shall have a similar consideration for double complexes over ${\rm Mod}\La$.

\medskip

Note that $Q^{\rm o}$ admits a grading $(Q^{\rm o})_0=\cup_{n\in \Z}\, (Q^{\rm o})^n$, where $(Q^{\rm o})^n=Q^{-n}$ for $n\in \Z.$ In particular, the quadratic dual $\La^!=kQ^{\rm o}/R^!=\{\gamma^! \mid \gamma\in kQ\}$, where $\gamma^!=\gamma^{\rm o}+R^!$, is also defined by a gradable quiver. Given $a\in Q_0$, the simple $\La^!$-module associated with $a$ will be denoted by $S_a^!$, while its projective cover and injective envelope will be written as $P^!_a$ and $I^!_a$, respectively.

\medskip

We start with defining two {\it Koszul functors}; compare \cite[page 489]{BGS}. Given a module $M \in {\rm Mod}\hspace{.5pt}\La^!$, as shown below, we shall obtain a complex
$F(M)^\cdt\in C({\rm Mod}\hspace{.5pt}\La)$ if, for any $n\in Z$, we define
$F(M)^n=\oplus_{x\in (Q^{\rm o})^n}\, P_x\otimes M(x)=\oplus_{x\in Q^{-n}}\, P_x\otimes M(x)$ \vspace{1pt}
and $d_{F(M)}^n=(d_{F(M)}^n(y,x))_{(y,x) \in Q^{-n-1}\times Q^{-n}}: F(M)^n \to F(M)^{n+1}$, where
$$
d^n_{F(M)}(y,x)={\textstyle\sum}_{\alpha\in Q_1(y,x)}\,P[\bar{\alpha}]\otimes M(\alpha^{\rm o}): P_x\otimes M(x)\to P_y\otimes M(y).
$$
Given a morphism $f: M\to N$ in ${\rm Mod}\hspace{.5pt}\La^!$, we shall obtain a complex morphism $F(f)^\cdt: F(M)^\cdt\to F(N)^\cdt$ if, for any $n\in Z$, we set
$$
F(f)^n=\oplus_{x\in Q^{-n}}\, \id_{P_x}\otimes f(x)\,: \oplus_{x\in Q^{-n}} \, P_x\otimes M(x) \to \oplus_{x\in Q^{-n}} \, P_x\otimes N(x).
$$

On the other hand, given a module $N\in \Mod\hspace{.5pt}\La$, we shall obtain a complex $G(N)^\cdt\in C({\rm Mod}\,\La^!)$ if, for any $n\in Z$, we put
$G(N)^n=\oplus_{x\in Q^n}\, I_x^! \otimes N(x)$ and
$$d_{G(N)}^n=(d_{G(N)}^n(y,x))_{(y,x) \in Q^{n+1}\times Q^n},$$ where
$d^n_{G(N)}(y,x)={\textstyle\sum}_{\alpha:x\to y}\, I[\alpha^!]\otimes N(\alpha): I^!_x \otimes M(x)\to I^!_y\otimes N(y).$ Given a morphism $g: M\to N$ in $\Mod\hspace{.5pt}\La$, we shall have a morphism $G(g)^\cdt: G(M)^\cdt\to G(N)^\cdt$\vspace{1pt} if we define
$
G(g)^n=\oplus_{x\in Q^n}\,\id_{I^!_x}\otimes g(x)\,: \oplus_{x\in Q^n} \, I^!_x\otimes M(x)\to \oplus_{x\in Q^n} \,  I^!_x\otimes N(x).
$

\medskip

\begin{Prop}\label{F-property}

Let $\La=kQ/R$, where $Q$ is a locally finite gradable quiver and $R$ is a quadratic ideal. The above construction yields two fully faithful exact functors

\begin{enumerate}[$(1)$]

\item $F: \Mod\hspace{.5pt}\La^! \to C({\rm Mod}\hspace{.5pt}\La): M\to F(M)^\cdt; f\mapsto F(f)^\cdt;$

\vspace{1pt}

\item $G: \Mod\hspace{.5pt}\La \to C(\Mod\hspace{.5pt}\La^!): N\to G(N)^\cdt; g\mapsto G(g)^\cdt.$

\end{enumerate}

\end{Prop}

\noindent{\it Proof.} We shall prove only Statement (1). Fix a module $M\in \Mod\hspace{.5pt}\La^!$. First, we claim that $F(M)^\cdt$ is a complex. Let $n\in \Z$. Given $z\in Q^{-n-2}$ and $x\in Q^{-n}$, we write $Q(z, x)=\{\alpha_1\beta_1, \ldots, \alpha_s\beta_s\}$, where $\alpha_i, \beta_i\in Q_1$. By definition, we obtain
$$(d^{n+1}_{F(M)}\circ d^n_{F(M)})(z, x)={\textstyle\sum}_{i=1}^s P[\bar{\alpha}_i\bar{\beta}_i]\otimes M(\beta_i^!\alpha_i^!): P_x\otimes M(x) \to P_z\otimes M(z).$$

As seen in the proof of Lemma \ref{perp}, we may find bases  $\{\rho_1, \ldots, \rho_r, \rho_{r+1}, \ldots, \rho_s\}$ and
$\{\eta_1, \ldots, \eta_r, \eta_{r+1}, \ldots, \eta_s\}$ of $kQ(z, x)$ such that $\{\rho_1, \ldots, \rho_r\}$ is a basis of $R_2(z, x)$ and $\{\eta_{r+1}^{\rm o}, \ldots, \eta_s^{\rm o}\}$ is a basis of $R^!_2(x, z)$, while $\{\eta_1^*, \ldots, \eta_r^*, \eta_{r+1}^*, \ldots, \eta_s^*\}$ is the dual basis of $\{\rho_1, \ldots, \rho_r, \rho_{r+1}, \ldots, \rho_s\}$. By Lemma \ref{orthogonal}(2), we obtain
$$\textstyle
\sum_{i=1}^s (\alpha_i\beta_i)\otimes (\alpha_i\beta_i)^*= {\sum}_{i=1}^s \rho_i \otimes \eta_i^* \in kQ(z, x) \otimes D(kQ(z, x)).$$
Considering the $k$-linear isomorphism $D(kQ(z, x))\to kQ^{\rm o}(x, z)$ and the projections $kQ(z, x)\to e_x\La e_z$ and $kQ^{\rm o}(x, z)\to e_z\La^!e_x$, we get $\sum_{i=1}^s\bar{\alpha}_i\bar{\beta}_i\otimes \beta_i^!\alpha_i^! = \sum_{i=1}^s \bar{\rho}_i \otimes \eta_i^!$ \vspace{1.5pt}
Applying the linear map $e_x\La e_z \otimes e_z \La^! e_x \to \Hom_{\it\Lambda}(P_x, P_z) \otimes \Hom_k(M(x), M(z)),$
as seen in Lemma \ref{Mor}(1) and (4), yields ${\textstyle\sum}_{i=1}^s P[\bar{\alpha}_i\bar{\beta}_i]\otimes M(\beta_i^!\alpha_i^!)=
{\textstyle\sum}_{i=1}^s P[\bar{\rho}_i]\otimes M(\eta_i^!),$ \vspace{1.5pt} which vanishes since $\rho_i\in R(z, x)$ for $1\le i\le r$ and $\eta_i^{\hspace{.4pt} \rm o}\in R^!(x, z)$ for $r<i\le s$. Therefore, $d^{n+1}_{F(M)}\circ d^n_{F(M)}=0$. Our claim is established.
It is then easy to verify that $F$ is a functor, and since the tensor product is over $k$, it is exact and faithful.

Let finally $f^\ydt: F(M)^\cdt\to F(N)^\cdt$ be a complex morphism, where $M, N\in {\rm Mod}\hspace{.5pt}\La^!.$ Write
$f^n=(f^n(y,x))_{(y,x)\in Q^{-n}\times Q^{-n}}: \oplus_{x\in Q^{-n}}\, P_x\otimes M(x)\to \oplus_{y\in Q^{-n}}\, P_y\otimes N(y),$ where $f^n(y,x): P_x\otimes M(x)\to P_y\otimes N(y)$ is $\La$-linear.
Since $x, y\in Q^{-n}$, we see that $e_y\La e_x\ne 0$ if and only if $y= x$; and in this case, $e_y\La e_x = ke_x$. By Lemma \ref{rqz-pm}, $f^n(y,x)=0$ if $y\ne x;$ and otherwise, $f^n(y, x)=\id_{P_x}\otimes g(x)$ for some $g(x)\in \Hom_k(M(x), N(x))$, and consequently, we obtain
$$f^n=\oplus_{x\in Q^{-n}}\,\id_{P_x}\otimes g(x): \oplus_{x\in Q^{-n}}\, P_x\otimes M(x)\to \oplus_{x\in Q^{-n}}\, P_x\otimes N(x).$$

Fix $(y, x)\in Q^{-n-1}\times Q^{-n}\vspace{1pt}$ with $e_x\La e_y\ne 0$. Since $f^{n+1} \circ d_{F(M)}^n=d_{F(N)}^n \circ f^n$,
we see that $\sum_{\alpha\in Q_1(y, x)}\, P[\bar{\alpha}]\otimes (g(y) \circ M(\alpha))
=\textstyle{\sum}_{\alpha\in Q_1(y, x)}\, P[\bar{\alpha}]\otimes (N(\alpha) \circ g(x)).$ \vspace{1pt} By the uniqueness stated in Lemma \ref{rqz-pm}, we obtain $g(y) \circ M(\alpha)=N(\alpha) \circ g(x)$, for every $\alpha\in Q_1(y,x)$. Thus the $k$-linear maps $g(x)$, $x\in Q_0$, form a $\La^!$-linear morphism $g: M\to N$ such that $F(g)=f$.
The proof of the proposition is completed.

\medskip

By Proposition \ref{F-extension}, the Koszul functors $F, G$ can be extended to two {\it complex Koszul functors} $F^C:C({\rm Mod}\hspace{.5pt}\La^!)\to C({\rm Mod}\hspace{.5pt}\La)$ and $G^C: C({\rm Mod}\hspace{.5pt}\La)\to C({\rm Mod}\hspace{.5pt}\La^!)$. However, they cannot descend to the derived categories. Thus, we shall consider some subcategories of the complex categories. For this purpose, we need some notation and terminology. Let $M$ be a module in ${\rm Mod}\hspace{.5pt}\La$. The grading $Q=\cup_{n\in \Z} Q^n$ induces a graduation $M=\oplus_{n\in \Z} M_n$ with $M_n=\oplus_{x\in Q^n} M(x)$, which makes $M$ into a graded $\La$-module. Observe that this graduation for $P_a$ with $a\in Q_0$ is different from its usual $J$-radical graduation. A module $M$ is called {\it bounded-above} if $M_n=0$ for almost all $n>0$, and {\it bounded-below} if $M_n=0$ for almost all $n<0$. We shall denote by $\Mod^-\hspace{-2pt}\La$ and $\Mod^+\hspace{-3pt}\La$ the full subcategories of $\Mod \La$ generated by the bounded-above modules and by the bounded-below modules, respectively. 
Keeping this graduation in mind, we can view a complex $M^\cdt$ over ${\rm Mod}\hspace{.5pt}\La$ as a bigraded $k$-module $M^i_j=\oplus_{x\in Q^j}\, M(x)$, with $i, j\in \Z$.

\medskip

\begin{Defn} Let $\La=kQ/R$, where $Q$ is a locally finite gradable quiver and $R$ is a quadratic ideal in $kQ$. Given $p, q\in \mathbb{R}$ with $p\ge 1$ and $q\ge 0$, we denote by

\begin{enumerate}[$(1)$]

\item $C^{\,\downarrow}_{p,q}({\rm Mod}\hspace{.5pt}\La)$ \vspace{1pt} the full abelian subcategory of $C({\rm Mod}\hspace{.5pt}\La)$ of the complexes $M^\cdt$ with $M^i_j=0$ for $i+pj>>0$ or $i-qj<< 0\,;$ in other words, $M^\cdt$ concentrates in the lower triangle formed by two lines of slopes $-\frac{1}{p}$ and $-\frac{1}{q}$ respectively;

\vspace{2pt}

\item $C^{\,\uparrow}_{p,q}({\rm Mod}\hspace{.5pt}\La)$ \vspace{1pt} the full abelian subcategory of $C({\rm Mod}\hspace{.5pt}\La)$ of the complexes $M^\cdt$ with $M^i_j=0$ for $i+pj<<0$ or $i-qj>>0;$ in other words, $M^\cdt$ concentrates in the upper triangle formed by two lines of slopes $-\frac{1}{p}$ and $-\frac{1}{q}$ respectively.

\end{enumerate}

\end{Defn}

\medskip

\noindent{\sc Remark.} (1) Taking $p=1$ and $q=0$, we obtain the categories $C^\downarrow({\rm Mod}\hspace{.5pt}\La)$ and $C^\uparrow({\rm Mod}\hspace{.5pt}\La)$ as defined in \cite[(2.12)]{BGS}.

\vspace{1pt}

\noindent (2) It is easy to see that the categories $C^{\,\downarrow}_{p,q}({\rm Mod}\hspace{.5pt}\La)$ and $C^{\,\uparrow}_{p,q}({\rm Mod}\hspace{.5pt}\La)$ are subcategories of $C({\rm Mod}^-\hspace{-3pt}\La)$ and $C({\rm Mod}^+\hspace{-3pt}\La)$, respectively.

\medskip

We shall denote by $K^{\,\downarrow}_{p,q}({\rm Mod}\hspace{.5pt}\La)$ the quotient of $C^{\,\downarrow}_{p,q}({\rm Mod}\hspace{.5pt}\La)$ \vspace{1pt} modulo the null-homotopic morphisms, and by $D^{\,\downarrow}_{p,q}({\rm Mod}\hspace{.5pt}\La)$ the localization of $K^{\,\downarrow}_{p,q}({\rm Mod}\hspace{.5pt}\La)$ with respect to the quasi-isomorphisms. Similarly, we shall have the homotopy category $K^{\,\uparrow}_{p,q}({\rm Mod}\hspace{.5pt}\La)$ \vspace{1pt} and the derived category $D^{\,\uparrow}_{p,q}({\rm Mod}\hspace{.5pt}\La)$. The following statement extends Proposition 20 in \cite{MOS} under our gradable setting.

\medskip

\begin{Theo}\label{F-diag}

Let $\La=kQ/R$, where $Q$ is a locally finite gradable quiver and $R$ is a quadratic ideal in $kQ$. If $p, q\in \mathbb{R}$ with $p\ge 1$ and $q\ge 0$, \vspace{1pt} then the Koszul functor $F\hspace{-1pt}: \Mod\hspace{.5pt}\La^! \to C(\Mod\hspace{.5pt}\La)$ induces a commutative diagram
$$\xymatrix{
C^{\,\downarrow}_{p,q}({\rm Mod}\hspace{.4pt}\La^!) \ar[d]_{F^C} \ar[r] & K^{\,\downarrow}_{p,q}({\rm Mod}\hspace{.4pt}\La^!)\ar[d]^{F^K} \ar[r] & D^{\,\downarrow}_{p,q}({\rm Mod}\hspace{.4pt}\La^!) \ar[d]^{F^D} \\
C^{\,\uparrow}_{q+1, p-1}({\rm Mod}\hspace{.4pt}\La)\ar[r] & K^{\,\uparrow}_{q+1, p-1}({\rm Mod}\hspace{.4pt}\La) \ar[r]  & D^{\,\uparrow}_{q+1, p-1}({\rm Mod}\hspace{.4pt}\La);
}\vspace{-3pt}$$  while the Koszul functor $G\hspace{-1pt}: \Mod\hspace{.5pt}\La \to C(\Mod\hspace{.5pt}\La^!)$ induces a commutative diagram
$$\xymatrix{
C^{\,\uparrow}_{p,q}({\rm Mod}\hspace{.4pt}\La)\ar[r] \ar[d]_{G^C} & K^{\,\uparrow}_{p,q}({\rm Mod}\hspace{.4pt}\La) \ar[r] \ar[d]^{G^K}  & D^{\,\uparrow}_{p,q}({\rm Mod}\hspace{.4pt}\La)\ar[d]^{G^D}\\
C^{\,\downarrow}_{q+1, p-1}({\rm Mod}\hspace{.4pt}\La^!)  \ar[r] & K^{\,\downarrow}_{q+1, p-1}({\rm Mod}\hspace{.4pt}\La^!)\ar[r] & D^{\,\downarrow}_{q+1, p-1}({\rm Mod}\hspace{.4pt}\La^!),}$$
where $F^D$ and $G^D$ are triangle-exact, called \em{derived Koszul functors}.

\end{Theo}

\noindent{\it Proof.} Let $p, q\in \mathbb{R}$ with $p\ge 1$ and $q\ge 0$. Consider $M^\ydt\in C^{\,\downarrow}_{p,q}({\rm Mod}\hspace{.5pt}\La^!)$. There exist $s, t\in \Z$ such that $M^i(x)=0$, for $x\in (Q^{\rm o})^j$ with $i+pj > s$ or $i-qj < t.$
Note that $F^C(M^\cdt)$ is the total complex of $F(M^\cdt)^\cdt$. Fix $m, n\in \Z$. Given $y\in Q^m$, we obtain $F^C(M^\ydt)^n(y)=\oplus_{i\in \Z; \, x\in (Q^{\rm o})^{n-i}}\, P_x(y)\otimes M^i(x)= \oplus_{i\le n+m; \, x\in Q^{i-n}}\, P_x(y)\otimes M^i(x).\vspace{1pt}$ Fix an integer $i\le n+m$. If $n+(q+1)m<s$, then $i-q(n-i)<s$; and if $n-(p-1)m>t$, then $i+p(n-i)>t$. In either case, $M^i(x)=0$ for all $x\in (Q^{\rm o})^{n-i}$, and thus, $F^C(M^\ydt)^n(y)=0$. That is, $F^C(M^\ydt)\in C^{\,\uparrow}_{p,q}({\rm Mod}\hspace{.4pt}\La).$ This yields a functor $F^C: C^{\,\downarrow}_{p,q}({\rm Mod}\hspace{.4pt}\La^!)\to C^{\,\uparrow}_{q+1, p-1}({\rm Mod}\hspace{.4pt}\La).$ By Proposition \ref{F-extension}, $F^C$ induces a triangle-exact functor $F^K: K^{\,\downarrow}_{p,q}({\rm Mod}\hspace{.4pt}\La^!)\to K^{\,\uparrow}_{q+1, p-1}({\rm Mod}\hspace{.4pt}\La)$ making the left square of the first diagram commute.

Assume next that $M^\cdt$ is acyclic. Since $F$ is exact, $F(M^\cdt)^\cdt$ has exact rows. Moreover, $F(M^i)^{n-i}=\oplus_{x\in (Q^{\rm o}){n-i}}\, P_x\otimes M^i(x)$, with $i\in \Z$, form the $n$-diagonal of $F(M^\cdt)^\cdt$. By the assumption on $M^\cdt$, we see that $M^i(x)=0$ for all $x\in (Q^{\rm o})^{n-i}$ with $i< (nq +t)(1+q)^{-1}$. That is, $F(M^\cdt)^\cdt$ is diagonally bounded-below. By Corollary \ref{AAL}, $F^C(M^\cdt)$ is acyclic. Thus, $F^K$ sends quasi-isomorphisms to quasi-isomorphisms, and then, it induces a triangle-exact functor $F^D: D^{\,\downarrow}_{p,q}({\rm Mod}\hspace{.4pt}\La^!)\to D^{\,\uparrow}_{q+1, p-1}({\rm Mod}\hspace{.4pt}\La)$ making the right square of first diagram commute.

\vspace{1pt}

Most of the arguments for the second part of the theorem will be similar to those used above. We shall only show, for an acyclic complex $N^\cdt\in C({\rm Mod}\hspace{.5pt}\La)$, that $G^C(N^\cdt)$ is acyclic. It suffices to show that
$G^C(N^\cdt)(y)$ is acyclic, for any $y\in Q^s$ with $s\in \Z$. Since $G^C(N^\cdt)$ is the total complex of $G(N^\cdt)^\cdt$, we see that $G^C(N^\cdt)(y)$ is the total complex of $G(N^\ydt)^\ydt(y)$. Since $G$ is exact, $G(N^\cdt)^\cdt$ has exact rows, and so does $G(N^\ydt)^\ydt(y)$. Given $n\in \Z$, the modules $G(N^i)^{n-i}(y)=\oplus_{x\in Q^{n-i}}\, I^!_x(y)\otimes N^i(x)$ with $i\in \Z$ form the $n$-diagonal of $G(N^\ydt)^\ydt(y)$. If $i<n-s$ and $x\in Q^{n-i}$, then $Q$ has no path from $x$ to $y$, and hence, $I^!_x(y)=D(e_y (\La^!)^{\rm o} e_x)=0$. Thus, $G(N^i)^{n-i}(y)=0$ for all $i<n-p$. That is, $G(N^\cdt)^\cdt(y)$ is diagonally bounded-below. By Corollary \ref{AAL}, $G^C\hspace{-1.5pt}(M^\cdt)(y)$ is acyclic. The proof of the theorem is completed.

\medskip

Next, we shall show that the two derived Koszul functors are mutually quasi-inverse in case $\La$ is Koszul. For this purpose, given a simple $\La$-module $S$, we shall denote by ${\mathcal P}_S^\pdt$ its minimal projective resolution and by ${\mathcal I}_S^\pdt$ its minimal injective co-resolution, which can be explicitly described as below; compare \cite[(1.2.6)]{BGS} 

\medskip

\begin{Lemma}\label{inj-im}

Let $\La=kQ/R$ be a Koszul algebra, where $Q$ is locally finite with a grading $Q_0=\cup_{n\in \Z}\hspace{.5pt} Q^n$. If $a\in Q^{\hspace{.5pt}s}$, then $F(I_a^!)^\cdt\cong \mathcal{P}_{S_a}^\ydt\hspace{-1pt}[s]$ and $G(P_a)^\cdt\cong \mathcal{I}_{S^!_a}^\ydt\hspace{-2pt}[-s].$

\end{Lemma}

\noindent{\it Proof.} Fix $a\in Q^s$. By Theorem \ref{Koz-proj-rls} and Lemma \ref{k-cplx-iso}, ${\mathcal P}_{S_a}^\pdt$ is isomorphic to
$$L^\cdt: \xymatrixrowsep{18pt}\xymatrix{\cdots \ar[r] & L^{-i} \ar[r]^-{d^{-i}} & L^{-i+1} \ar[r] &
\cdots \ar[r] & L^{-1} \ar[r]^{d^{-1}} & L^0 \ar[r] & 0\ar[r] & \cdots,}$$ where $L^{-i}=\oplus_{x\in Q_0} P_x \otimes D(e_a\La^!_ie_x)$ for $i\ge 0,$ and $d^{-i}=(d^{-i}(y,x))_{(y,x)\in Q_0\times Q_0}\vspace{1pt}$ with $d^{-i}(y,x)={\textstyle\sum}_{\alpha\in Q_1(y,x)} P[\bar{\alpha}]\otimes DP[\alpha^!]: P_x \otimes  D(e_a\La^!_ie_x)\to P_y \otimes  D(e_a\La^!_{i-1}e_y)$ for $i>0.$ Fix an integer $n\ge 0$. Since $Q$ is gradable, $e_a \La^!_n e_x=0$ in case $x\not\in Q^{n+s}$; and otherwise, $e_a \La^!_n e_x= e_a \La^! e_x$. Thus,
$L^{-n}=\oplus_{x\in Q^{n+s}} P_x \otimes D(e_a \La^! e_x).$ Note that the $k$-linear isomorphism $e_x (\La^!)^{\rm o} e_a\to e_a \La^! e_x$ induces a $k$-linear isomorphism $\theta_{a,x}: D(e_a \La^! e_x)\to D(e_x (\La^!)^{\rm o} e_a)=I_a^!(x)$ such that the diagram
$$\xymatrix{
\oplus_{x\in Q^{n+s}} P_x \otimes D(e_a \La^! e_x)\ar[d]_{\oplus \, (\id\otimes \theta_{a,x})}  \ar[rrrr]^----{\textstyle\sum_{\alpha\in Q_1(y,x) }P[\bar\alpha]\otimes DP[\alpha^!]} &&&&
\oplus_{y\in Q^{n+s-1}} P_y \otimes D(e_a \La^! e_y)\ar[d]^{\oplus \,(\id\otimes \theta_{a,y})}   \\
\oplus_{x\in Q^{n+s}} P_x \otimes I^!_a(x)
\ar[rrrr]^----{\textstyle\sum_{\alpha\in Q_1(y,x)} P[\bar\alpha]\otimes \hspace{-2pt} I_a^!(\alpha^{\rm o})} &&&&
\oplus_{y\in Q^{n+s-1}} P_y \otimes I^!_a(y)
}$$ is commutative with vertical isomorphisms. Since $F(I_a^!)^{-n-s}=0$ for $n<0$, we conclude that $L^\cdt\cong \mathfrak{t}^s(F(I_a^!)[-s])\cong F(I_a^!)[-s]$. This establishes the first part of the lemma. Next, observe that
$\La^!$ is Koszul with $(\La^!)^!=\La$; see (\ref{Dual-Koszul}) and (\ref{q-dual}). In view of Theorem \ref{Opp-Koszul}(3), we see that $\mathcal{I}_{S_a^!}^\pdt$ is isomorphic to
$$T^\pdt: \quad \xymatrix{ 0 \ar[r] & T^0 \ar[r]^{d^0} & T^1 \ar[r] & \cdots \ar[r] & T^n \ar[r]^{d^n} & T^{n+1} \ar[r] & \cdots
}$$
where $T^n=\oplus_{x\in Q_0} \,I_x^! \otimes e_x \La_n e_a $ and $d^n=(d^n(y,x))_{(y,x)\in Q_0\times Q_0}: T^n\to T^{n+1}$ with
$d^n(y,x)=\textstyle\sum_{\alpha\in Q_1(x,y)} I[\alpha^!]\otimes P_a(\alpha): I_x^! \otimes  e_x\La_ne_a\to I_y^! \otimes  e_y\La_{n+1}e_a,$\vspace{1pt} for $n\ge 0$.
Fix an integer $n\ge 0$. Since $Q$ is gradable, $e_x \La_n e_a=0$ in case $x\not\in Q^{n+s}$; and otherwise, $e_x\La_n e_a = e_x\La e_a.$ Thus,
$T^n= \oplus_{x\in Q^{n+s}}\, I_x^! \otimes e_x\La e_a=G(P_a)^{n+s}$ and \vspace{1pt} $d^n=d^{n+s}_{G(P_a)}$, for $n\ge 0$. Since $G(P_a)^{n+s}=0$ for $n<0$, we conclude that $\mathcal{I}_{S_a^!}^\pdt\cong \mathfrak{t}^s(G(P_a)^\pdt[s])\cong G(P_a)^\pdt[s].$
The proof of the lemma is completed.

\medskip

More generally, we can obtain an explicit projective resolution for any module.

\medskip

\begin{Prop}\label{natural-tran}

Let $\La=kQ/R$ be Koszul, with $Q$ a locally finite gradable quiver. Given $M\in {\rm Mod}\hspace{.4pt}\La$, we have a natural quasi-isomorphism $\eta_{_M}^\ydt\hspace{-2pt}: (F^C\circ G)(M)^\ydt \to M$.

\end{Prop}

\noindent{\it Proof.} Let $M\in {\rm Mod}\hspace{.4pt}\La$. By definition, $(F^C\circ G)(M)^\ydt=F^C(G(M)^\cdt)$, which is the total complex of the double complex $F(G(M)^\ydt)^\ydt$. Given $i, n\in \Z$, we obtain $G(M)^i=\oplus_{x\in Q^i}\, I_x^!\otimes M(x)$ and
$$(F^C \hspace{-2.5pt} \circ \hspace{-1.5pt} G)(M)^n=\oplus_{i\in \mathbb{Z}; a\in Q^{i-n}} P_a\otimes G(M)^i(a)=\oplus_{i\in \mathbb{Z}; \hspace{.7pt} a\in Q^{i-n}; \hspace{.7pt} x\in Q^i}\hspace{.5pt} P_a\otimes I_x^!(a)\otimes M(x).$$

Suppose that $n>0$. Given any $(a, x)\in  Q^{i-n}\times Q^i$ with $i\in \Z$, since $Q$ has no path from $x$ to $a$, we obtain $I_x^!(a)=0$. Therefore, $(F^C \hspace{-2.5pt} \circ \hspace{-1.5pt} G)(M)^n=0$.

\vspace{2pt}

Suppose that $n<0.$ We claim that ${\rm H}^n((F^C \hspace{-2.5pt} \circ \hspace{-1.5pt} G)(M)^\ydt)=0$, or equivalently, ${\rm H}^n((F^C \hspace{-2pt} \circ \hspace{-1.5pt} G)(M)^\ydt(y))=0,$ for any $y\in Q^p$ with $p\in \Z.$ Now, $(F^C \hspace{-1.5pt} \circ \hspace{-.5pt} G)(M)^\ydt(y)$ is the total complex of $F(G(M)^\ydt)^\ydt(y)$, whose $n$-diagonal consists of
the modules $F(G(M)^i)^{n-i}(y)=\oplus_{a\in Q^{i-n}; \hspace{.7pt} x\in Q^i} P_a(y)\otimes I_x^!(a)\otimes M(x),$ $i\in \Z.$ \vspace{1pt} If $i>n+p$, then $P_a(y)=0$ for all $a\in Q^{i-n}$. Hence, $F(G(M)^i)^{n-i}(y)=0$. That is, $F(G(M)^\ydt)^\ydt(y)$ is $n$-diagonally bounded-above. For any $i\in \Z$, the $i$-th column of $F(G(M)^\ydt)^\ydt$ is $\frak{t}^i(F(G(M)^i)^\ydt)=\oplus_{x\in Q^i}\, \frak{t}^i(F(I_x^!)^\cdt)\otimes M(x),$ where $F(I_x^!)^\ydt\cong {\mathcal P}_{S_x}^\ydt[i]$; see (\ref{inj-im}). Thus,

\vspace{1pt}
\noindent ${\rm H}^{n-i}(\frak{t}^i(F(G(M)^i)^\ydt))\cong \oplus_{x\in Q^i} {\rm H}^{n-i}({\mathcal P}_{S_x}^\ydt[i]) \otimes M(x)=\oplus_{x\in Q^i} {\rm H}^n({\mathcal P}_{S_x}^\ydt) \otimes M(x)=0,\vspace{1pt} $ and consequently,  ${\rm H}^{n-i}(\frak{t}^i(F(G(M)^i)^\ydt)(y))=0.$ Since $\frak{t}^i(F(G(M)^i)^\ydt)(y)$ is the $i$-column of $F(G(M)^\ydt)^\ydt(y)$, by Lemma \ref{Homology-zero}(2), we see that ${\rm H}^n((F^C \hspace{-2.5pt} \circ \hspace{-1.5pt} G)(M)^\ydt(y))=0$.

\vspace{1pt}

It remains to show that ${\rm H}^0((F^C \hspace{-3pt} \circ \hspace{-1.5pt} G)(M)^\ydt)$ is naturally isomorphic to $M$.
For this purpose, observing that the $1$-diagonal of the double complex $F(G(M)^\ydt)^\ydt$ is zero, we illustrate its $(-1)$-diagonal and $0$-diagonal as follows:
$$\xymatrixcolsep{35pt}\xymatrix{
\oplus_{b\,\in Q^i}\, P_b\otimes I_b^!(b)\otimes M(b)\\
\oplus_{(a, x)\,\in Q^{i+1}\times Q^i}\, P_a\otimes I_x^!(a)\otimes M(x) \ar[u]^-{v^{i, -i-1}} \ar[r]^-{h^{i,-i-1}} &
\oplus_{c\, \in Q^{i+1}}  P_{c}\otimes I_c^!(c)\otimes M(c),}$$
where $v^{i,-i-1}=(v^{i,-i-1}(b, a,x))_{(b, a, x)\in Q^i \times  Q^{i+1} \times Q^i},\vspace{1.5pt}$ with
$v^{i,-i-1}(b,a, x)=0$ in case $b\ne x$, and otherwise, $v^{i,-i-1}(x,a,x)={\textstyle\sum}_{\alpha\in Q_1(x, a)}(-1)^i P[\bar{\alpha}]\otimes I_x^!(\alpha^{\rm o})\otimes \id_{M(x)};$ and $h^{i,-i-1}=(h^{i,-i-1}(c,a, x))_{(c,a, x)\in Q^{i+1} \times Q^{i+1}\times Q^i}$, \vspace{1pt} with $h^{i,-i-1}(c,a, x)=0$ in case $c\ne a$, and otherwise, $h^{i,-i-1}(a,a,x)= {\textstyle\sum}_{\alpha\in Q_1(x, a)} \id_{P_a} \otimes I[\alpha^!]_a\otimes M(\alpha).$

Recall that $(\La^!)^{\rm o}=kQ/(R^!)^{\rm o}=\{\hat{\gamma\hspace{2pt}} \mid \gamma\in kQ\}$, where $\hat{\gamma\hspace{2pt}}=\gamma+(R^!)^{\rm o}$. Given $(x, y)\in Q^i\times Q^{i+1}$ with $i\in \Z$, in view of Lemma \ref{Special-func}, $I_x^!(x)$ has a $k$-basis $\{\hat{e}_x^\star\}$ and $I_x^!(y)$ has a $k$-basis $\{\hat\alpha^\star \mid \alpha\in Q_1(x, y)\}$.

\vspace{3pt}

{\sc Sublemma.} {\it Let $d^{-1}$ be the differential of degree $-1$ of $(F^C \hspace{-2.5pt} \circ \hspace{-1.5pt} G)(M)^\ydt$. Consider
$(x, a)\in Q^i\times Q^{i+1}$ for some $i\in \Z$. If $\bar{\gamma\hspace{1.6pt}}\hspace{-1pt}\in P_a$, $\beta\in Q_1(x, a)$ and $u\in M(x)$, then
$d^{-1}(\bar{\gamma\hspace{1.6pt}}\hspace{-1pt} \otimes \hat{\beta}^\star \otimes u)=(-1)^i \,\bar{\gamma\hspace{1.6pt}}\hspace{-1pt} \hspace{.4pt} \bar\beta \otimes \hat{e}_x^\star \otimes u +  \bar{\gamma\hspace{1.6pt}}\hspace{-1pt} \otimes \hat{e}_a^\star \otimes \bar\beta u.$}

{\it Proof.} Given $\alpha\in Q_1(x, a)$, we see that $I_x^!(\alpha^{\rm o})(\hat{\beta}^\star)=0$ if $\alpha\ne \beta$, and otherwise, $I_x^!(\alpha^{\rm o})(\hat{\beta}^\star)=\hat{e}_x^\star$. On the other hand, $I[\alpha^!]_a(\hat{\beta}^\star)=0$ if $\alpha\ne \beta$, and otherwise, $I[\alpha^!]_a(\hat{\beta}^*)=\hat{e}_a^\star.$ This yields

$$\begin{array}{rcl}
d^{-1} (\bar{\gamma\hspace{1.6pt}}\hspace{-1pt}\otimes \hat{\beta}^\star \otimes u)
                 & = & (-1)^{i}{\textstyle\sum}_{\alpha\in Q_1(x, a)} (P[\bar\alpha]\otimes I_x^!(\alpha^{\rm o})\otimes \id_{M(x)}) (\bar{\gamma\hspace{1.6pt}}\hspace{-1pt} \otimes \hat{\beta}^\star\otimes u ) \vspace{2.5pt}\\
                 &   & + {\textstyle\sum}_{\alpha\in Q_1(x, a)}\, (\id_{P_a} \otimes I[\alpha^!]_a\otimes M(\alpha))
                 (\bar{\gamma\hspace{1.6pt}}\hspace{-1pt} \otimes \hat{\beta}^\star\otimes u )  \vspace{2.5pt}\\
                 & = & (-1)^i \hspace{.4pt} \bar{\gamma\hspace{1.6pt}}\hspace{-1pt} \hspace{.4pt} \bar\beta \otimes \hat{e}_x^\star\otimes u +  \gamma \otimes \hat{e}_a^\star \otimes \bar\beta u.
\end{array}$$

Now, it is easy to see that we have a natural $\La$-linear map
$$\eta^0_{_M}: (F^C \hspace{-3pt} \circ \hspace{-1.3pt} G)(M)^0 \to M: {\textstyle\sum}_{{(i,x)\in \Z\times Q^i}}\, \bar{\gamma\hspace{1.6pt}}\hspace{-1pt}_x\otimes \hat{e}^\star_x\otimes u_x\mapsto {\textstyle\sum}_{(i,x)\in \Z\times Q^i}\, (-1)^{\frac{i(i+1)}{2}} \bar{\gamma\hspace{1.6pt}}\hspace{-1pt}_x u_x,$$
where $\bar{\gamma\hspace{1.6pt}}\hspace{-1pt}_x\in P_x$ and $u_x\in M(x)$. We claim that $\eta_{_M}^0 \hspace{-2.5pt} \circ \hspace{-1.5pt} d^{-1}=0$. Indeed, consider an element $\omega\in (F^C \hspace{-3pt} \circ \hspace{-1.3pt} G)(M)^{-1}$. We may assume that $\omega \in P_a\otimes I_x^!(a)\otimes M(x),$ for some $(a, x)\in Q^{i+1} \times Q^i$ with $i\in \Z$. In this case, we may assume further that $\omega=\bar{\gamma\hspace{1.6pt}}\hspace{-1pt} \otimes \hat\beta^\star \otimes u,$ where $\bar{\gamma\hspace{1.6pt}}\hspace{-1pt}\in P_a$, $\beta\in Q_1(x, a)$, and $u\in M(x)$. In view of the sublemma, we obtain
$$
\begin{array}{rcl}
(\eta^0_{_M} \hspace{-2pt} \circ \hspace{-1pt} d^{-1}) (\omega)
& = & \eta_{_M}^0 \left( (-1)^i \hspace{.4pt} \bar{\gamma\hspace{1.6pt}}\hspace{-1pt} \hspace{.4pt} \bar\beta \otimes \hat{e}_x^\star \otimes u +  \bar{\gamma\hspace{1.6pt}}\hspace{-1pt} \otimes \hat{e}_a^\star \otimes \bar\beta u \right) \vspace{3pt} \\
& = & (-1)^{\frac{i(i+1)}{2}+i} (\bar{\gamma\hspace{1.6pt}}\hspace{-1pt} \hspace{.4pt} \bar\beta u) +  (-1)^{\frac{(i+1)(i+2)}{2}} (\bar{\gamma\hspace{1.6pt}}\hspace{-1pt} \bar\beta u) \vspace{3pt} \\
& = & 0.
\end{array} $$

Given $\omega\in {\rm Ker}(\eta_{_M}^0)$, we shall define an integer $n_\omega$ as follows. If $\omega=0$, set $n_\omega=0;$ and in this case, $\omega\in {\rm Im}(d^{-1})$. Otherwise, let $n_\omega$ be minimal for which $\omega=\sum_{i=1}^s \bar{\gamma\hspace{1.6pt}}\hspace{-1pt}_i \otimes \hat{e}_{x_i}^\star\otimes u_i$, where $x_i\in Q_0\hspace{.6pt};$ $\gamma_i\in kQ_{\le n_\omega}(x_i, -)\hspace{.6pt};$ the $u_i$ are linearly independent in $M(x_i)$. For $1\le i\le s$, write $\gamma_i=\lambda_i \varepsilon_{x_i}+ \sigma_{i1} \alpha_{i1} + \cdots+ \sigma_{i,t_i} \alpha_{i,t_i}$, where $\lambda_i\in k \hspace{.6pt};$ $\alpha_{ij}\in Q_1(x_i, a_{ij})\hspace{.6pt};$
$\sigma_{ij}\in kQ_{\le n_\omega-1}(a_{ij}, -)$. \vspace{1pt} Setting  $|x|=i$ for $x\in Q^i$, we obtain ${\textstyle\sum}_{i=1}^s\, (-1)^{\frac{|x_i|(|x_i|+1)}{2}}\bar{\gamma\hspace{1.6pt}}\hspace{-1pt}_i u_i=0\vspace{1.5pt}$. Then, $\sum_{i=1}^s \lambda_iu_i=0$, and hence, $\lambda_i=0$, that is, $\gamma_i=\sigma_{i1} \alpha_{i1} + \cdots+ \sigma_{i,t_i} \alpha_{i,t_i}$, for $i=1, \ldots, s.$ Setting
$$\sigma={\textstyle\sum}_{i=1}^s {\textstyle\sum}_{j=1}^{t_i}\, (-1)^{|x_i|} \, \bar\sigma_{i,j} \otimes \hat\alpha_{ij}^\star \otimes u_i\in (F^C \hspace{-3pt} \circ \hspace{-1.3pt} G)(M)^{-1},$$
we deduce from the sublemma that
$$\begin{array}{rcl}
d^{-1}(\sigma) &=& {\textstyle\sum}_{i=1}^s {\textstyle\sum}_{j=1}^{t_i}\, \left(\bar\sigma_{ij} \, \bar\alpha_{ij} \otimes \hat{e}_{x_i}^\star \otimes u_i + (-1)^{|x_i|} \, \bar\sigma_{ij} \otimes \hat{e}_{a_{ij}}^\star \otimes \bar\alpha_{ij} \, u_i\right) \vspace{3pt}\\
&=& {\textstyle\sum}_{i=1}^s \left( \bar{\gamma\hspace{1.6pt}}\hspace{-1pt}_i \otimes \hat{e}_{x_i}^\star \otimes u_i + {\textstyle\sum}_{j=1}^{t_i} (-1)^{|x_i|} \, \bar\sigma_{ij} \otimes \hat{e}_{a_{ij}}^\star \otimes \bar\alpha_{ij} \, u_i\right)\vspace{3pt}\\
&=& \omega+\omega',
\end{array}$$
where $\omega' =  {\textstyle\sum}_{i=1}^s {\textstyle\sum}_{j=1}^{t_i} (-1)^{|x_i|}\,\bar\sigma_{ij} \otimes \hat{e}_{a_{ij}}^\star\otimes \bar\alpha_{ij} \, u_i$. By definition, $n_{\omega'}<n_\omega$, and
$$\eta_{_M}^0(\omega')= {\textstyle\sum}_{i=1}^s  {\textstyle\sum}_{j=1}^{t_i}(-1)^{|x_i|+\frac{|x_i|(x_i|+1)}{2}} \bar\sigma_{ij} \bar{\alpha}_{ij}\,u_i=- {\textstyle\sum}_{i=1}^s (-1)^{\frac{|x_i|(x_i|+1)}{2}} \bar{\gamma\hspace{1.6pt}}\hspace{-1pt}_i\,u_i=0.$$
By induction, $\omega\in {\rm Im}(d^{-1})$. Thus, ${\rm Im}(d^{-1})={\rm Ker}(\eta_{_M}^0)$. This yields a natural quasi-isomorphism $\eta_{_M}^\cdt: (F\circ G)(M)^\cdt\to M.$ The proof of the proposition is completed.

\medskip

For the injective co-resolution, our result will be slightly more restrictive.

\medskip

\begin{Prop}\label{natural-tran-2}

Let $\La=kQ/R$ be Koszul, with $Q$ a locally finite gradable quiver. Given $N\in {\rm Mod}^-\hspace{-3pt}\La^!$, we have a natural quasi-isomorphism $\theta_{\hspace{-1pt}{_N}}^{\hspace{.7pt}\cdt}: N \to (G^C\circ F)(N)^\ydt$.

\end{Prop}

\noindent{\it Proof.} Fix $N\in {\rm Mod}^-\hspace{-3pt}\La^!$. Let $r$ be such that $N(a)=0$ for all $a\in (Q^{\rm o})^i=Q^{-i}$ with $i>r$. By definition, $(G^C \hspace{-1.5pt} \circ \hspace{-.5pt} F)(N)^\ydt$ is the total complex of the double complex $G(F(N)^\ydt)^\ydt$. We split our proof into several statements.

\vspace{2pt}

{\sc Statement 1.} {\it Given $i\in \Z$, the $i$-th column of the double complex $G(F(N)^\ydt)^\ydt$ is $\frak{t}^i(G(F(N)^i)^\ydt)=\oplus_{a\in Q^{-i}}\, \frak{t}^i(G(P_a)^\cdt)\otimes N(a)\cong \oplus_{a\in Q^{-i}}\, \frak{t}^i({\mathcal I}_{S^!_a}^\ydt[i])\otimes N(a).$}

\vspace{1pt}

Indeed, $F(N)^i=\oplus_{a\in Q^{-i}}\, P_a\otimes N(a)$, and in view of by Lemma \ref{inj-im}, we see that $G(P_a)^\cdt\cong {\mathcal I}_{S^!_a}^\ydt[i]$ for all $a\in Q^{-i}$. 

\vspace{2pt}

{\sc Statement 2.} {\it Given $n\in \Z$, we obtain $(G^C \hspace{-2.5pt} \circ \hspace{-1.5pt} F)(N)^n=0$ in case $n<0;$ and ${\rm H}^n((G^C \hspace{-2.5pt} \circ \hspace{-1.5pt} F)(N)^\ydt)=0$ in case $n>0$.}

\vspace{1pt}

Indeed, given any $n\in \Z$, we obtain $(G^C \hspace{-2.5pt} \circ \hspace{-1.5pt} F)(N)^n=\oplus_{i\in \Z} \hspace{.8pt} G(F(N)^i)^{n-i}$, where
$$G(F(N)^i)^{n-i}=\oplus_{x\in Q^{n-i}} I^!_x\otimes F(N)^i(x)=\oplus_{\hspace{.7pt} x\in Q^{n-i}; \hspace{.7pt} a\in Q^{-i}}\; I^!_x\otimes P_a(x)\otimes N(a).$$ If $n<0$, then $P_a(x)=0$ for all $(x, a)\in  Q^{n-i}\times Q^{-i}$ with $i\in \Z$, and therefore, $(G^C \hspace{-2.5pt} \circ \hspace{-1.5pt} F)(N)^n=0$. Suppose that $n>0.$ Since $N(a)=0$ for every $a\in Q^{-i}$ with $i>r$, we see that $G(F(N)^\ydt)^\ydt$ is $n$-diagonally bounded-above. Moreover, in view of Statement 1, the $(n-i)$-th cohomology of the $i$-th column of $G(F(N)^\ydt)^\ydt$ is given by
$${\rm H}^{n-i}(\frak{t}^i(G(F(N)^i)^\ydt))\cong \oplus_{a\in Q^{-i}} {\rm H}^{n-i}({\mathcal I}_{S^!_a}^\pdt[i]) \otimes N(a)=\oplus_{a\in Q^{-i}} {\rm H}^n({\mathcal I}_{S^!_a}^\pdt) \otimes N(a)=0. $$
By Lemma \ref{Homology-zero}(2), ${\rm H}^n((G^C \hspace{-2.5pt} \circ \hspace{-1.5pt} F)(N)^\ydt)=0$.

\vspace{2pt}

It remains to construct a natural isomorphism $N\to {\rm H}^0((G^C \hspace{-3pt} \circ \hspace{-1.5pt} F)(N)^\ydt)$. Observe that
the modules $G(F(N)^i)^{-i}=\oplus_{\hspace{.7pt} x\in Q^{n-i}; \hspace{.7pt} a\in Q^{-i}}\; I^!_x\otimes P_a(x)\otimes N(a),$ with $i\in \Z$ form the $0$-diagonal of $G(F(N)^\cdt)^\cdt$. Moreover, $\La^!=kQ^{\rm o}/R^!=\{\gamma^! \mid \gamma\in kQ\}$, where $\gamma^!=\gamma^{\rm o}+R^!$, while $(\La^!)^{\rm o}=kQ/(R^!)^{\rm o}=\{\hat{\gamma\hspace{2pt}} \mid \gamma\in kQ\}$, where $\hat{\gamma\hspace{2pt}}=\gamma+(R^!)^{\rm o}$. Given $a,y\in Q_0$, there exists a linear map $$N_{a,y}: N(y)\to \Hom_k(e_y (\La^!)^{\rm o} e_a, P_a(a)\otimes N(a)): u\mapsto N_{a,y}(u),$$ where $N_{a,y}(u)$ maps $\hat{\gamma\hspace{2pt}}$ to $e_a\otimes\gamma^! u,$ for all $\gamma\in kQ(a, y).$ By Corollary \ref{Cor 1.2}, there exists a $k$-isomorphism
$$\theta_{a,y}: \Hom_k(e_y (\La^!)^{\rm o} e_a, k)\otimes P_a(a)\otimes N(a) \to \Hom_k(e_y (\La^!)^{\rm o} e_a, P_a(a)\otimes N(a)).$$
This yields a linear map $f^a_y=\theta_{a,y}^{-1}\circ N_{a,y}: N(y)\to I^!_a(y)\otimes P_a(a)\otimes N(a).$

\vspace{2pt}

{\sc Statement 3.} {\it If $\{\hat{\gamma\hspace{2pt}}\hspace{-1.5pt}_1,\cdots, \hat{\gamma\hspace{2pt}}\hspace{-2pt}_s\}$ is a basis of $e_y(\La^!)^{\rm o}e_a$ with dual basis $\{\hat{\gamma\hspace{2pt}}\hspace{-1.5pt}_1^\star,\cdots, \hat{\gamma\hspace{2pt}}\hspace{-2pt}_s^\star\}$, where $a, y\in Q_0$, then $f^a_y(u)=\sum_{i=1}^s \hat{\gamma\hspace{2pt}}\hspace{-1.5pt}^\star_i\otimes e_a\otimes \gamma_i^! u,$ for all $u\in N(y)$.}

\vspace{2pt}

Indeed, every $\hat{\gamma\hspace{2pt}}\hspace{-1.5pt}\in e_z(\La^!)^{\rm o}e_a$ is written as $\hat{\gamma\hspace{2pt}}\hspace{-1.5pt}=\sum_{j=1}^s\, \lambda_j \hat{\gamma\hspace{2pt}}\hspace{-1.5pt}_j$, for some $\lambda_j\in k$. Given $u\in N(z)$, by the definition of $\theta_{a, y}$, we obtain
$$\textstyle\theta_{a,y}(\sum_{i=1}^s \hat{\gamma\hspace{2pt}}\hspace{-1.5pt}^\star_i\otimes e_a\otimes \gamma_i^! u)(\hat{\gamma\hspace{2pt}}\hspace{-1.5pt})=\sum_{1\le i, j\le s} (\lambda_j \hspace{.4pt} \hat{\gamma\hspace{2pt}}\hspace{-1.5pt}^\star_i(\hat{\gamma\hspace{2pt}}\hspace{-1.5pt}_j)) \cdot ( e_a\otimes \gamma_i^! u)=e_a\otimes \gamma^!u=N_{a,y}(u)(\hat{\gamma\hspace{2pt}}\hspace{-1.5pt}).$$
Thus, $\theta_{a,z}(\sum_{i=1}^s \hat{\gamma\hspace{2pt}}\hspace{-1.5pt}^\star_i\otimes e_a\otimes \gamma_i^! u)=N_{a,y}(u)$, and hence, $f^a_y(u)=\sum_{i=1}^s \hat{\gamma\hspace{2pt}}\hspace{-1.5pt}^\star_i\otimes e_a\otimes \gamma_i^! u.\vspace{.5pt}$

\vspace{2pt}

{\sc Statement 4.} \vspace{1pt} {\it Given any $a\in Q_0$, there exists a natural $\La^!$-linear morphism $f^a=(f^a_y)_{y\in Q_0}: N\to I^!_a\otimes P_a(a)\otimes N(a)$.}

\vspace{2pt}

Indeed, for any $\alpha: z\to y$ in $Q_1$, it is easy to verify that commutativity of
$$\xymatrixcolsep{35pt}\xymatrixrowsep{30pt}\xymatrix{
N(y) \ar[r]^-{N_{a,y}} \ar[d]_{N(\alpha^{\rm o})} &  \Hom_k(P^{!, \hspace{.4pt}{\rm o}}_a(y), P_a(a) \otimes N(a)) \ar[d]^{\Hom(P_a^{!, \hspace{.4pt}{\rm o}}(\alpha), P_a(a) \otimes N(a))}
&  I_a^!(y) \otimes P_a(a) \otimes N(a) \ar[l]_-{\theta_{a,y}} \ar[d]^{I_a^!(\alpha^{\rm o})\otimes \id \otimes \id}\\
N(z) \ar[r]^-{N_{a,z}} &  \Hom_k(P^{!, \hspace{.4pt}{\rm o}}_a(z), P_a(a) \otimes N(a))
&  I_a^!(z) \otimes P_a(a) \otimes N(a). \ar[l]_-{\theta_{a,z}}
}$$
Thus, $f^a$ is $\La^!$-linear. Given a $\La^!$-linear morphism $g: N\to M$, we have a diagram

$$\xymatrixcolsep{35pt}\xymatrixrowsep{25pt}\xymatrix{
N(y) \ar[r]^-{N_{a,y}} \ar[d]_{g_y} &  \Hom_k(P^{!, \hspace{.4pt}{\rm o}}_a(y), P_a(a) \otimes N(a)) \ar[d]^{\Hom(P_a^{!, \hspace{.4pt}{\rm o}}(y), \id \otimes g_a)}
&  I_a^!(y) \otimes P_a(a) \otimes N(a) \ar[l]_-{\theta_{a,y}} \ar[d]^{\id \otimes \id \otimes g_a }\\
M(y) \ar[r]^-{M_{a,y}} &  \Hom_k(P^{!, \hspace{.4pt}{\rm o}}_a(y), P_a(a) \otimes M(a))
&  I_a^!(y) \otimes P_a(a) \otimes M(a), \ar[l]_-{\theta_{a,y}}
}$$ where the left square is easily verified to be commutative, while the commutativity of the right square follows from the naturality stated in Lemma \ref{Cor 1.2}(1). 

\vspace{2pt}

Given $a\in Q^{-i}$, in view of Statement (3), we obtain a natural $\La^!$-linear morphism $g^a=(g^a_y)_{y\in Q_0}: N\to I^!_a\otimes P_a(a)\otimes N(a)$, where $g^a_y=(-1)^{\frac{(i-1)i}{2}}f^a_y$.

\vspace{2.5pt}

{\sc Statement 5.} {\it Setting $g=(g^a)_{a\in Q_0}$, we obtain a natural $\La^!$-linear monomorphism
$g: N\to 
(G^C \hspace{-2.5pt} \circ \hspace{-1.5pt} F)(N)^0.$}

\vspace{1pt}

Indeed, $g$ is a $\La^!$-linear monomorphism if and only if, for any $y\in Q_0,$ the linear morphism $g_y=(g^a_y): N(y) \to (G^C \hspace{-2.5pt} \circ \hspace{-1.5pt} F)(N)^0=\oplus_{a\,\in Q_0}\, I^!_a(y)\otimes P_a(a)\otimes N(a)$ is injective. Assume now that $g_y(u)=0,$ for some $u\in N(y)$. Then
$g^a_y(u)=0$, for every $a\in Q_0$. In particular, $g^y_y(u)=0$, and hence, $f^y_y(u)=0$. Since $\{e_y\}$ is a basis of $e_y(\La)^{\rm o}e_y$, in view of Statement 3, we see that $e_y^\star\otimes e_y \otimes u=0$, and hence, $u=0$.

\vspace{2pt}

For the rest of the proof, observing that the $(-1)$-diagonal of $G(F(M)^\ydt)^\ydt$ contains only zero objects, we illustrate its $0$-diagonal and $1$-diagonal as follows:
$$\xymatrixrowsep{16pt}\xymatrix{
\oplus_{b\,\in Q^{-i}}\, I^!_b\otimes P_b(b)\otimes N(b) \ar[r]^-{h^{i,-i}} & \oplus_{(a, x)\,\in Q^{-i-1}\times Q^{-i}}\, I^!_x\otimes P_a(x)\otimes N(a) \\
& \oplus_{c\, \in Q^{-i-1}}  I^!_{c}\otimes P_c(c)\otimes N(c) \ar[u]_-{v^{i+1, -i-1}},}$$

\noindent where $h^{i,-i}=(h^{i,-i}(a,x,b))_{(a,x,b)\in Q^{-i}\times Q^{-i-1}\times Q^{-i}}$, \vspace{1pt} with $h^{i,-i}(a,x,b)=0$ for $b\ne x$, and $h^{i,-i}(a,x,x)= {\textstyle\sum}_{\alpha\in Q_1(a, x)} \id_{I^!_x} \otimes P[\bar\alpha]\otimes N(\alpha^{\rm o}),$ and on the other hand, we have $v^{i+1,-i-1}\hspace{-2pt}=\hspace{-2pt}(v^{i+1,-i-1}(a,\hspace{-1pt}x,\hspace{-1pt}c))_{(a,x,c)\in Q^{-i} \times Q^{-i-1} \times Q^{-i-1}}\vspace{1.5pt}$ with
$v^{i+1,-i-1}(a,\hspace{-1pt}x,\hspace{-1pt}c)\hspace{-2pt}=\hspace{-2pt}0$ for $c\ne a$, and $v^{i+1,-i-1}(a,x,a)={\textstyle\sum}_{\alpha\in Q_1(a, x)}(-1)^{i+1} I[\alpha^!]\otimes P_a(\alpha)\otimes \id_{N(a)}$.

\smallskip

{\sc Statement 6}. {\it If $d^{\hspace{.7pt}0}$\hspace{-1pt} is the $0$-degree differential of $(G^C \circ F)(N)^\ydt$, then $d^{\hspace{.6pt}0} \hspace{-1pt} \circ g=0$.}

\vspace{1pt}

Indeed, it amounts to show, for any $p\in \Z$, that the diagram
$$\xymatrixcolsep{25pt}\xymatrixrowsep{22pt}\xymatrix{
\oplus_{x\,\in Q^{-p}}\, I^!_x\otimes P_x(x)\otimes N(x) \ar[rr]^-{\oplus h^{p,-p}(a,x,x)} && \oplus_{(a, x)\,\in Q^{-p-1}\times Q^{-p}}\, I^!_x\otimes P_a(x)\otimes N(a) \\
N\ar[rr]^-{(g^a)_{a\in Q^{-p-1}}} \ar[u]^-{(g^x)_{x\in Q^{-p}}} && \oplus_{a\, \in Q^{-p-1}}  I^!_{a}\otimes P_a(a)\otimes N(a), \ar[u]_-{\oplus v^{p+1,-p-1}(a,x,a)}}$$ is anti-commutative, or equivalently, we have an anti-commutative diagram
$$\xymatrixcolsep{25pt}\xymatrixrowsep{22pt}\xymatrix{
\oplus_{x\,\in Q^{-p}}\, I^!_x(y)\otimes P_x(x)\otimes N(x) \ar[rr]^-{\oplus h^{p,-p}(a,x,x)(y)} && \oplus_{(a, x)\,\in Q^{-p-1}\times Q^{-p}}\, I^!_x(y)\otimes P_a(x)\otimes N(a) \\
N(y)\ar[rr]^-{(g^a_y)_{a\in Q^{-p-1}}} \ar[u]^-{(g^x_y)_{x\in Q^{-p}}} && \oplus_{a\, \in Q^{-p-1}}  I^!_a(y)\otimes P_a(a)\otimes N(a), \ar[u]_-{\oplus v^{p+1,-p-1}(a,x,a)(y)}}$$ for all $y\in Q_0$. Fix $u\in N(y)$ for some $y\in Q_0$. Consider $\alpha\in Q_1(a, x)$ with $(a, x)\,\in Q^{-p-1}\times Q^{-p}$. Choosing a $k$-basis $\{\hat\delta_1, \ldots, \hat\delta_s\}$ of $e_y(\La^!)^{\rm o} e_x$, we deduce from Statement 3 that
$\textstyle (\id\otimes P[\bar\alpha]\otimes N(\alpha^{\rm o})) \left(g^x_y(u)\right)=(-1)^{\frac{(p-1)p}{2}} \sum_{i=1}^s
\hat\delta_i^\star \otimes \bar \alpha \otimes \alpha^! \delta_i^! u.$
On the other hand, for any $k$-basis $\{\hat{\gamma\hspace{2pt}}\hspace{-1.5pt}_1, \ldots, \hat{\gamma\hspace{2pt}}\hspace{-1.5pt}_t\}$ of $e_y(\La^!)^{\rm o} e_a$, we obtain
$$\textstyle (I[\alpha^!]\otimes P_a(\alpha)\otimes \id\,) \left(g^a_y(u)\right)=(-1)^{\frac{p(p+1)}{2}}\sum_{i=1}^tI[\alpha^!](\hat{\gamma\hspace{2pt}}\hspace{-1.5pt}_i^\star) \otimes \bar \alpha \otimes \gamma_i^! u.$$

Let $\theta: I^!_x(y)\otimes P_a(x)\otimes N(a)\to \Hom_k(e_y(\La^!)^{\rm o} e_x, P_a(x)\otimes N(a))$ be a $k$-linear isomorphism as stated in Corollary \ref{Cor 1.2}. Given any $1\le j\le s$, it is easy to see that
$$\textstyle \theta[(\id\otimes P[\bar\alpha]\otimes N(\alpha^{\rm o})) \left(g^x_y(u)\right)](\hat\delta_j)
=(-1)^{\frac{(p-1)p}{2}}(\bar \alpha \otimes \alpha^! \delta_j^! u),$$
and $\textstyle\theta[(I[\alpha^!]\otimes P_a(\alpha)\otimes \id\,) \left(g^a_y(u)\right)](\hat\delta_j) =
(-1)^{\frac{p(p+1)}{2}} \sum_{i=1}^t \hat{\gamma\hspace{2pt}}\hspace{-1.5pt}_i^\star(\hat\delta_j \hat\alpha) \left(\bar \alpha \otimes \gamma_i^! u\right).$

Fix some $1\le j\le s$. If $\hat\delta_j \, \hat\alpha=0$, then $\alpha^! \delta_j^!=0$, and hence,
$$\textstyle\theta[(I[\alpha^!]\otimes P_a(\alpha)\otimes \id\,)\left(g^a_y(u)\right)](\hat\delta_j)=0= (-1)^p\theta [(\id\otimes P[\bar\alpha]\otimes N(\alpha^{\rm o})) \left(g^x_y(u)\right)](\hat\delta_j).$$
If $\hat\delta_j \, \hat\alpha\ne 0$, then it extends to a $k$-basis $\{\hat{\gamma\hspace{2pt}}\hspace{-1.5pt}_1, \ldots, \hat{\gamma\hspace{2pt}}\hspace{-1.5pt}_t\}$ with $\hat{\gamma\hspace{2pt}}\hspace{-1.5pt}_1=\hat\delta_j \, \hat\alpha$ of $e_y(\La^!)^{\rm o} e_a$. Under this assumption, we obtain
$$\begin{array}{rcl}
\theta[(I[\alpha^!]\otimes P_a(\alpha)\otimes \id\,) \left(g^a_y(u)\right)](\hat\delta_j)
&=& (-1)^{\frac{p(p+1)}{2}} \sum_{i=1}^t \hat{\gamma\hspace{2pt}}\hspace{-1.5pt}_i^\star(\hat{\gamma\hspace{2pt}}\hspace{-1.5pt}_1) (\bar \alpha \otimes \gamma_i^! u ) \vspace{1pt} \\
&=& (-1)^{\frac{p(p+1)}{2}}  (\bar \alpha \otimes \gamma_1^! u ) \\
&=& (-1)^{\frac{p(p+1)}{2}} (\bar \alpha \otimes \hat\delta_j  \hat\alpha \, u)  \\
&=& (-1)^p\theta [(\id\otimes P[\bar\alpha]\otimes N(\alpha^{\rm o})) \left(g^x_y(u)\right)](\hat\delta_j).
\end{array}$$
Thus, $\theta[(I[\alpha^!]\otimes P_a(\alpha)\otimes \id\,) \left(g^a_y(u)\right)]= (-1)^p\theta [(\id\otimes P[\bar\alpha]\otimes N(\alpha^{\rm o})) \left(g^x_y(u)\right)],$ and hence, $(I[\alpha^!]\otimes P_a(\alpha)\otimes \id\,) \left(g^a_y(u)\right)= (-1)^p (\id\otimes P[\bar\alpha]\otimes N(\alpha^{\rm o})) \left(g^x_y(u)\right).\vspace{1pt}$ This yields that $(h^{p,-p}(a,x,x)(y)\circ g^x_y)(u)+(v^{p+1,-p-1}(a,x,a)(y)\circ g^a_y)(u)=0$, \vspace{1pt} and hence, $h^{p,-p}(a,x,x)(y)\circ g^x_y+v^{p+1,-p-1}(a,x,a)(y)\circ g^a_y=0$. This in turn implies the required anti-commutativity. 

\smallskip

Let $\omega=(\omega^i)_{i\in \Z}\in {\rm Ker}(d^{\hspace{.6pt}0})$, with $\omega^i\in G(F(N)^i)^{-i}=\oplus_{a\in Q^{-i}}\, I^!_a\otimes P_a(a)\otimes N(a)$. Recall that $G(F(N)^i)^{-i}=0$, for $i> r$. We shall define $n_\omega (\hspace{.7pt}\le r)$ as follows. If $\omega=0,$ then $n_\omega=r$; and otherwise, $n_\omega$ is minimal for which $w^{n_\omega}\neq 0$. If $n_\omega=r$, then clearly $\omega\in {\rm Im}(g)$. Assume that $n_\omega<r$. Since $\omega\in {\rm Ker}(d^0)$, we have $v^{n_\omega,-n_\omega}(\omega^{n_\omega})=-h^{n_\omega-1,1-n_\omega}(\omega^{n_\omega-1})=0$. By Statement 1, the $n_\omega$-th column of $G(F(N)^\cdt)^\cdt$ is, 
up to a twist, the shift by $n_\omega$ of the minimal injective co-resolution of $\oplus_{a\in Q^{-n_\omega}}\, S^!_a \otimes P_a(a) \otimes N(a)$. Thus, $w^{n_\omega}\in S(\oplus_{a\in Q^{-n{\hspace{-.5pt}_\omega}}}\, I^!_a\otimes P_a(a)\otimes N(a))$, and by Lemma \ref{Special-func}, $\omega^{n_\omega}=\sum_{a\in Q^{-n{\hspace{-.5pt}_\omega}}} \hat{e}^\star_a\otimes e_a\otimes u_a$, where $u_a\in N(a)$.

\vspace{1pt}

By Statement 3, we obtain $g(\sum_{a\in Q^{-n_\omega}} u_a)=\sum_{a\in Q^{-n_\omega}} \hat{e}^\star_a\otimes e_a\otimes u_a=\omega^{n_\omega}$, and by Statement 6, $\nu=\omega-g(\sum_{a\in Q^{-n_\omega}} u_a)\in {\rm Ker}(d^{\hspace{.7pt}0})$. Writing $\nu=(\nu^i)_{i\in \Z}$ with $\nu^i\in G(F(N)^i)^{-i}$, we see that $\nu^{n_\omega}=\omega^{n_\omega}-g(\sum_{a\in Q^{-n_\omega}} u_a)=0$, and $\nu^i=\omega^i=0$ for all $i < n_\omega$. Therefore, $n_\nu<n_\omega$. Assuming inductively that
$\nu\in {\rm Im}(g)$, we obtain $\omega\in {\rm Im}(g)$. This shows that ${\rm Ker}(d^{\hspace{.7pt}0})={\rm Im}(g)$. Setting $\theta_N^0=g$ and $\theta_N^i=0$ for all $i\ne 0$, we obtain a quasi-isomorphism $\theta^\pdt_N: N\to (G^C\circ F)(N)^\cdt$, which is natural in $N$ by Statement 4. The proof of the proposition is completed.

\medskip

The following statement is one of the main results of this paper. It extends under the gradable setting Theorem 2.12.1 in \cite{BGS} and Theorem 30 in \cite{MOS}, both deal with only one pair of subcategories of the graded derived categories, consisting of some bounded-above complexes and of some bounded-below complexes respectively.

\medskip

\begin{Theo}\label{Main}

Let $\La=kQ/R$ be a Koszul algebra, where $Q$ is a locally finite gradable quiver. Given any $p, q\in \mathbb{R}$ with $p\ge 1$ and $q\ge 0$, \vspace{1.5pt} the derived Koszul functor $F^D: D^{\;\downarrow}_{p,q}({\rm Mod}\hspace{.4pt}\La^!)\to D^{\,\uparrow}_{q+1, p-1}({\rm Mod}\hspace{.4pt}\La)\vspace{1pt}$ is a triangle-equivalence with quasi-inverse
$G^D: D^{\,\uparrow}_{q+1, p-1}({\rm Mod}\hspace{.4pt}\La)\to D^{\;\downarrow}_{p,q}({\rm Mod}\hspace{.4pt}\La^!)$.

\end{Theo}

\vspace{2pt}

\noindent{\it Proof.} Let $p, q\in \mathbb{R}$ with $p\ge 1$ and $q\ge 0$. In view of Theorem \ref{F-diag}, we obtain a commutative diagram
$$\xymatrixrowsep{16pt}\xymatrix{
C^{\,\downarrow}_{p,q}({\rm Mod}\hspace{.4pt}\La^!) \ar[d]_{F^C} \ar[r]^{P_{\hspace{-1pt}\it\Lambda^!}} & K^{\,\downarrow}_{p,q}({\rm Mod}\hspace{.4pt}\La^!)\ar[d]^{F^K} \ar[r]^{L_{\hspace{-.8pt}\it\Lambda^!}} & D^{\,\downarrow}_{p,q}({\rm Mod}\hspace{.4pt}\La^!) \ar[d]^{F^D} \\
C^{\,\uparrow}_{q+1, p-1}({\rm Mod}\hspace{.4pt}\La)\ar[r]^{P_{\hspace{-.6pt}\it\Lambda}} \ar[d]_{G^C} & K^{\,\uparrow}_{q+1, p-1}({\rm Mod}\hspace{.4pt}\La) \ar[r]^{L_{\hspace{-.6pt}\it\Lambda}} \ar[d]^{G^K} & D^{\,\uparrow}_{q+1, p-1}({\rm Mod}\hspace{.4pt}\La)\ar[d]^{G^D}\\
C^{\,\downarrow}_{p,q}({\rm Mod}\hspace{.4pt}\La^!) \ar[r]^{P_{\hspace{-1pt}\it\Lambda^!}} & K^{\,\downarrow}_{p,q}({\rm Mod}\hspace{.4pt}\La^!) \ar[r]^{L_{\hspace{-.8pt}\it\Lambda^!}} & D^{\,\downarrow}_{p,q}({\rm Mod}\hspace{.4pt}\La^!).}$$

First, we shall show that $G^D \circ F^D \cong \id_{D^{\,\downarrow}_{p,q}({\rm Mod}\hspace{.4pt}\it\Lambda)}.$ Observe that $C^{\,\downarrow}_{p,q}({\rm Mod}\hspace{.4pt}\La^!)$ is a full subcategory of
$C({\rm Mod}^-\hspace{-2.5pt}\La^!)$. Consider $G^C\circ F: {\rm Mod}^-\hspace{-3pt}\La^!\to C({\rm Mod}\La^!)$ and the canonical embedding functor $\kappa: {\rm Mod}^-\hspace{-3pt}\La^! \rightarrow C({\rm Mod}\La^!)$. By Proposition \ref{natural-tran-2}, we obtain a functorial morphism $\theta=(\theta_{\hspace{-1pt}_N}^\pdt)_{N\in {\rm Mod}^-\hspace{-2pt}\it\Lambda^!}:  \kappa \to G^C \hspace{-2pt} \circ F ,\vspace{1pt}$ where $\theta_{\hspace{-1pt}_N}^\pdt: \kappa(M)^\cdt \to (G^C \hspace{-2pt} \circ F)(N)^\cdt$ is a quasi-isomorphism. By Lemmas \ref{F-composition} and \ref{main-lemma1}, $\theta$ induces a functorial morphism 
$\theta^C: \kappa^C\to (G^C \hspace{-2pt} \circ F)^C=G^C \circ F^C.$

Given $N^{\ydt}\in C^{\,\downarrow}_{p,q}({\rm Mod}\hspace{.4pt}\La^!)$, we claim that $\theta^C_{\hspace{-.5pt}_{N^\ydt}}: N^\ydt\to (G^C\circ F)^C(N^\cdt)$ is a quasi-isomorphism. Since $N^i\in {\rm Mod}^-\hspace{-3pt}\La^!$, by Proposition \ref{natural-tran-2}, $\theta_{\hspace{-1pt}{N^i}}^\pdt: \kappa(N^i)^\cdt \to (G^C \hspace{-2pt} \circ F)(N^i)^\cdt$ is a quasi-isomorphism, for every $i\in \Z$.
It is evident that the double complex $\kappa(N^\cdt)^\cdt$ is diagonally bounded-above. Given any integer $n$, the $n$-diagonal of the double complex $(G^C\circ F)(N^\cdt)^\cdt$ consists of the modules
$$(G^C\circ F)(N^i)^{n-i} =\oplus_{j\in \Z; \hspace{.4pt} x\in Q^{-j}\hspace{-1pt}; \hspace{.4pt}y\in Q^{n-i-j} }\, I_y^!\otimes P_x(y)\otimes N^i(x); \, i\in \Z.$$

If $i>n$, then $P_x(y)=0$, for all $x\in Q^{-j}$ and $y\in Q^{n-i-j}$ with $j\in \Z$, and hence, $(G^C\circ F)(N^i)^{n-i}=0$. That is, $(G^C\circ F)(N^\cdt)^\cdt$ is diagonally bounded-above. By Lemma \ref{main-lemma2}, $\theta^{C}_{_{N^\ydt}}$ is indeed a quasi-isomorphism in $C^{\,\downarrow}_{p,q}({\rm Mod}\hspace{.4pt}\La^!)$. As a consequence, $\theta^D_{_{N^\ydt}}=L_{\it\Lambda^!}(P_{\it\Lambda^!}(\theta^C_{_{N^\ydt}})): N^\ydt\to (G^D\circ F^D)(N^\cdt)$
is a natural isomorphism in $D^{\,\downarrow}_{p,q}({\rm Mod}\hspace{.4pt}\La^!).\vspace{1pt}$ This yields a functorial isomorphism $\theta^D:
G^D \hspace{-1.5pt} \circ \hspace{-.5pt} F^D\to \id_{D^{\,\downarrow}_{p,q}({\rm Mod}\hspace{.4pt}\it\Lambda^!)}$.

Next, we shall show that $F^D\circ G^D \cong \id_{D^{\,\uparrow}_{q+1, p-1}({\rm Mod}\hspace{.4pt}\it\Lambda)}.$ For this purpose, consider the canonical embedding $\kappa: {\rm Mod}\hspace{.5pt}\La \rightarrow C({\rm Mod}\hspace{.5pt}\La)$ and $F^C\circ G: {\rm Mod}\hspace{.5pt}\La \rightarrow C({\rm Mod}\hspace{.5pt}\La)$. By Lemma \ref{natural-tran}, we obtain a functorial morphism $\eta=(\eta_{_M}^\ydt)_{M\in {\rm Mod}\it\Lambda}: F^C \hspace{-2pt} \circ G \to \kappa,\vspace{1pt}$ where $\eta_{_M}^\ydt: (F^C \hspace{-2pt} \circ G)(M)^\cdt\to \kappa(M)^\cdt$ is a quasi-isomorphism. By Lemmas  \ref{F-composition} and \ref{main-lemma1}, $\eta$ induces a functorial morphism $\eta^C: (F^C \hspace{-2pt} \circ G)^C\to \kappa^C$, that is a functorial morphism $\eta^C: F^C \circ G^C\to \id_{\,C({\rm Mod}\it\Lambda)}$. As argued above, it suffices to show, for any $M^{\ydt}\in C^{\,\uparrow}_{q+1, p-1}({\rm Mod}\hspace{.4pt}\La)$, that $\eta^C_{_{\hspace{-1pt}M^\pdt}}: (F^C \hspace{-2pt} \circ G)^C(M^\cdt)\to \kappa^C(M^\cdt)$ is a quasi-isomorphism, or equivalently, it has an acyclic cone.

By definition, $\eta^C_{_{\hspace{-1pt}M^\pdt}}=\eta^{\,\cdt}_{_{\hspace{-1pt}M^\pdt}}=\T(\eta^{\cdt {} \cdt}_{_{\hspace{-1pt}M^\pdt}}),\vspace{1pt}$ where $\eta^{\cdt {} \cdt}_{_{\hspace{-1pt}M^\pdt}}: (F^C\circ G)(M^\cdt)^\cdt\to \kappa(M^\cdt)^\cdt$ is defined by $\eta^{i, j}_{_{\hspace{-1pt}M^\pdt}}=\eta_{_{\hspace{-1pt}M^i}}^j: (F^C\circ G)(M^i)^j\to \kappa(M^i)^j$. By Lemma \ref{Cone}, $C_{\eta^{\,\cdt}_{_{\hspace{-1pt}M^\pdt}}}=\T(V_{\eta^{\cdt {} \cdt}_{_{\hspace{-1pt}M^\pdt}}}).\vspace{1pt}$ To show that $C_{\eta^{\,\cdt}_{_{\hspace{-1.5pt}M^\pdt}}}\vspace{1pt}$ is acyclic, it amounts to show that $C_{\eta^{\,\cdt}_{_{\hspace{-1pt}M^\pdt}}}(z)=
\T(V_{\eta^{\cdt {} \cdt}_{_{\hspace{-1pt}M^\pdt}}}(z))$ is acyclic for any $z\in Q_0$, where $V_{\eta^{\cdt {} \cdt}_{_{\hspace{-1pt}M^\pdt}}}(z)$ is the vertical cone of the double complex morphism $\eta^{\cdt {} \cdt}_{_{\hspace{-1pt}M^\pdt}}(z): (F^C\circ G)(M^\cdt)^\cdt(z)\to \kappa(M^\cdt)^\cdt(z)$.
Assume that $z\in Q^s$ with $s\in \Z$. Since the $\eta_{_{\hspace{-.6pt}M^i}}^\cdt: (F^C\circ G)(M^i)^\cdt\to \kappa(M^i)^\cdt$ are quasi-isomorphisms, by Lemma \ref{V-cone}(1), $V_{\eta^{\cdt {} \cdt}_{_{\hspace{-1pt}M^\pdt}}}$ has acyclic columns, and so does $V_{\eta^{\cdt {} \cdt}_{_{\hspace{-1pt}M^\pdt}}}(z)$. Given any $n\in \Z$, the $n$-diagonal of the double complex $(F^C\circ G)(M^\ydt)^\cdt(z)$ consists of the modules
$$\begin{array}{rcl}
(F^C\circ G)(M^i)^{n-i}(z) & = &\oplus_{j\in \mathbb{Z}; \, x\in Q^j; \, y\in Q^{i+j-n}}\, P_y(z)\otimes I_x^!(y)\otimes M^i(x) \vspace{1pt} \\
&=&
\oplus_{j\le n+s-i; \, x\in Q^j; \, y\in Q^{i+j-n}}\, P_y(z)\otimes I_x^!(y)\otimes M^i(x),
\end{array}$$
with $i\in \Z.$ Since $M^{\ydt}\in C^{\,\uparrow}_{q+1, p-1}({\rm Mod}\hspace{.4pt}\La)$, there exists some $t$ for which $M^i(x)=0$ for $x\in Q^j$ with $i-(p-1)j>t$. Let $x\in Q^j$ with \vspace{1pt} $j\le n+s-i$.  If $i> \frac{(p-1)(n+s)+t}{p}$, then
$i-(p-1)j\ge i-(p-1)(n+s-i)=pi-(p-1)(n+s)>t,$
and hence, $M^i(x)=0$. This shows that $(F^C\circ G)(M^\ydt)^\cdt(z)$ is diagonally bounded-above. Since $\kappa(M^\cdt)^\cdt(z)$ is clearly diagonally bounded-above, we see that $V_{\eta^{\cdt {} \cdt}_{_{\hspace{-1pt}M^\pdt}}}(z)$ is diagonally bounded-above. By Corollary \ref{AAL},
$\T(V_{\eta^{\cdt {} \cdt}_{_{\hspace{-1pt}M^\pdt}}}(z))$ is acyclic, that is, $C_{\eta^{\,\cdt}_{_{\hspace{-1pt}M^\pdt}}}(z)$ is acyclic.
The proof of the theorem is completed.

\medskip

\noindent{\sc Remark.} In case $\La$ is Koszul, by Lemma \ref{inj-im} and Theorem \ref{Main}, the Koszul derived functors send respective indecomposable injective modules and indecomposable projective modules.

\medskip

The following statement is a generalization of the result stated in \cite[(2.12.6)]{BGS}, which deals with the derived category of finite dimensional modules in $\La$ is finite dimensional.

\medskip

\begin{Theo}\label{Main-2}

Let $\La=kQ/R$ be a Koszul algebra, where $Q$ is a locally finite gradable quiver. If $\La$ is right $($respectively, left$\hspace{.5pt})$ locally bounded and $\La^!$ \vspace{1pt} is left $($respectively, right$\hspace{.5pt})$ locally bounded, then we obtain $D^b({\rm Mod}^{\hspace{.5pt}b\hspace{-2.5pt}}\La^!)\cong D^b({\rm Mod}^{\hspace{.5pt}b\hspace{-2pt}}\La)$ \vspace{1pt} and $D^b({\rm mod}^{\hspace{.5pt}b\hspace{-2.5pt}}\La^!)\cong D^b({\rm mod}^{\hspace{.5pt}b\hspace{-2.5pt}}\La).$

\end{Theo}

\noindent{\it Proof.} First, assume that $\La$ is right locally bounded and $\La^!$ is left locally bounded. Then, $P_x\in {\rm mod}^{\hspace{.5pt}b\hspace{-2pt}}\La$ and $I_x^!\in {\rm mod}^{\hspace{.5pt}b\hspace{-3pt}}\La^!$, for all $x\in Q_0$. Therefore, we obtain Koszul functors $F: {\rm Mod}^{\hspace{.5pt}b\hspace{-2pt}}\La^!\to C^b({\rm Mod}^{\hspace{.5pt}b\hspace{-2pt}}\La)$ and $F: {\rm mod}^{\hspace{.5pt}b\hspace{-2pt}}\La^!\to C^b({\rm mod}^{\hspace{.5pt}b\hspace{-2pt}}\La)$, and Koszul inverse functors
$G: {\rm Mod}^{\hspace{.5pt}b\hspace{-2pt}}\La \to C^b({\rm Mod}^{\hspace{.5pt}b\hspace{-2pt}}\La^!)$ and $G: {\rm mod}^{\hspace{.5pt}b\hspace{-2pt}}\La\to C^b({\rm mod}^{\hspace{.5pt}b\hspace{-2pt}} \La^!)$. Using the same argument as in the proof of Theorem \ref{Main}, \vspace{1pt} we obtain triangle-equivalences $F^D: D^b({\rm Mod}^{\hspace{.5pt}b\hspace{-2pt}}\La^!)\to D^b({\rm Mod}^{\hspace{.5pt}b\hspace{-2pt}}\La)$ and $F^D: D^b({\rm mod}^{\hspace{.5pt}b\hspace{-2pt}}\La^!)\to D^b({\rm mod}^{\hspace{.5pt}b\hspace{-2pt}}\La)$. Suppose that $\La$ is left locally bounded and $\La^!$ \vspace{1pt} is right locally bounded. Since $\La^!$ is Koszul with $(\La^!)^!=\La$, we have triangle-equivalences
$D^b({\rm Mod}^{\hspace{.5pt}b\hspace{-2.5pt}}\La)\cong D^b({\rm Mod}^{\hspace{.5pt}b\hspace{-2pt}}\La^!)$ and $D^b({\rm mod}^{\hspace{.5pt}b\hspace{-2.5pt}}\La)\cong D^b({\rm mod}^{\hspace{.5pt}b\hspace{-2.5pt}}\La)^!.$ The proof of the theorem is completed.

\medskip

An infinite path in a quiver is called {\it right infinite} if it has no ending point and {\it left infinite} if it has no starting point. 

\medskip

\begin{Cor}\label{Main-3}

Let $\La=kQ/R$ be a Koszul algebra, where $Q$ is a locally finite gradable quiver. If $Q$ has no right infinite path or no left infinite path, then $D^b({\rm Mod}^{\hspace{.5pt}b\hspace{-2.5pt}}\La^!)\cong D^b({\rm Mod}^{\hspace{.5pt}b\hspace{-2pt}}\La)$ and $D^b({\rm mod}^{\hspace{.5pt}b\hspace{-2.5pt}}\La^!)\cong D^b({\rm mod}^{\hspace{.5pt}b\hspace{-2.5pt}}\La).$

\end{Cor}

\noindent{\it Proof.} If $Q$ has no right infinite path, then $Q^{\rm o}$ has no left infinite path, and in particular, $\La$ is right locally bounded and $\La^!$ is left locally bounded. If $Q$ has no left infinite path, then $\La$ is left locally bounded and $\La^!$ is right locally bounded. The proof of the corollary is completed.

\bigskip

{\sc Acknowledgement.} The third named author is supported in part by the Natural Sciences and Engineering Research Council of Canada.

\end{document}